\documentclass[a4paper, 10pt, oneside, reqno]{amsart}

\usepackage[utf8]{inputenc}
\usepackage{amsfonts,amssymb,amsthm,amsmath}
\usepackage{mathrsfs}
\usepackage{geometry}
\usepackage{graphicx,graphics}
\usepackage{color}
\usepackage{tikz}
\usepackage{tikz-cd}
\usepackage{comment}
\usepackage{subfiles}
\usepackage{stmaryrd}
\usepackage{hyperref}
\usepackage{todonotes}

\hypersetup{colorlinks = true, linkbordercolor = {white}, linkcolor = {blue}, citecolor = {blue}}

\theoremstyle{plain}
\newtheorem{thm}{Theorem}[subsection]
\newtheorem{lem}[thm]{Lemma}

\newtheorem{corol}[thm]{Corollary}
\theoremstyle{definition}
\newtheorem{defi}[thm]{Definition}

\theoremstyle{remark}
\newtheorem{rque}[thm]{Remark}

\def \F {\mathbb{F}}

\def \Fq {\F_q}
\def \Fl {\F_{\ell}}

\def \Fqb {\overline{\F}_q}
\def \Flb {\overline{\F}_{\ell}}

\def \Z {\mathbb{Z}}
\def \Zl {\Z_{\ell}}
\def \Zlb {\overline{\Z}_{\ell}}

\def \Q {\mathbb{Q}}
\def \Ql {\Q_{\ell}}
\def \Qlb {\overline{\Q}_{\ell}}
\def \RR {\mathrm{R}}

\def \Ind {\mathrm{Ind}}

\def \Gm {\mathbb{G}_m}

\def \Spec {\mathrm{Spec}}

\def \A {\mathbb{A}}

\def \Gbf {\mathbf{G}}
\def \Bbf {\mathbf{B}}
\def \Ubf {\mathbf{U}}
\def \Tbf {\mathbf{T}}
\def \Wbf {\mathbf{W}}
\def \Nbf {\mathbf{N}}
\def \BwB {\Bbf w\Bbf}
\def \Pbf {\mathbf{P}}
\def \Vbf {\mathbf{V}}
\def \Lbf {\mathbf{L}}
\def \Zbf {\mathbf{Z}}

\def \qfrak {\mathfrak{q}}
\def \rfrak {\mathfrak{r}}
\def \chfrak {\mathfrak{ch}}
\def \hcfrak {\mathfrak{hc}}

\def \Hom {\mathrm{Hom}}
\def \D {\mathbb{D}}
\def \id {\mathrm{id}}

\def \Tr {\mathrm{Tr}}

\def \1 {\mathbb{1}}

\def \End {\mathrm{End}}
\def \Mod {\mathrm{Mod}}

\def \ind {\mathrm{ind}}

\def \Cat {\mathrm{Cat}}

\def \Frob {\mathrm{Frob}}

\def \Fun {\mathrm{Fun}}

\def \Frob {\mathrm{F}}
\def \Ad {\mathrm{Ad}}
\def \Av {\mathrm{Av}}

\def \For {\mathrm{For}}

\def \op {\mathrm{op}}

\def \mon {\mathrm{mon}}

\def \Spr {\mathrm{Spr}}
\def \triv {\mathrm{triv}}
\def \pr {\mathrm{pr}}
\def \RGamma {\mathrm{R\Gamma}}

\def \DD {\mathrm{D}}
\def \pt {\mathrm{pt}}

\def \Pro {\mathrm{Pro}}
\def \unip {\mathrm{unip}}
\def \Perv {\mathrm{Perv}}

\def \CH {\mathrm{CH}}

\def \pt {\mathrm{pt}}

\def \CH {\mathrm{CH}}
\def \can {\mathrm{can}}

\def \Ext {\mathrm{Ext}}

\def \Rep {\mathrm{Rep}}
\def \RHom {\mathrm{RHom}}

\def \ho {\mathrm{ho}}

\def \DDc {\DD_{\cons}}
\def \cons {\mathrm{c}}
\def \DDic {\DD_{\mathrm{ic}}}

\def \lis {\mathrm{lis}}
\def \wild {\mathrm{wild}}
\def \tame {\mathrm{tame}}
\def \Perf {\mathrm{Perf}}

\def \inv {\mathrm{inv}}

\def \HH {\mathbb{H}}

\def \rev {\mathrm{rev}}
\def \Irr {\mathrm{Irr}}
\def \coh {\mathrm{coh}}

\def \proet {\mathrm{proet}}
\def \Corr {\mathrm{Corr}}

\def \Stk {\mathrm{Stk}}

\def \lef {\mathrm{left}}
\def \rig {\mathrm{right}}
\def \HHo {\HH^{\omega}}

\def \Rcal {\mathcal{R}}
\def \Ccal {\mathcal{C}}

\def \O {\mathcal{O}}

\def \Irr {\mathrm{Irr}}

\def \res {\mathrm{res}}

\def \Rcal {\mathcal{R}}
\def \tame {\mathrm{tame}}
\def \wild {\mathrm{wild}}

\def \Lcal {\mathcal{L}}
\def \Zcal {\mathcal{Z}}
\def \CS {\mathrm{CS}}

\def \Pre {\mathrm{Pre}}

\def \Ch {\mathrm{Ch}}

\title{Free monodromic Hecke categories and their categorical traces}
\author{Arnaud Eteve}

\begin{document}

\maketitle

\textbf{Abstract.} The goal of this paper is to give a new construction of the free monodromic categories defined by Yun. We then use this formalism to give simpler constructions of the free monodromic Hecke categories and then compute the trace of Frobenius and of the identity on them. As a first application of the formalism, we produce some new proofs of the key theorems in Deligne--Lusztig theory.  

\tableofcontents

\section{Introduction}

Let $p > 0$ be a prime number, $\Fq$ a finite field with $q = p^r$ elements and $\Fqb$ an algebraic closure of $\Fq$. Let $\Gbf$ be a reductive group over $\Fqb$ together with a Frobenius endomorphism $\Frob : \Gbf \to \Gbf$. The first goal of this paper is to rework the formalism of monodromic sheaves of Verdier \cite{Verdier} and free monodromic sheaves of Yun \cite{BezYun}. The second goal of this paper is to use this formalism to produce simpler definitions of the free monodromic Hecke categories and to study their categorical traces.
 As explained in \cite{BenZviNadlerComplexGroup}, in the setting of $D$-modules, these traces are closely related to the theory character sheaves of Lusztig \cite{LusztigCS1}. In the étale context, taking the trace of Frobenius, instead of the identity, yields representations of the finite group $\Gbf^{\Frob}$, furthermore this formalism interacts particularly well with Deligne--Lusztig theory and offers a way to produce simple and conceptual proofs of key results in the theory. Taking the categorical trace (or categorical center) of Frobenius on categories coming from geometric representation theory is a particularly fruitful idea for `arithmetic' problems, that is, constructions where a Frobenius is present. We refer to \cite{HowToInventShtukas} for a general discussion. In the context of representations of finite groups of Lie type, Lusztig \cite{LusztigCatCenter1} \cite{LusztigCatCenter2} has shown such categorical statements. The key difference with our work is that the results of \emph{loc. cit.} require fine knowledge of unipotent representations, character sheaves, their classification and the theory of cells in Coxeter groups while ours requires very little knowledge. In a companion paper \cite{EteveJordanDecomp}, we construct a categorical form of Jordan decomposition, generalizing the construction of \cite{LusztigYun}, for representations of $\Gbf^{\Frob}$ which remains valid for modular coefficients. In the rest of this introduction we give a more detailed overview of the content of this paper.

\subsection{The free monodromic formalism}

Let $X$ be a scheme of finite type over an algebraically closed field $k$ and let $\Tbf$ be a torus over $k$ acting on $X$. Verdier, in \cite{Verdier}, had defined a full subcategory of $\DD^b_c(X, \Lambda)$,  the category of $\Lambda \in \{\Zl, \Fl, \Ql\}$-linear étale sheaves on $X$, called monodromic sheaves. It is composed of all sheaves that are locally constant and tamely ramified along $\Tbf$-orbits. For some applications to representation theory, Yun in \cite[Appendix A]{BezYun} has defined a certain pro-completion of this category of monodromic sheaves under some restriction on the action of $\Tbf$. In the first section of this paper, we propose a different construction of this category which has several advantages, most notably, it remains valid without any hypothesis on the action of $\Tbf$. In \cite{EteveEberhardt}, the Betti version of this construction is worked out.

Let us now highlight our construction. For simplicity we will work over $\Lambda = \Flb$ in this introduction. The case $\Lambda = \Zlb$ is treated in the same way and to treat the case of $\Qlb$-sheaves, we use some condensed mathematics. 

We denote by $\pi_1^{\ell}(\Tbf)$ the maximal pro-$\ell$ quotient of the étale fundamental group of $\Tbf$ at the point $1$. It is known that this group is isomorphic to $X_*(\Tbf) \otimes \Zl(1)$, where  $X_*$ denotes the lattice of cocharacters of $\Tbf$.

\begin{defi}[Section \ref{sectionTheRings}]
We define 
$$\RR_{\Tbf} = \Fl\llbracket \pi_1^{\ell}(\Tbf) \rrbracket \otimes_{\Fl} \Flb.$$
\end{defi}
This ring was first considered by \cite{GabberLoeser}. There is a natural morphism 
\begin{equation*}
\pi_1(\Tbf) \to \RR_{\Tbf}^{\times}.
\end{equation*}
Corresponding to this morphism there is a rank one locally constant $\RR_{\Tbf}$-sheaf on $\Tbf$ which we call the free monodromic unipotent local system on $\Tbf$ and denote by $L_{\Tbf}$. We show in Lemma \ref{lemMultiplicative} that this sheaf is multiplicative in the sense of Appendix \ref{subsectionTwistedEquivSheaves}. Using the formalism of \cite{GaiWhittCats} we can form the corresponding category of twisted equivariant sheaves on $X$. 

\begin{defi}
We define the category of free monodromic unipotent sheaves on $X$ to be 
$$\DDic(X, \RR_{\Tbf})_{\unip}$$
the category of sheaves that are $({\Tbf}, L_{\Tbf})$-equivariant on $X$.
\end{defi}

We define $\Ch_{\Lambda}(\Tbf)$ the set of finite order characters of $\pi_1(\Tbf)$ that are tame and of order prime to $p$ (if $\Lambda = \Qlb$) and $p\ell$ (if $\Lambda = \Zlb, \Flb$). Corresponding to a character $\chi$, there is a rank one Kummer sheaf denoted by $\Lcal_{\chi}$. The sheaf $L_{\Tbf} \otimes_{\Lambda} \Lcal_{\chi}$ is again a multiplicative sheaf and we denote by 
\begin{equation*}
\DDic(X, \RR_{\Tbf})_{\chi}
\end{equation*}
the category of $({\Tbf}, L_{\Tbf} \otimes \Lcal_{\chi})$-equivariant sheaves on $X$. We call the category 
\begin{equation*}
\bigoplus_{\chi \in \Ch_{\Lambda}({\Tbf})}\DDic(X, \RR_{\Tbf})_{\chi},
\end{equation*}
the category of free monodromic sheaves on $X$. 

To motivate our construction, we show the following theorem which compares to the original construction of Yun. 
\begin{thm}[Theorem \ref{ThmEquivalenceBR}]
There is a natural equivalence of categories 
\begin{equation*}
\ho(\DDc(X, \RR_{\Tbf})_{\unip}) = \DD(X\! \fatslash {\Tbf}, \Lambda),
\end{equation*}
where the notation on the RHS is from \cite{BezrukRiche} for the free monodromic unipotent category.
\end{thm} 

The proof holds for all coefficients and for twisted sheaves if we replace the completed category by the one constructed in the thesis of Gouttard \cite{Gouttard}.

\subsection{Hecke categories}

We study the free monodromic Hecke categories in Section \ref{sectionFreeMonoHecke}, these categories were firstly defined in \cite{BezYun} and also studied in \cite{BezrukRiche} and \cite{Gouttard}, where the last two papers focus on the case of modular coefficients. One of the goals of this construction is to generalize their construction to the integral setting (i.e. over $\Zlb$) and to simplify the construction of the main structures on this category. We fix a (possibly disconnected) reductive group $\Gbf$ over $k$ with a Borel pair $\Bbf = \Tbf\Ubf$ and we consider the stack $\Ubf \backslash \Gbf/\Ubf$. 

On this stack there are three tori that acts 
\begin{enumerate}
\item the torus $\Tbf$ acting by left translations, 
\item the torus $\Tbf$ acting by right translations,
\item the torus $\Tbf \times \Tbf$ acting by simultaneous left and right translations.
\end{enumerate}
Taking monodromic sheaves with respect to these actions defines three possible versions of the Hecke category. We show that they are all canonically isomorphic in \ref{sectionDefiHecke}. We denote the resulting category by $\HH$. 

We then equip this category with a monoidal structure given by convolution product. Consider the correspondence 
\[\begin{tikzcd}
	& {\Ubf \backslash \Gbf \times^{\Ubf} \Gbf/\Ubf} & {\Ubf \backslash \Gbf/\Ubf} \\
	{\Ubf \backslash \Gbf/\Ubf} && {\Ubf \backslash \Gbf/\Ubf}
	\arrow["{p_1}", from=1-2, to=2-1]
	\arrow["{p_2}"', from=1-2, to=2-3]
	\arrow["m", from=1-2, to=1-3]
\end{tikzcd}\]
and let $A, B \in \HH$, using the model of $\HH$ as sheaves that are (twisted) equivariant for the action of $\Tbf \times \Tbf$, we define
\begin{equation*}
A * B = \For_{\RR_{\Tbf^2}}m_!(A \widehat{\boxtimes}_{\Lambda} B)[\dim \Tbf]
\end{equation*}
where $A \widehat{\boxtimes} B$ is obtained from $A \boxtimes_{\Lambda} B$ by completing along an ideal. We then push it forward along $m$ and forget two out of the four copies of $\Tbf$. The details of this construction are given in Section \ref{sectionMonoidalStructures}. The key difference with the construction of \cite{BezYun} is that we have restored the symmetry between the left and right $\Tbf$-actions and the proof that this construction defines a monoidal structure is now immediate.
We conclude our study by showing the following properties.
\begin{thm}
\begin{enumerate}
\item Lemma \ref{lemCompactGeneration} :  The category $\HH$ is compactly generated.
\item Theorem \ref{thmQuasiRigidity} :  The category $\HH$ is quasi-rigid and is a direct sum of rigid categories. 
\item Corollary \ref{corolPivotalStructure} :  The category $\HH$ is equipped with a canonical pivotal structure. 
\end{enumerate}
\end{thm}
These structure results are shown as in \cite{BenZviNadlerComplexGroup} with modifications to accommodate the free monodromic sheaves.

\subsection{Categorical trace, representations and character sheaves}

From now on, we specialize $k$ to be $\Fqb$. Once we have at our disposal the free monodromic Hecke category, we compute the categorical traces and centers of both the identity and the Frobenius on them, see Definition \ref{defi:CategoricalTrace}. These $\infty$-categorical versions have classical analogs, see \cite{ZhuTrace}. Note that we will focus on the trace of this category as, by a theorem of \cite{BenZviNadlerComplexGroup}, the center and trace (resp. $\Frob$-center and $\Frob$-trace) on $\HH$ are isomorphic. 

In the case of the Frobenius, we show the following theorem.
\begin{thm}[Theorem \ref{thmTraceFrobenius}]
There is a canonical isomorphism
\begin{equation*}
\Tr(\Frob_*, \HH) = \DD(\frac{\Gbf}{\Ad_{\Frob}\Gbf}, \Lambda).
\end{equation*}
\end{thm}
The category appearing in the RHS is the category of ind-constructible sheaves of $\Lambda$-modules on the quotient stack of $\Gbf$ acting on itself by 
$$g.x = gx\Frob(g)^{-1}.$$
This statement holds for all reductive groups $\Gbf$ (even if $\Gbf$ is disconnected). If $\Gbf$ is connected, the Lang map provides an isomorphism $\frac{\Gbf}{\Ad_{\Frob}\Gbf} = \pt/\Gbf^{\Frob}$. By étale descent, the category of sheaves on this finite stack is equivalent to the (unbounded, derived) category of representations of $\Gbf^{\Frob}$ on $\Lambda$-modules. 

Lusztig, in \cite{LusztigCatCenter1} \cite{LusztigCatCenter2}, has shown some analogous statement after taking certain subquotients on the Hecke category side (given by truncating by a two sided cell in the Weyl group). The main differences are that our result remains valid for modular and integral representations. Our proof is closely related to \cite{BenZviNadlerComplexGroup} and does not require the classification of irreducible representations established in \cite{LusztigBook}. Similarly, Zhu \cite{ZhuTrace} has shown some analogous statement using a $1$-categorical definition of categorical traces. Our approach is conceptually much closer to the one of Zhu.

 We expect many applications of this result. In our companion paper \cite{EteveJordanDecomp}, we plan to use this construction to prove a general version of the Jordan decomposition for representations of finite reductive groups, extending the result of \cite{LusztigYun} to the modular setting. 

In the case of the trace of the identity, there are many results in the literature that assert `the center of the Hecke category is the category of character sheaves', see \cite{BenZviNadlerComplexGroup}, \cite{BezrukCatCenter}, \cite{BFOCharSheaves}, \cite{LusztigCSCatCenter}, \cite{ZhuTrace}, where the category of character sheaves is a certain full subcategory of $\DD(\frac{G}{\Ad(G)}, \Lambda)$ first constructed by Lusztig in \cite{LusztigCS1}. Since we work with free monodromic categories, we expect, as in \cite{BezrukavnikovTomalchov}, to deal with `free monodromic character sheaves'. The problems coming from the free monodromic aspects did not arise in the case of the trace of the Frobenius as the monodromy is essentially killed by the Frobenius. In Section \ref{sectionFreeMonocharacterSheaves}, we setup the requirements to define and study the category of free monodromic character sheaves. 
\begin{defi}
We define the category of free-monodromic character sheaves as the categorical center of $\HH$, namely,
\begin{equation*}
\CS^{\wedge} = \Zcal(\HH),
\end{equation*}
where $\CS^{\wedge}$ stands for the category of free-monodromic character sheaves.
\end{defi}
We have chosen here to first define the category of character sheaves as the categorical center of $\HH$ and then to identify it with a certain subcategory of sheaves on $\frac{\Gbf}{\Ad(\Gbf)}$. We then show the following identification theorem. 
\begin{thm}[Theorem \ref{thmTraceIdentity}]
There is a canonical isomorphism 
\begin{equation*}
\CS^{\wedge} = A\--\mathrm{mod}(\Pre\CS^{\wedge}),
\end{equation*}
where $A$ is an algebra in $\Pre\CS^{\wedge}$ which is a certain full subcategory of $\DD(\frac{\Gbf}{\Ad(\Gbf)}, \RR_{\Tbf \times \Tbf})$ which we define in Section \ref{sectionFreeMonocharacterSheaves}. 
\end{thm}
In the case where $G$ is connected and has connected center, we then identify the algebra $A$ in Theorem \ref{thmUnitCenter} and we then recover the classical story of character sheaves, as well as the completed version of \cite{BezrukavnikovTomalchov}. 

There are natural maps from a category to its trace and from the center of a category to the category. That is, we have 
$$\Zcal(\HH) \to \HH \to \Tr(\id,\HH)$$
and 
$$\Zcal_{\Frob}(\HH) \to \HH \to \Tr(\Frob_*, \HH)$$
(we refer to Section \ref{sectionCategoricalTraces} for the notations). We show that these maps are naturally described in terms of pull-push along the following correspondences 
\begin{equation*}
\frac{\Gbf}{\Ad(\Gbf)} \leftarrow \frac{\Gbf}{\Ad(\Bbf)} \to \frac{\Ubf \backslash \Gbf/\Ubf}{\Ad(\Tbf)} \leftarrow \Ubf \backslash \Gbf/\Ubf, 
\end{equation*}
and 
\begin{equation*}
\frac{\Gbf}{\Ad_{\Frob}(\Gbf)} \leftarrow \frac{\Gbf}{\Ad_{\Frob}\Bbf} \to \frac{\Ubf \backslash \Gbf/\Ubf}{\Ad_{\Frob}(\Tbf)} \leftarrow \Ubf \backslash \Gbf/\Ubf. 
\end{equation*}
These correspondences are known as the horocycle correspondence and they were introduced by Lusztig in \cite{LusztigCS1} (though in non stacky form). We review their formalism in Section \ref{sectionHorocycle}. Note that a $\Bbf$-version of this construction was also considered in \cite{BonnafeRouquiuerDudas}.

Let us finally highlight how Deligne--Lusztig theory fits into this formalism. We consider the stack $\frac{\Ubf \backslash \Gbf/\Ubf}{\Ad_{\Frob}\Tbf}$ and we equip it with its Bruhat stratification. For an element $w \in \Wbf$, there is an isomorphism (depending on some lift $\dot{w}$ of $w$)
\begin{equation*}
\frac{\Ubf \backslash \BwB/\Ubf}{\Ad_{\Frob}\Tbf} = \pt/(\Tbf^{w\Frob} \rtimes (\Ubf \cap {^w}\Ubf)). 
\end{equation*}
As $\Ubf \cap {^w}\Ubf$ is connected, the category of sheaves on this stratum is equivalent to the category of representations of the finite group $\Tbf^{w\Frob}$. Moreover the collection of all these finite tori are exactly the ones appearing in the original construction of \cite{DeligneLusztig}. In Section \ref{sectionHorocycle}, we explain how to recover the usual Deligne--Lusztig induction and restriction functors out of this geometry. It should be noted that the category of sheaves on the stack $\frac{\Ubf \backslash \Gbf/\Ubf}{\Ad_{\Frob}\Tbf}$ is glued in a nontrivial way from the categories of representations of the finite tori $\Tbf^{w\Frob}$ contained in $\Gbf$. It does not seem to us that this construction was previously studied in the literature and we plan to come back to this in future work. As a first application, we are able to give a very short proof of the fundamental theorem of Deligne and Lusztig, that all irreducible representations of $\Gbf^{\Frob}$ appear in the cohomology of the Deligne--Lusztig varieties, see Section \ref{subsectionConvolutionPatterns}. 

\subsection{Convention and notations}

\subsubsection{$\infty$-categories} We will essentially work with the formalism of $\infty$-categories of \cite{LurieA} and most of our categories will be stable. We denote by $\Pr$ and $\Pr_{\Lambda}$ the category of presentable categories and the category of $\Lambda$-linear presentable linear categories, see \cite{LurieB} for the definition. We denote by $\Delta$ the simplex category. For a (dg)-ring $A$, we denote by $\DD(A)$ the $\infty$-derived category of $A$. We denote by $\ho$ the homotopy category of an $\infty$-category, so that $\ho(\DD(A))$ is the usual derived category of $A$. For a ring object $A$ in some monoidal category $\Ccal$, we denote by $A\--\mathrm{mod}(\Ccal)$ the category of modules over it. We will repeatedly use the Barr-Beck-Lurie theorem \cite[Theorem 4.7.0.3]{LurieB} and the formalism of monadic functors, hence whenever we invoke this theorem, we mean the one of \emph{loc. cit}. For $\Lambda$-linear categories, unless specified, the functor $\Hom$ always denotes the mapping complex in $\DD(\Lambda)$ and $\Hom^i$ is its $i$-th cohomology group. 

\subsubsection{Geometry} We fix an algebraically closed field $k$ and all (unless specified) schemes and stacks are of finite type over $k$. In Sections \ref{sectionHorocycle} and \ref{sectionFreeMonocharacterSheaves}, we assume that $k$ is algebraic closure of finite field.

\subsubsection{Sheaves} We fix once and for all prime $\ell \neq p = \mathrm{char}(k)$. We will denote by $\Lambda \in \{\Flb, \Zlb, \Qlb\}$ a coefficient ring. We denote by $\DDic(X, \Lambda)$ and $\DDc(X, \Lambda)$ the category of ind-constructible sheaves and constructible sheaves of $\Lambda$-modules respectively on the stack $X$. In particular, for a finite group $\Gamma$, there is a canonical isomorphism $\DDic(\pt/\Gamma, \Lambda) = \DD(\Rep_{\Lambda}\Gamma)$. 

\subsubsection{Group} We fix once and for all a (possibly disconnected) reductive group $\Gbf$ over $k$ and when $k$ is the algebraic closure of finite field, we fix a Frobenius endomorphism $\Frob : \Gbf \to \Gbf$, that is a purely inseparable isogeny such that a power of it is a Frobenius coming from some $\Fq$-structure. We denote by $\Gbf^{\circ}$ the neutral component of $\Gbf$. We fix a Borel pair $\Bbf = \Tbf\Ubf$ of $\Gbf^{\circ}$. If $\Gbf$ is equipped with a Frobenius we assume that $(\Bbf, \Tbf)$ is $\Frob$-stable. 

\subsubsection{Condensed mathematics} We will consider condensed rings. These are sheaves of rings on the proétale site of the point $\pt_{\proet}$. They will be used in the definition of categories of free monodromic sheaves sheaves, we will systematically refer to either \cite{BhattScholze} or \cite{HemoRicharzScholbach} for their use.

\subsection{Acknowledgments}

Part of the results of this paper where in the thesis of the author \cite{eteveThesis} but several results have been added : mostly the $\Qlb$-version of the free monodromic sheaves and the part on free monodromic character sheaves. This second set of results was worked out while the author was a guest at the Max Planck Institute for Mathematics Bonn. We would like to thank our advisor Jean-François Dat for the many discussions and continuous support over the years. We thank Maxime Ramzi, Olivier Dudas, Simon Riche, Cédric Bonnafé, Colton Sandvik, Jens Eberhardt and Tom Gannon for discussions and advice. We thank Thibaud van der Hove, Rizacan \c Cilo\v glu and Xinwen Zhu for pointing out some typos and mistake in an earlier version of this paper.

\section{The free monodromic formalism}\label{sectionFreeMonoFormalism}

For the rest of this section, we fix $k$ an algebraically closed field of characteristic $p > 0$ and $\Tbf$ a torus over $k$. 

\subsection{The rings}\label{sectionTheRings}

We denote by $X^*(\Tbf)$ and $X_*(\Tbf)$ the character and cocharacter lattices of $\Tbf$. We denote by $\pi_1(\Tbf)$ the étale fundamental group of $\Tbf$ at the base point $1 \in \Tbf$. We denote by $\pi_1^t(\Tbf)$ the tame quotient, that is, the largest quotient of pro-order prime to $p$. It is known that 
\begin{equation*}
\pi_1^t(\Tbf) = \varprojlim_{(n,p) = 1} \Tbf[n] = X_*(\Tbf) \otimes \pi_1^t(\Gm).
\end{equation*}
We also denote by $\pi_1(\Tbf)^{\wild}$ the kernel of the projection $\pi_1(\Tbf) \rightarrow \pi_1^t(\Tbf)$. 

Let $\ell$ be a prime different from $p$, as in \cite{GabberLoeser} we define $\pi_1(\Tbf)_{\ell}$ the largest pro-$\ell$ quotient of $\pi_1^t(\Tbf)$. We define the following condensed rings 
\begin{enumerate}
\item if $\Lambda \in \{\Fl, \Zl\}$, then we set 
\begin{equation*}
\RR_{\Tbf, \Lambda} = \Lambda \llbracket \pi_1(\Tbf)_{\ell} \rrbracket = \varprojlim_n \Lambda[\Tbf[\ell^n]].
\end{equation*}
\item if $\Lambda = \Ql$, then let $I_{\Zl}$ be the augmentation ideal of $\RR_{\Tbf, \Zl}$, we set
\begin{equation*}
\RR_{\Tbf, \Ql} = (\RR_{\Tbf, \Zl}[\frac{1}{\ell}])^{\wedge} = \varprojlim_n (\RR_{\Tbf, \Zl}[\frac{1}{\ell}])/I_{\Zl}^n,
\end{equation*}
where $(-)^{\wedge}$ denotes the completion at the ideal $I_{\Zl}$. 
\end{enumerate}

\begin{rque}
For $\Lambda \in \{\Fl, \Zl\}$, the ring $\RR_{\Tbf, \Lambda}$ is the condensed ring attached to the topological ring defined in the same way. For $\Lambda = \Ql$, the underlying ring $\RR_{\Tbf, \Ql}(*)$ is $(\RR_{\Tbf, \Zl}[\frac{1}{\ell}])^{\wedge}$ where $\RR_{\Tbf, \Zl}[\frac{1}{\ell}]$ is considered as a discrete ring. 
\end{rque}

\begin{rque}\label{rkDescription}
If $\Lambda \in \{\Fl, \Zl \}$, after fixing a topological generator of $\pi_1^t(\Gm)$ and a basis of $X_*(\Tbf)$, we get isomorphisms
\begin{equation*}
\RR_{\Tbf, \Lambda} = \Lambda\llbracket t_1, \dots, t_n \rrbracket. 
\end{equation*}
\end{rque}

If $\Lambda$ is a $\Lambda_0 \in \{\Fl, \Zl, \Ql\}$-algebra, then we denote by $\RR_{\Tbf, \Lambda} = \RR_{\Tbf, \Lambda_0} \otimes_{\Lambda_0} \Lambda$. 

\begin{defi}[\protect{\cite[Definition 6.1]{HemoRicharzScholbach}}]
A condensed ring $R$ is t-admissible if $R(*)$ is regular coherent and for all extremally disconnected sets $S$ the maps $R(*) \rightarrow \Gamma(S, R)$ is flat. 
\end{defi}

\begin{lem}\label{lemtAdmissible}
For $\Lambda \in \{\Fl, \Zl, \Ql, \Flb, \Qlb, \Zlb\}$ the rings $\Lambda$ and $\RR_{\Tbf, \Lambda}$ are t-admissible.
\end{lem}

\begin{proof}
The case of $\Lambda$ is done in \cite{HemoRicharzScholbach}. For $\RR_{\Tbf, \Lambda}$, we first consider the case $\Lambda \in \{\Fl, \Zl\}$. By Remark \ref{rkDescription}, the ring $\RR_{\Tbf,\Zl}$ is clearly regular and noetherian, in particular coherent. Let $S$ be an extremally disconnected set and write it $S = \varprojlim_i S_i$ as a limit of finite sets. Then we have
\begin{equation*}
\Gamma(S, \RR_{\Tbf, \Zl}) = \varprojlim_{n,m} \varinjlim_i \Gamma(S_i, \Z/\ell^n\Z[T[m]]).
\end{equation*}
Each $\Gamma(S_i, \Z/\ell^n\Z[T[m]])$ is flat over $\Z/\ell^n\Z[T[m]]$ and therefore so is $\varinjlim_i \Gamma(S_i, \Z/\ell^n\Z[T[m]])$. By \cite[Tag 0912]{Stacks}, we have that $\RR_{\Tbf,\Zl} \rightarrow \Gamma(S, \RR_{\Tbf,\Zl})$ is flat. The same argument holds for $\RR_{\Tbf, \Lambda}$ for $\Lambda \in \{\O_E, \F_{\ell^n}\}$, where $E/\Ql$ is a finite extension with ring of integers $\O_E$. We then have $\RR_{\Tbf, \Zlb} = \varinjlim_E \RR_{\Tbf, \O_E}$ and $\RR_{\Tbf, \Flb} = \varinjlim_n \RR_{\Tbf, \F_{\ell^n}}$. Let $S$ be extremally disconnected, we then have 
\begin{equation*}
\Gamma(S, \RR_{\Tbf, \Flb}) = \varinjlim_n \Gamma(S,\RR_{\Tbf, \F_{\ell^n}})
\end{equation*}
and since each $\Gamma(S,\RR_{\Tbf, \F_{\ell^n}})$ is flat over $\RR_{\Tbf, \F_{\ell^n}}$, we deduce that $\Gamma(S, \RR_{\Tbf, \Flb})$ is flat over $\RR_{\Tbf, \Flb}$. The same argument holds for $\RR_{\Tbf, \Zlb}$. 

We now consider the case of $\Lambda = \Ql$. By Remark \ref{rkDescription}, it is clear that the underlying ring of $\RR_{\Tbf, \Ql}$ is regular. We only need to check the flatness condition. Let $S = \lim_i S_i$ be extremally disconnected. We have 
\begin{align*}
\Gamma(S, \RR_{\Tbf, \Ql}) &=  \Gamma(S, \varprojlim_n \RR_{\Tbf, \Zl}[\frac{1}{\ell}]/I^n) \\
&= \varprojlim_n \Gamma(S, \RR_{\Tbf, \Zl}[\frac{1}{\ell}]/I^n)\\
&= \varprojlim_n \Gamma(S, (\RR_{\Tbf, \Zl}/I^n)[\frac{1}{\ell}])\\
&= \varprojlim_n \Gamma(S, \RR_{\Tbf, \Zl}/I^n)[\frac{1}{\ell}].
\end{align*}
The second line comes from the fact that $\Gamma(S, -)$ commutes with limits  as it is a right adjoint. The third one comes from the commutation of quotient and localizations and the last one comes from the fact that taking $[\frac{1}{\ell}]$ is a filtered colimit and $S$ is compact. Now each $\Gamma(S, \RR_{\Tbf, \Zl}/I^n)$ is flat over $\RR_{\Tbf, \Zl}/I^n$ and $\Gamma(S, \RR_{\Tbf, \Zl}/I^n)[\frac{1}{\ell}]$ is flat over $(\RR_{\Tbf, \Zl}/I^n)[\frac{1}{\ell}]$. Again, by \cite[Tag 0912]{Stacks}, we conclude that $\Gamma(S, \RR_{\Tbf, \Ql})$ is flat over $\RR_{\Tbf, \Ql}$. The same argument holds if we replace $\Ql$ by a finite extension $E$ and $\RR_{\Tbf, \Qlb} = \varinjlim_{E} \RR_{\Tbf, E}$. 
\end{proof}

\subsection{Categories of sheaves}\label{subsectionCategoriesOfScheaves}

We fix $\Lambda \in \{\Fl, \Zl, \Ql, \Flb, \Qlb, \Zlb\}$. We recall the definition of ind-constructible sheaves \cite{HemoRicharzScholbach}.

\begin{defi}[Constructible and Lisse sheaves]
Let $R$ be a condensed ring a sheaf $A \in \DD(X_{\proet}, R)$ is lisse if it is dualizable. The sheaf $A$ is constructible if there is a stratification $X = \sqcup X_i$ such that the restriction of $A$ to each strata is lisse. We denote by $\DD_{\lis}(X, R) \subset \DDc(X, R) \subset \DD(X_{\proet}, R)$ the full subcategories of lisse and constructible sheaves.
\end{defi}

\begin{defi}
Let $R$ be a condensed ring, we denote by $\DDic(X, R)$ the full subcategory of $\DD(X_{\proet}, R)$ of sheaves that are isomorphic to a filtered colimit of constructible sheaves. There is a natural map $\Ind(\DDc(X, R)) \rightarrow \DD(X, R)$. 
\end{defi}

\begin{lem}[\protect{\cite[Corollary 8.3]{HemoRicharzScholbach}}]
Assume that $X$ is proétale locally uniformly $R$-cohomologically bounded then the map $\Ind(\DDc(X, R)) \rightarrow \DD(X, R)$ is an isomorphism. 
\end{lem}

\begin{lem}
All schemes of finite type over $k$ are $\RR_{\Tbf, \Lambda}$ uniformly cohomologically bounded for $\Lambda \in \{\Fl, \Zl, \Ql, \Flb, \Zlb, \Qlb\}$. 
\end{lem}

\begin{proof}
The case where $\Lambda \in \{\Fl, \Zl\}$ follows from classical limit arguments as in \cite[Lemma 8.6]{HemoRicharzScholbach}. The case of $\Lambda \in \{\Flb, \Zlb\}$ follows from a filtered colimit argument. The case of $\Lambda = \Qlb$ follows from the case $\Lambda = \Ql$ using the same colimit argument. If $\Lambda = \Ql$, then if $U \rightarrow X$ is an affine pro-étale scheme over $X$, then we have 
\begin{align*}
\RGamma(U, \RR_{\Tbf, \Ql}) &= \varprojlim_{n} \RGamma(U, \RR_{\Tbf, \Ql}/I^n) \\
&=\varprojlim_{n}\RGamma(U, \RR_{\Tbf, \Zl}/I^n[\frac{1}{\ell}]) \\
&=\varprojlim_{n}\RGamma(U, \RR_{\Tbf, \Zl}/I^n)[\frac{1}{\ell}]. 
\end{align*}
Since $\RR_{\Tbf, \Zl}/I^n$ is a finite flat $\Zl$-algebra, there exists $d \geq 0$ such that $\RGamma(U, \RR_{\Tbf, \Zl}/I^n) \in \DD^{[0, d]}(\Zl)$. Finally, as $\varprojlim_n$ has cohomlogical dimension at most $1$, the resulting complex is in degrees $[0, d+1]$.
\end{proof}

\begin{thm}
For $\DD \in \{\DDc, \DDic\}$ there exists a $6$-functors formalism on schemes such that 
\begin{equation*}
\DD(X) = \DD(X, \RR_{\Tbf, \Lambda}). 
\end{equation*}
Moreover this $6$-functors formalism satisfies $!$ and $*$-descent along smooth maps, see Definition \ref{defiDescentFor6Functors}. In particular $X \mapsto \DDic(X, \RR_{\Tbf, \Lambda})$ defines a $6$-functors formalism on algebraic stacks. 
\end{thm}

\begin{proof}
The case of étale maps is done in \cite[Corollary 8.7.]{HemoRicharzScholbach}. It is then enough to do the case of a projection $f : X \times \A^1 \to X$. Denote by $f^n : X \times \A^n \to X$ the $n$-fold product of this map, which is again a projection. For all objects in $A \in \DD(X, \Lambda)$, the obvious maps $A \to f_{n,*}f_n^*A$ and $f_{n,!}f_n^!A \to A$ are isomorphisms. Hence by \cite[Proposition 6.18]{SixFunctors}, $f$ is of universal $!$ and $*$-descent. 
\end{proof}

\begin{rque}\label{rqueFiniteExpansion} 
It is not necessary to restrict ourselves to the category of schemes of finite type over $k$. By \cite{SixFunctors}, this six functors formalism is defined on the larger category of qcqs schemes over $k$ and where the maps for which the $!$-functors are defined are the maps of `finite expansion'. These are the maps $f : X \to Y$ such that on open affines $\Spec(A) \subset X \to \Spec(B) \subset Y$ the exists a finite set of elements $X_1, \dots, X_n \in B$ such that $B[X_1, \dots, X_n] \to A$ is integral. This includes morphisms that are not of finite type. 
\end{rque} 

\begin{thm}
Let $X$ be an algebraic stack, there exists two t-structures on $\DDc(X, \RR_{\Tbf, \Lambda})$ called the standard and perverse t-structure. They are characterized by 
\begin{enumerate}
\item for the standard one : for all geometric points $x : \Spec(k) \rightarrow X$ the functor $x^* : \DDc(X, \RR_{\Tbf}) \rightarrow \Perf_{\RR_{\Tbf,\Lambda}(*)}$ is $t$-exact.
\item for the perverse one : $A$ is in $\DDc(X, \RR_{\Tbf, \Lambda})^{\geq^p 0}$ if and only if
\begin{equation*}
\forall i \ \dim supp H^{-i} \leq i .
\end{equation*}
\end{enumerate}
\end{thm}

\begin{proof}
The case of the standard t-structure follows from \cite{HemoRicharzScholbach} and Lemma \ref{lemtAdmissible}. The passage from the standard t-structure to the perverse one is standard and done as in \cite{BBD}.
\end{proof}

The natural map $\RR_{\Tbf, \Zlb} \rightarrow \RR_{\Tbf, \Qlb}$ induces a functor for all stacks $X$, 
\begin{align*}
\DDic(X, \RR_{\Tbf, \Zlb}) &\rightarrow \DDic(X, \RR_{\Tbf, \Qlb}) \\
M &\mapsto M \otimes_{\RR_{\Tbf, \Zlb}} \RR_{\Tbf, \Qlb}. 
\end{align*}
This functor commutes with the $6$-functors and preserves constructibility. 

Since the ring $\RR_{\Tbf, \Lambda}$ is regular coherent, the categories $\Perf(\RR_{\Tbf, \Lambda}(*))$ and $\DD_{\lis}(X, \RR_{\Tbf, \Lambda})$ are equipped with a $t$-structure coming from the restriction the standard $t$-structure by \cite{HemoRicharzScholbach}. We call the heart $\DD_{\lis}(X, \RR_{\Tbf, \Lambda})^{\heartsuit}$ of this category the category of local systems. Let $X$ be a scheme and $\bar{x}\rightarrow X$ a geometric point of $X$.
\begin{lem}
Assume that $X$ is a qcqs normal connected scheme, then $\DD_{\lis}(X, \RR_{\Tbf, \Lambda})^{\heartsuit}$ is equivalent to the category of continuous representations of $\pi_1^{et}(X, \bar{x})$ on finite type $\RR_{\Tbf, \Lambda}(*)$-modules.
\end{lem} 

\begin{proof}
By \cite[Proposition 5.2]{HemoRicharzScholbach}, since $X$ is qcqs, $\DD_{\lis}(X, -)$ commutes with filtered colimits of rings. Hence we can assume that $\Lambda \in \{\Zl, \Fl, \Ql\}$. The case of $\Lambda \in \{\Fl, \Zl\}$ follows from a standard limit argument and the case of a discrete ring \cite[Corollary 5.1.5]{BhattScholze}. Similarly, since inverting $\ell$ is a $t$-exact localization by \cite[Proposition 5.5]{HemoRicharzScholbach} we deduce from the $\RR_{\Tbf, \Zl}$-case the $\RR_{\Tbf, \Zl}[\frac{1}{\ell}]$-case. The $\RR_{\Tbf, \Ql}$-case then follows from a limit argument as before. 
\end{proof}

\subsection{The universal multiplicative character}

We fix $\Lambda \in \{\Qlb, \Zlb, \Flb\}$ and from now on, if the context is clear, we set $\RR_{\Tbf, \Lambda} = \RR_{\Tbf}$. 

There is a canonical morphism 
\begin{equation}\label{canon}
can : \pi_1(\Tbf) \rightarrow \RR_{\Tbf}^{\times}
\end{equation}
which defines an $\RR_{\Tbf}$-rank one local system on $\Tbf$ which we denote by $L_{\Tbf}$. We will denote $L_{\Tbf}$ by $L_{\Tbf, \Lambda}$ if we want to put emphasis on the coefficient. Note that $L_{\Tbf, \Ql} = L_{\Tbf, \Zl} \otimes_{\RR_{\Tbf, \Zl}} \RR_{\Tbf, \Ql}$.

\begin{lem}[\cite{GabberLoeser} 3.1]\label{lemGL}
Let $\alpha : \Tbf \rightarrow \Tbf'$ be a morphism of tori, it induces a morphism $p_* : \RR_{\Tbf} \rightarrow \RR_{\Tbf'}$. 
\begin{enumerate}
\item $\alpha^*L_{\Tbf'} = L_{\Tbf} \otimes_{\RR_{\Tbf}} \RR_{\Tbf'}$.
\item Assume $\alpha$ is a quotient map of relative dimension $d$ then we have $\alpha_!L_{\Tbf} = L_{\Tbf'}[-2d](-d)$.
\end{enumerate}
\end{lem}

\begin{proof}
The proof in the case $\Lambda \in \{\Flb, \Zlb\}$ comes from \cite{GabberLoeser}. The case $\Lambda = \Qlb$ follows from compatibility with change of coefficients $\otimes_{\RR_{\Tbf, \Zl}} \RR_{\Tbf, \Ql}$.
\end{proof}

\begin{rque}
In \emph{loc. cit.}, a quotient map means the projection on a direct factor.
\end{rque}

\begin{defi}
We denote by $\Ch_{\Lambda}(\Tbf)$ 
\begin{enumerate}
\item if $\Lambda \in \{\Flb, \Zlb\}$, the set of characters $\pi_1^t(\Tbf) \rightarrow \Lambda^{\times}$ of finite order prime to $\ell$.
\item if $\Lambda = \Qlb$, the set of all characters $\pi_1^t(\Tbf) \rightarrow \Lambda^{\times}$ of finite order.
\end{enumerate}
\end{defi} 

To each $\chi \in \Ch_{\Lambda}(\Tbf)$ there is a corresponding rank one local system $\mathcal{L}_{\chi}$. These local systems are usually called Kummer local systems. We denote by $L_{\Tbf, \chi} = L_{\Tbf} \otimes_{\Lambda} \mathcal{L}_{\chi}$. Each of those is an $\RR_{\Tbf, \Lambda}$-local system locally free of rank one. 

\begin{rque}
Let $\Lambda \in \{\Flb, \Zlb, \Qlb\}$ and let $\Lambda_{0} \in \{\Fl, \Zl, \Ql\}$ be the corresponding subring. Let $\Tbf^{\vee}_{\Lambda}$ denote the torus dual to $\Tbf$ over $\Lambda$, then $\Ch(\Tbf, \Lambda) \subset \Tbf^{\vee}_{\Lambda}(\Lambda)$ is isomorphic to a subset of the set of $\Lambda$-points of $T^{\vee}_{\Lambda}$. The ring $\RR_{\Tbf, \Lambda}$ is isomorphic to $\Tbf^{\vee, \wedge}_{\Lambda_0, 1} \otimes_{\Lambda_0} \Lambda$, the completion along the unit section of $\Tbf^{\vee}_{\Lambda_0}$ then scalar extended to $\Lambda$. 
\end{rque}

\begin{lem}\label{lemMultiplicative}
For $\chi \in \Ch_{\Lambda}(\Tbf)$, the $\RR_{\Tbf, \Lambda}$-local system $L_{\Tbf, \chi}$ is mutiplicative in the sense of Definition \ref{defiMultiplicativeSheaf}.
\end{lem}

\begin{proof}
It is well known that $\Lcal_{\chi}$ is multiplicative. The fact that $L_{\Tbf}$ is multiplicative follows from Lemma \ref{lemGL} and the following computation 
\begin{align*}
\pr_1^*L_{\Tbf} \otimes_{\RR_{\Tbf}} \pr_2^*L_{\Tbf} &= (L_{\Tbf \times T} \otimes_{\RR_{\Tbf \times \Tbf, \pr_{1,*}}} \RR_{\Tbf}) \otimes_{\RR_{\Tbf}} (L_{\Tbf \times \Tbf} \otimes_{\RR_{\Tbf \times \Tbf, \pr_{2,*}}} \RR_{\Tbf}) \\
&= L_{\Tbf \times \Tbf} \otimes_{\RR_{\Tbf \times \Tbf}} (\RR_{\Tbf} \otimes_{\RR_{\Tbf}} \RR_{\Tbf}) \\
&= L_{\Tbf \times \Tbf} \otimes_{\RR_{\Tbf \times \Tbf}, m_*} \RR_{\Tbf},
\end{align*}
where $\pr_i, m : \Tbf \times \Tbf \to \Tbf$ denote the projections and the multiplication respectively.
\end{proof}

\begin{lem}\label{lemFunctorialityPrimeToL}
Let $\alpha : \Tbf \rightarrow \Tbf'$ be a finite étale isogeny. The following holds 
\begin{enumerate}
\item if $\Lambda \in \{\Flb, \Zlb\}$ and $\alpha$ has order prime to $\ell$ then the map $\RR_{\Tbf, \Lambda} \rightarrow \RR_{\Tbf', \Lambda}$ is an isomorphism, 
\item if $\Lambda = \Qlb$ then the map $\RR_{\Tbf, \Qlb} \rightarrow \RR_{\Tbf', \Qlb}$ is an isomorphism. 
\end{enumerate} 
\end{lem}

\begin{proof}
For $(i)$, if $\alpha$ has order prime to $\ell$ then $\alpha$ induces an isomorphism $\Tbf[\ell^n] \rightarrow \Tbf'[\ell^n]$ for all $n \geq 0$ hence the induced map $\RR_{\Tbf, \Lambda} \rightarrow \RR_{\Tbf', \Lambda}$ is an isomorphism for $\Lambda \in \{\Flb, \Zlb\}$. 

For $(ii)$, we denote by $I$ and $I'$ the augmentation ideals of $\RR_{\Tbf, \Zlb}[\frac{1}{\ell}]$ and $\RR_{\Tbf', \Zlb}[\frac{1}{\ell}]$. There exists $n > 0$ such that 
$$I'^n \subset I\RR_{\Tbf', \Zlb}[\frac{1}{\ell}] \subset I'$$ 
and we have 
$$I \subset \RR_{\Tbf, \Zlb}[\frac{1}{\ell}] \cap I\RR_{\Tbf', \Zlb}[\frac{1}{\ell}] \subset I'\cap \RR_{\Tbf', \Zlb}[\frac{1}{\ell}] = I.$$
As a consequence of the first inclusion, it follows that the completion along the $I'$-adic topology or the $I$-adic topology of $\RR_{\Tbf', \Zlb}[\frac{1}{\ell}]$ coincide. By the Artin-Rees Lemma \cite[Tag 00IN]{Stacks}, the completion of $\RR_{\Tbf, \Zlb}[\frac{1}{\ell}]$ along the $I$-adic or $\RR_{\Tbf, \Zlb}[\frac{1}{\ell}] \cap I\RR_{\Tbf', \Zlb}[\frac{1}{\ell}]$-adic topology coincide. All inclusions in the second line show that the map $\RR_{\Tbf, \Qlb}/I \rightarrow \RR_{\Tbf', \Qlb}/I'$ is an isomorphism, hence the map $\RR_{\Tbf,\Qlb} \rightarrow \RR_{\Tbf',\Qlb}$ is an isomorphism. 
\end{proof}

\begin{rque}
Note that $\RR_{\Tbf, \Qlb}$ is, up to scalar extension from $\Ql$ to $\Qlb$, the completion of $\RR_{\Tbf, \Zlb}$ at the point corresponding to $1 \in \Spec(\RR_{\Tbf, \Zlb}[\frac{1}{\ell}]) = \Tbf^{\vee, \wedge}_{\Zlb, 1}[\frac{1}{\ell}]$. More generally for any point $\chi \in  \Tbf^{\vee, \wedge}_{\Zlb, 1}[\frac{1}{\ell}]$, we can define the completion at $\chi$ of $\RR_{\Tbf, \Zlb}[\frac{1}{\ell}]$, let us denote this ring by $\RR_{\Tbf, \Qlb, \chi}$, all the rings $\RR_{\Tbf, \Qlb, \chi}$ are isomorphic as topological rings and there is a natural map $\iota_{\chi} : \RR_{\Tbf, \Zlb} \rightarrow \RR_{\Tbf, \Qlb, \chi}$. Upon identifying $\RR_{\Tbf, \Qlb, \chi}$ with $\RR_{\Tbf, \Qlb}$ by translation by $\chi$, there is an isomorphism of sheaves on $\Tbf$ : $$ L_{\Tbf, \Zlb} \otimes_{\RR_{\Tbf, \Zlb}, \iota_{\chi}} \RR_{\Tbf, \Qlb} = L_{\Tbf, \Qlb} \otimes_{\Qlb} \Lcal_{\chi}.$$ 
\end{rque}

\begin{lem}\label{lemFunctorialityOrderL}
Let $f : \Tbf \to \Tbf'$ be a finite isogeny whose kernel has order a power of $\ell$. Let $\Lambda \in \{\Zlb, \Flb\}$. The map $f$ induces a map of algebras $f : \RR_{\Tbf} \to \RR_{\Tbf'}$. There is a canonical isomorphism of algebras 
$$\RR_{\Tbf'} \otimes_{\RR_{\Tbf}} \Lambda = \Lambda[\ker(f)].$$
\end{lem}

\begin{proof}
The algebra on the LHS can be interpreted as the fiber at $1$ of $f_*(L_{\Tbf} \otimes_{\RR_{\Tbf}} \Lambda) = f_*\Lambda$. 
\end{proof}

\subsection{Categories of free monodromic sheaves}

Let $X$ be a scheme with a $\Tbf$ action. Recall that from Definition \ref{defiTwistedEquivSheaves}, we have a notion of twisted equivariant sheaves. 

\begin{defi}
We define $\DDic(X, \RR_{\Tbf})_{\chi}$ to be the category of $(\Tbf, L_{\Tbf, \chi})$-equivariant sheaves on $X$. We denote by $\DDc(X, \RR_{\Tbf})_{\chi}$ the fullsubcategory of sheaves that are constructible on $X$. 
\end{defi}

We call this category the category of free monodromic sheaves on $X$. 

\begin{rque}
Recall from Section \ref{subsectionTwistedEquivSheaves}, the natural forgetful functor $\DDic(X, \RR_{\Tbf})_{\chi} \to \DDic(X, \RR_{\Tbf})$ has a left adjoint $\Av_{\chi}$ given by $a_!(L_{\Tbf, \chi}[2\dim \Tbf] \boxtimes - ): \DDic(X, \RR_{\Tbf}) \to \DDic(X, \RR_{\Tbf})_{\chi}$ where $a : \Tbf \times X \to X$ is the action map. 
\end{rque}

\begin{lem}
For $f : X \rightarrow Y$ a morphism of schemes with $\Tbf$-actions and $\chi \in \Ch(\Tbf)$ there are well defined functors 
\begin{equation*}
f_!, f_* : \DDic(X, \RR_{\Tbf})_{\chi} \leftrightarrows \DDic(Y, \RR_{\Tbf})_{\chi} : f^*, f^!.
\end{equation*}
\end{lem}

\begin{proof}
This is immediate from the functoriality of the formation of $\Tbf$-invariants. 
\end{proof}

\begin{rque}
The category $\DDic(X, \RR_{\Tbf})_{\chi}$ is a priori not monoidal, or at least not in a way that is compatible with the forgetful functor $\DDic(X, \RR_{\Tbf})_{\chi} \rightarrow \DDic(X, \RR_{\Tbf})$. 
\end{rque}

\begin{lem}
The functor $X \mapsto \DDic(X, \RR_{\Tbf})_{\chi}$ satisfies $!$ and $*$-descent along $\Tbf$-equivariant smooth maps.
\end{lem}

\begin{proof}
Let $Y \rightarrow X$ be a $\Tbf$-equivariant smooth map. Since $\DDic$ satisfies $!$ and $*$ descent along smooth maps we have 
$$ \DDic(X, \RR_{\Tbf}) = \varprojlim_{n, ?} \DDic(Y^{\times_X n}, \RR_{\Tbf})$$ where $? \in \{!,*\}$ and the transition maps in the limit is taken with respect to $?$-pullback. Since $\DDic(X, \RR_{\Tbf})_{\chi}$ is defined in term of categorical invariants, its formation commutes with limits. Hence we have
$$ \DDic(X, \RR_{\Tbf})_{\chi} = \varprojlim_{n,?} \DDic(Y^{\times_X n}, \RR_{\Tbf})_{\chi},$$ 
which is the desired statement. 
\end{proof}

Let $\chi_{\ell} \in \Ch_{\Qlb}(\Tbf)$ be a character of order a power of $\ell$. Consider the map $\iota_{\chi} : \RR_{\Tbf, \Zlb} \rightarrow \RR_{\Tbf, \Qlb}$. Then as $L_{\Tbf, \Zlb} \otimes_{\RR_{\Tbf, \Zlb}, \iota_{\chi_{\ell}}} \RR_{\Tbf, \Qlb} = L_{\Tbf, \Qlb} \otimes_{\Qlb} \Lcal_{\chi_{\ell}}$, the functor $-\otimes_{\RR_{\Tbf, \Zlb}, \iota_{\chi_{\ell}}} \RR_{\Tbf, \Qlb}$ induces a well defined functor 
\begin{equation*}
\DDic(X, \RR_{\Tbf, \Zlb})_{\chi} \rightarrow \DDic(X, \RR_{\Tbf, \Qlb})_{\chi\chi_{\ell}}
\end{equation*}
for all $\chi \in \Ch_{\Zlb}(\Tbf)$.

\subsection{Comparison with Verdier's monodromy action}

Let $X$ be a scheme with a $\Tbf$-action and let $\Lambda \in \{\Flb, \Zlb, \Qlb, \RR_{\Tbf, \Flb}, \RR_{\Tbf, \Zlb}, \RR_{\Tbf, \Qlb}\}$ be a coefficient ring. 

\begin{defi}
Let $Y$ be a connected scheme and $\bar{y}$ a geometric point of $Y$. We say $Y$ is a categorical $K(\pi,1)$ if the realization functor 
\begin{equation*}
\DD^b(\Rep_{\Lambda}(\pi_1(Y, \bar{y}))) \rightarrow \DD_{\lis}(Y, \Lambda),
\end{equation*}
is an equivalence, where $\Rep_{\Lambda}(\pi_1(Y, \bar{y}))$ is the category of continuous representations of $\pi_1(Y,\bar{y})$ on finite dimension $\Lambda$-modules.
\end{defi}

\begin{lem}[\cite{Achinger}]\label{lemCatPi1}
The torus $\Tbf$ is a categorical $K(\pi,1)$. 
\end{lem}

Let $A \in \Rep_{\Lambda}(\pi_1(\Tbf))$. Recall that $\pi_1(\Tbf)^{\wild}$ is normal a normal subgroup in $\pi_1(\Tbf)$. We denote by $\DD_{\lis}^{\tame}(\Tbf, \Lambda)$ the full subcategory of $\DD_{\lis}(\Tbf, \Lambda)$ of sheaves $A$ for which $\pi_1(\Tbf)^{\wild}$ acts trivially on $H^i(A)$. 

\begin{lem}\label{lemSplit}
The functor 
$$\DD^b(\Rep_{\Lambda}(\pi_1^t(\Tbf))) \to \DD_{\lis}(\Tbf, \Lambda)$$
is fully faithful with essential image $\DD_{\lis}^{\tame}(\Tbf, \Lambda)$.
\end{lem}

\begin{proof}[Proof of Lemma \ref{lemSplit}]
By Nakayama, we may assume that $\Lambda$ is a field. The source category is then generated by its irreducible representations. It is therefore enough to prove that for $\chi_1, \chi_2$ two irreducible representations, we have 
$$\Ext^i(\chi_1, \chi_2) = \Hom^i_{\DD_{\lis}(\Tbf,\Lambda)}(\chi_1, \chi_2).$$
We may therefore assume that $\chi_1 = \chi_2 = \Lambda$. In this case, both sides identify with $H^*(\Tbf, \Lambda)$. 
\end{proof}

\begin{defi}
Let $\tilde{\Tbf}^t = \varprojlim_{n} \Tbf$ the inverse limits of copies of $\Tbf$ where the transition maps $[n] : \Tbf \to \Tbf$ are $t \mapsto t^n$ for $n$ prime to $p$. This is a group scheme over $k$. 
\end{defi}

The map $\tilde{\Tbf}^t \to \Tbf$ is a morphism of group schemes with kernel equal to $\pi_1^t(\Tbf)$.  

\begin{rque}
We note that the multiplication map $\tilde{m} : \tilde{\Tbf}^t \times \tilde{\Tbf}^t \to \tilde{\Tbf}^t$ is not of finite type (since $\tilde{\Tbf}^t$ is not of finite type) but it is of finite expansion, see Remark \ref{rqueFiniteExpansion}. In particular $\tilde{m}_!$ is well defined and fits into the $6$-functors formalism given by $\DDic(-, \Lambda)$. Therefore for $X$ a scheme with a $\Tbf$-action, we can consider the action of $\tilde{\Tbf}^t$ on $X$ induced by the natural projection $\tilde{\Tbf}^t \to \Tbf$ and the category of $\tilde{\Tbf}^t$-equivariant sheaves on $X$ is well defined. 
\end{rque} 

\begin{lem}[\cite{Verdier}]\label{lemVerdierIdempotency}
We have $\RGamma_c(\tilde{\Tbf}^t, \Lambda) = \Lambda[-2\dim \Tbf]$.
\end{lem}

\begin{thm}\label{thmDefiMonodromic}
Let $X$ be $\Tbf$-scheme, the forgetful functor $\DDic(X/\tilde{\Tbf}^t, \Lambda) \to \DDic(X, \Lambda)$ is fully faithful and its image is the full subcategory of sheaves $A$ on $X$ satisfying one of the following equivalent conditions
\begin{enumerate}
\item for all $x \in X$, denote by $a_x : \Tbf \to X, t \mapsto t.x$ the orbit map of $x$, the sheaf $a_x^*A$ is ind-lisse and tame on $\Tbf$,
\item for all $x \in X$, the sheaf $a_x^! A$ is ind-lisse and tame on $\Tbf$. 
\end{enumerate}
These sheaves are call $\Lambda$-monodromic sheaves. 
\end{thm}

\begin{rque}
This notion was first introduced by Verdier \cite{Verdier} and a weak version of the $\tilde{\Tbf}^t$-equivariance was discussed. The use of pro-étale cohomology allows us to make sense of this equivariance condition. 
\end{rque} 

\begin{rque}
Note that since $\pi_1^t(\Tbf) = \ker(\tilde{\Tbf}^t \to \Tbf)$ acts trivially on $X$ any monodromic sheaf admits an action of $\pi_1^t(\Tbf)$ called the canonical monodromy. 
\end{rque}

\begin{proof}
By Remark \ref{rqueEquivAsModules}, the category of $\tilde{\Tbf}^t$-equivariant sheaves on $X$ is the category of modules over the algebra $\Lambda_{\tilde{\Tbf}^t} \in \DDic(\tilde{\Tbf}^t, \Lambda)$. By Lemma \ref{lemVerdierIdempotency}, this algebra is idempotent and the forgetful functor is thus fully faithful. The two conditions in the theorem are equivalent by Verdier duality. It remains to identify the subcategory of $\tilde{\Tbf}^t$-equivariant sheaves with the sheaves satisfying the category of sheaves satisfying these conditions. Since the orbit maps are $\Tbf$-equivariant for the translation action of $\Tbf$ on itself and that any $\tilde{\Tbf}^t$-equivariant sheaf on $\Tbf$ is lisse and tame, any sheaf in the essential image of the forgetful functor satisfy condition $(i)$ of the theorem. Conversely since averaging by $\tilde{\Tbf}^t$ is an idempotent functor, it is enough to show that for any $A \in \DDic(X, \Lambda)$ satisfying condition $(i)$, the natural map $A \to \Av_{\tilde{\Tbf}^t}(A)$ is an isomorphism. This last claim can be checked after pulling back to orbits where the claim is now clear by condition $(i)$. 
\end{proof}

\begin{defi}
Let $\chi \in \Ch_{\Lambda}(\Tbf)$ be a multiplicative sheaf on $\Tbf$, we say that a sheaf $A \in \DDic(X, \Lambda)$ is $\chi$-monodromic if is in the subcategory of $\DDic(X, \Lambda)$ generated by $\chi$-equivariant sheaves. For $\chi = 1$ we call them unipotent monodromic sheaves and denote the corresponding category $\DDic(X, \Lambda)_{\unip-\mon}$. 
\end{defi}

\begin{rque}
The condition $\chi$-monodromic is local in the $\Tbf$-equivariant smooth topology. 
\end{rque} 

Consider now the full subcategory $\DDc(X, \RR_{\Tbf})^{\Lambda-\cons}$ of $\DDc(X, \RR_{\Tbf})$ of objects such that their image under the forgetful functor $\DDc(X, \RR_{\Tbf}) \rightarrow \DD(X_{\proet}, \Lambda)$ is in $\DDc(X, \Lambda)$. This then induces a well defined functor $\For : \DDc(X, \RR_{\Tbf})^{\Lambda-\cons}_{\unip} \rightarrow \DDc(X, \Lambda)$.

\begin{lem}
The forgetful functor induces an equivalence $\DDc(X, \RR_{\Tbf})^{\Lambda-\cons}_{\unip} \rightarrow \DDc(X, \Lambda)_{\mon, \unip}$.
\end{lem}

\begin{proof}
Consider the Bar resolution $X \times \Tbf^{\bullet+1} \rightarrow X$ of $X$, since both sides satisfy smooth descent, we have a commutative diagram 
\[\begin{tikzcd}
	{\DDc(X, \RR_{\Tbf})_{\unip}^{\Lambda-\cons} } & {\varprojlim_{\Delta}\DDc(X \times \Tbf^{\bullet+1}, \RR_{\Tbf})_{\unip}^{\Lambda-\cons} } \\
	{\DDc(X, \Lambda)_{\mon,\unip} } & {\varprojlim_{\Delta}\DDc(X \times \Tbf^{\bullet+1}, \Lambda)_{\mon, \unip} }
	\arrow[Rightarrow, no head, from=1-1, to=1-2]
	\arrow[Rightarrow, no head, from=2-1, to=2-2]
	\arrow[from=1-2, to=2-2]
	\arrow[from=1-1, to=2-1]
\end{tikzcd}\]
Hence, it is enough to show the statement for $X \times \Tbf^{n+1}$. More generally assume that $X = Y \times \Tbf$ splits $\Tbf$-equivariantly as a product where $\Tbf$ acts trivially on $Y$. 

We first show that the forgetful functor is fully faithful. Let $A' = A \boxtimes_{\Lambda} \Lambda_{\Tbf} \in \DDc(Y \times \Tbf, \Lambda)_{\mon, \unip}$, since $\Lambda_{\Tbf} = L_{\Tbf}/I$,  where $I$ denotes the augmentation ideal of $\RR_{\Tbf}$, we have  $A' = \For(A_0)$, where $A_0 = (A \otimes_{\Lambda} \RR_{\Tbf}) \boxtimes_{\RR_{\Tbf}} L_{\Tbf}/I$. Let $A' = A \boxtimes_{\Lambda} \Lambda_{\Tbf}, B' = B \boxtimes_{\Lambda} \Lambda_{\Tbf}$ be two such objects and denote by $A_0$ and $B_0 \in \DDc(X, \RR_{\Tbf})_{\unip}^{\Lambda-\cons}$ the corresponding lifts, then the Künneth formula implies that 
\begin{equation*}
\Hom(A_0, B_0) = \Hom_{\Lambda}(A,B) \otimes_{\Lambda} \Hom_{\DDc(\Tbf, \RR_{\Tbf})_{\unip}}(\Lambda_{\Tbf}, \Lambda_{\Tbf}). 
\end{equation*}
Let us evaluate $\Hom_{\DDc(\Tbf, \RR_{\Tbf})_{\unip}}(\Lambda_{\Tbf}, \Lambda_{\Tbf})$. We have 
\begin{align*}
\Hom_{\DDc(\Tbf, \RR_{\Tbf})_{\unip}}(\Lambda_{\Tbf}, \Lambda_{\Tbf}) &= \Hom_{\RR_{\Tbf}}(\RR_{\Tbf}/I, \RR_{\Tbf}) \otimes_{\RR_{\Tbf}} \Hom_{\DDc(\Tbf, \RR_{\Tbf})_{\unip}}(L_{\Tbf}, L_{\Tbf}) \otimes_{\RR_{\Tbf}} \RR_{\Tbf}/I \\
&= \Hom_{\RR_{\Tbf}}(\RR_{\Tbf}/I, \RR_{\Tbf}) \otimes_{\RR_{\Tbf}} \RR_{\Tbf} \otimes_{\RR_{\Tbf}} \RR_{\Tbf}/I \\
&= \End_{\RR_{\Tbf}}(\RR_{\Tbf}/I) \\
&= \RGamma(\Tbf, \Lambda).
\end{align*}
The second line comes from the fact that $\Hom_{\DDc(\Tbf, \RR_{\Tbf})_{\unip}}(L_{\Tbf}, L_{\Tbf}) = \RR_{\Tbf}$ which can be seen through the equivalence $\DDc(T, \RR_{\Tbf})_{\unip} = \DD^b_{\coh}(\RR_{\Tbf})$. The last line comes from Lemma \ref{lemCohomTorus}. On the other hand after applying the forgetful functor, we get that 
\begin{equation*}
\Hom_{\DDc(Y \times \Tbf, \Lambda)_{\mon-\unip}}(A',B') = \Hom_{\Lambda}(A, B) \otimes_{\Lambda} \RGamma(\Tbf, \Lambda),
\end{equation*}
as $\RGamma(\Tbf, \Lambda) = \End(\Lambda_ {\Tbf})$. Therefore the forgetful functor is fully faithful on objects of the form $A_0$. We show that these objects generate the category $\DDc(Y \times \Tbf, \RR_{\Tbf})_{\unip}^{\Lambda-\cons}$. Let $A_0 \in \DDc(Y \times \Tbf, \RR_{\Tbf})_{\unip}^{\Lambda-\cons}$, as $A_0$ is $(\Tbf, L_{\Tbf})$-equivariant $A_0$ can be written as $A_0 = A' \boxtimes L_{\Tbf}$, where $A'$ is an $\RR_{\Tbf}$-constructible sheaf on $Y$. But as $A_0$ is also $\Lambda$-constructible, the $\RR_{\Tbf}$-structure of $A'$ factors is of $I^{\infty}$-torsion. As such $A'$ is in the full subcategory of $\DDc(Y, \RR_{\Tbf})$ generated by the essential image of the functor $\DDc(Y, \Lambda) \rightarrow \DDc(Y, \RR_{\Tbf})$ induced by the forgetful functor along the augmentation $\RR_{\Tbf} \rightarrow \Lambda$. Let $A' \in \DDc(Y, \RR_{\Tbf})$ be in the essential image of $ \DDc(Y, \Lambda)$, then we have 
\begin{equation*}
A' \boxtimes_{\RR_{\Tbf}} L_{\Tbf} = A' \boxtimes_{\Lambda} L_{\Tbf}/I = A' \boxtimes_{\Lambda} \Lambda_{\Tbf}.
\end{equation*}
Let $C,D \in \DDc(X, \RR_{\Tbf})^{\Lambda-\cons}_{\unip}$, since objects of the form $A' \boxtimes \Lambda_{\Tbf}$ generate the category we can write $C = \varinjlim_i (A'_i \boxtimes \Lambda_{\Tbf})$ and $D = \varinjlim_j (B'_j \boxtimes \Lambda_{\Tbf})$ where both colimits are finite. As the forgetful functor is a right adjoint it commutes with limits hence it also commutes with finite colimits \cite[Proposition 1.1.4.1]{LurieB}. We then have
\begin{align*}
\Hom(C,D) &= \varprojlim_i \varinjlim_j \Hom((A'_i \boxtimes \Lambda_{\Tbf}), (B'_j \boxtimes \Lambda_{\Tbf})) \\
&= \varprojlim_i \varinjlim_j \Hom(\For(A'_i \boxtimes \Lambda_{\Tbf}), \For(B'_j \boxtimes \Lambda_{\Tbf})) \\
&= \Hom(\For(A), \For(B)).
\end{align*}
Hence $\For$ is fully faithful. 
Since the objects of the form $\For(A' \boxtimes_{\Lambda} \Lambda_{\Tbf})$ generate the category $\DDc(Y \times \Tbf, \Lambda)_{\mon, \unip}$ under finite colimits, the essential surjectivity is clear. Indeed let $A \in \DDc(Y \times \Tbf, \Lambda)_{\mon, \unip}$, we can then write $A = \varinjlim_i \For(A' \boxtimes_{\Lambda} \Lambda_{\Tbf}) = \For(\varinjlim_i A' \boxtimes_{\Lambda} \Lambda_{\Tbf})$. 
\end{proof}

\begin{lem}\label{lemCohomTorus}
There is a canonical isomorphism $\End_{\RR_{\Tbf}}(\RR_{\Tbf}/I) = \RGamma(\Tbf, \Lambda)$. 
\end{lem}

\begin{proof}
Consider the following functors $\DDc(\pt, \RR_{\Tbf})^{\Lambda-\cons} \xrightarrow{\alpha^*} \DDc(\Tbf, \RR_{\Tbf})^{\Lambda-\cons} \xrightarrow{\For} \DDc(\Tbf, \Lambda)$, where $\alpha : \Tbf \rightarrow \pt$ is the structure map. By functoriality we get a map 
\begin{equation*}
\End_{\RR_{\Tbf}}(\RR_{\Tbf}/I) \rightarrow \End_{\Tbf}(\Lambda_{\Tbf}) = \RGamma(\Tbf, \Lambda).
\end{equation*}
It remains to check that this is an isomorphism. This can be done after taking cohomology, namely, we want that the induced map 
\begin{equation*}
\Ext^*_{\RR_{\Tbf}}(\RR_{\Tbf}/I, \RR_{\Tbf}/I) \rightarrow H^*(\Tbf, \Lambda),
\end{equation*}
is an isomorphism. Since it is deduced from functoriality, this map is a map of algebras. It is known that both sides are exterior algebras on their degree one parts. For the left hand this is $(I/I^2)^{\vee}$ (where the $(-)^{\vee}$ is the $\Lambda$-linear dual) and for the right hand side this is $H^1(\Tbf, \Lambda)$. But those two are canonically isomorphic to $\Hom^0(\pi_1^{\ell}(\Tbf), \Zl) \otimes_{\Zl} \Lambda$. 
\end{proof}

Let $A \in \DDc(X, \RR_{\Tbf})_{\unip}^{\Lambda-\cons}$ and consider the object $A' \in \DDc(X, \Lambda)_{\unip, \mon}$ be the image of $A$ under the forgetful functor. Consider the $\Lambda[\pi_1^t(\Tbf)]$-module structure on $A'$, since $A'$ is unipotent monodromic the morphism $\Lambda[\pi_1^t(\Tbf)] \rightarrow \End(A')$ factors through $\RR_{\Tbf}$.

\begin{lem}\label{lemMono}
The two $\RR_{\Tbf}$-structures on $A'$, one coming from Verdier's monodromy and one coming from the forgetful functor, coincide.
\end{lem}

\begin{proof}
The object $A \in \DDc(X, \RR_{\Tbf})^{\Lambda-\cons}$ is an $\RR_{\Tbf}$-unipotent monodromic, its canonical monodromy is a morphism $\RR_{\Tbf}[\pi_1^t(\Tbf)] \rightarrow \End(A)$. As $A$ comes from an equivariant sheaf this morphism factors through $\RR_{\Tbf}[\pi_1^t(\Tbf)]/I$ where $I$ is the ideal generated by elements $(t - \can(t))$ for $t \in \pi_1^t(\Tbf)$ and $\can$ is the canonical map (\ref{canon}). 
Consider now the following diagram 
\[\begin{tikzcd}
	& {\Lambda[\pi_1^t(\Tbf)]} \\
	{\RR_{\Tbf}} & {\RR_{\Tbf}[\pi_1^t(\Tbf)]} \\
	& {\RR_{\Tbf}[\pi_1^t(\Tbf)]/I} \\
	& {\End(A')}
	\arrow[from=2-1, to=2-2]
	\arrow[from=1-2, to=2-2]
	\arrow[from=2-2, to=3-2]
	\arrow[from=1-2, to=2-1]
	\arrow[from=3-2, to=4-2]
\end{tikzcd}\]
where the map $\Lambda[\pi_1^t(\Tbf)] \rightarrow \RR_{\Tbf}$ is induced by the morphism $\can$ and the other morphisms are the natural ones. The triangle does not commute but it commutes after projecting in $\RR_{\Tbf}[\pi_1^t(\Tbf)]/I$. The canonical monodromy is the $\RR_{\Tbf}$-structure coming from the vertical composition while the $\RR_{\Tbf}$-structure on the forgetful functor is the map $\RR_{\Tbf} \rightarrow \End(A')$.
\end{proof}

We now denote by $\Phi_{\unip} : \DDc(X, \Lambda)_{\unip, \mon} \rightarrow \DDc(X, \RR_{\Tbf})_{\unip}$ the inverse of the forgetful functor. 

\begin{rque}
The previous construction naturally extends to the non unipotent setting, for $\chi \in \Ch(T)$ we get a fully faithful functor
\begin{equation*}
\Phi_{\chi} : \DDc(X, \Lambda)_{\mon, \chi} \rightarrow \DDc(X, \RR_{\Tbf})_{\chi}.
\end{equation*}
\end{rque}

\subsection{Unipotent monodromic sheaves}

Let $X$ be a scheme with a free $\Tbf$ action, let $Y = X/\Tbf$ and $\pi : X \rightarrow Y$ be the projection.

\begin{lem}\label{lemPushForwardUnipMono}
There exists a functorial isomorphism $\pi_! = \pi_*[-\dim \Tbf] : \DDc(X, \RR_{\Tbf})_{\unip} \rightarrow \DDc(Y, \RR_{\Tbf})$.
\end{lem}

\begin{proof}
First assume we are given a section of the $\Tbf$-torsor $X \to Y$ so that we have a splitting $X = Y \times \Tbf$. Then for all sheaves $A$ there is a canonical isomorphism $A = A_0 \boxtimes_{\RR_{\Tbf}} L_{\Tbf}$ where $A_0 = 1*A$ where $1 : Y \to Y \times \Tbf$ is induced by the unit of $\Tbf$. Then the isomorphism of functor comes from the isomorphism $\RGamma_c(\Tbf,L_{\Tbf}) = \RGamma(\Tbf,L_{\Tbf})[-\dim \Tbf]$ from \cite{GabberLoeser}. In general, let $X \times \Tbf^{\bullet + 1}$ be the Bar resolution of $X$, which is also the Cech nerve of $X \to Y$. By descent, we have 
$$\DDc(X, \RR_{\Tbf})_{\unip} = \varprojlim_{n} \DDc(X \times \Tbf^{\bullet + 1}, \RR_{\Tbf})_{\unip}$$
where the action of $\Tbf$ on $X \times \Tbf^{\bullet + 1}$ is the action on the last copy of $T$, hence $X \times \Tbf^{\bullet + 1}$ is equipped with a canonical splitting of the $\Tbf$-action. Taking pushforward yields the lemma. 
\end{proof}

\begin{rque}
The isomorphism $\RGamma_c(\Tbf,L_{\Tbf}) = \RGamma(\Tbf,L_{\Tbf})[-\dim \Tbf]$ is not canonical up to a $\Lambda$-module of rank one which is trivialized after choosing a basis of $X^*(\Tbf)$. 
\end{rque} 

\subsection{Comparison with Yun's definition of free monodromic sheaves}

We now give a comparison between our categories of monodromic sheaves and the completion of the categories of monodromic sheaves of \cite{BezYun}, \cite{BezrukRiche} and \cite{Gouttard}. In this section, let $Y$ be a scheme and let $X \rightarrow Y$ be a $\Tbf$-torsor. 

\begin{thm}\label{ThmEquivalenceBR}
We have a natural equivalence 
\begin{equation*}
\ho(\DDc(X, \RR_{\Tbf, \Fl})_{\unip}) \simeq \widehat{\mathrm{D}}^b_c(X \fatslash \Tbf)
\end{equation*}
where the category on the right is the completed monodromic category of \cite{BezrukRiche}. In the non-unipotent case there is an equivalence
\begin{equation*}
\ho(\DDc(X, \RR_{\Tbf, \F_{\ell^n}})_{\chi}) \simeq \widehat{\mathrm{D}}^b_c(X \fatslash \Tbf)_{\mathcal{L}_{\chi}},
\end{equation*}
which holds after passing from $\Fl$ to a finite extension $\F_{\ell^n}$ where $\chi$ is defined. 
\end{thm}

\begin{proof}
We only show the version for $\chi = 1$, the generalization to other $\chi$ is straightforward. To define the desired functor, first recall that the category $\widehat{\mathrm{D}}^b_c(X \fatslash \Tbf)$ is a full subcategory of the category $\Pro (\mathrm{D}^b_c(X \fatslash \Tbf))$ of monodromic objects on $X$, see \cite[3.1 and 10.1]{BezrukRiche}. We first define a functor 
\begin{equation*}
\Psi : \ho(\DDc(X, \RR_{\Tbf, \Fl})_{\unip}) \rightarrow \Pro (\mathrm{D}^b_c(X \fatslash \Tbf))
\end{equation*}
by $A \mapsto ``\varprojlim" A \otimes_{\RR_{\Tbf, \Fl}} \RR_{\Tbf, \Fl}/\mathfrak{m}^n$ where $\mathfrak{m}$ is the maximal ideal of $\RR_{\Tbf, \Fl}$. The ring $\RR_{\Tbf, \Fl}/\mathfrak{m}^n$ is an Artin ring over $\Fl$ and thus is finite dimensional and $A \otimes_{\RR_{\Tbf, \Fl}} \RR_{\Tbf, \Fl}/\mathfrak{m}^n$ lives in $\DDc(X, \Fl)$ after forgetting the $\RR_{\Tbf, \Fl}/\mathfrak{m}^n$-structure. For any $y \in Y$ the restriction to the fiber $X_y = \pi^{-1}(y)$ of $A$ is isomorphic to $M \otimes_{\RR_{\Tbf,\Fl}} L_{\Tbf}$ for some $\RR_{\Tbf, \Fl}$-module $M$ and therefore the restriction of $A \otimes_{\RR_{\Tbf, \Fl}} \RR_{\Tbf, \Fl}/\mathfrak{m}^n$ is isomorphic to the sheaf denoted by $M \otimes_{\RR_{\Tbf,\Fl}} \mathcal{L}_{T,n}$ in \cite[3.2]{BezrukRiche}  and in particular is monodromic on $\Tbf$. This implies that $A \otimes_{\RR_{\Tbf, \Fl}} \RR_{\Tbf, \Fl}/\mathfrak{m}^n$ is monodromic as an $\Fl$-sheaf on $X$ and that $\Psi$ is well defined.

We prove that it factors through the category $\widehat{\mathrm{D}}^b_c(X \fatslash \Tbf)$, which means checking the two properties of \cite[Definition 3.1.]{BezrukRiche}. The pro-object $\Psi(A)$  is $\pi$-constant, indeed we have the following computation 
\begin{align*}
\varprojlim \pi_! (A \otimes_{\RR_{\Tbf, \Fl}} \RR_{\Tbf,\Fl}/\mathfrak{m}^n) &= \varprojlim (\pi_! A) \otimes_{\RR_{\Tbf, \Fl}} \RR_{\Tbf,\Fl}/\mathfrak{m}^n \\
&= \pi_! A.
\end{align*}
The first line follows from the compatibility of the formation of $\pi_!$ with change of coefficients and the second comes from the fact that all constructible $\RR_{\Tbf,\Fl}$-sheaves on $X$ are derived complete. The object $\pi_! A$ a priori lives in $\DDc(Y, \RR_{\Tbf})$, but fiberwise it is isomorphic to $M \otimes_{\RR_{\Tbf, \Fl}} \RGamma_c(\Tbf, L_{\Tbf, \Fl}) \simeq M \otimes_{\RR_{\Tbf, \Fl}} \Fl[2\dim \Tbf]$ and therefore its stalks are perfect $\Fl$-complexes hence after forgetting the $\RR_{\Tbf,\Fl}$-structure down to an $\Fl$-structure, the sheaf $\pi_!A$ is constructible. 

The pro-object $\Psi(A)$ is also uniformly bounded in degrees. There exists $a \leq b$ two integers such that $F \in \DDc(X)^{[a,b]}$, the functor $-\otimes_{\RR_{\Tbf, \Fl}} \RR_{\Tbf,\Fl}/\mathfrak{m}^n$ is of cohomological dimension $[-\dim \Tbf, 0]$ and the forgetful functor $\DDc(X, \RR_{\Tbf,\Fl}/\mathfrak{m}^n) \rightarrow \DDc(X, \Fl)$ is $t$-exact. This implies that $A \otimes_{\RR_{\Tbf, \Fl}} \RR_{\Tbf,\Fl}/\mathfrak{m}^n$ lives in cohomological degree $[a - \dim \Tbf, b]$ and the functor $\Psi$ factors through $\widehat{\mathrm{D}}^b_c(X \fatslash \Tbf)$.

The functor $\Psi$ is fully faithful. Let $A,B \in \DDc(X,\RR_{\Tbf,\Fl})_{\unip}$, we have
\begin{align*}
\Hom(A,B) &= \varprojlim_n \Hom(A, B \otimes_{\RR_{\Tbf,\Fl}} \RR_{\Tbf,\Fl}/\mathfrak{m}^n) \\
&= \varprojlim_n \varinjlim_m \Hom(A \otimes_{\RR_{\Tbf,\Fl}} \RR_{\Tbf,\Fl}/\mathfrak{m}^m , B \otimes_{\RR_{\Tbf,\Fl}} \RR_{\Tbf,\Fl}/\mathfrak{m}^n).
\end{align*}
The first equality comes from the isomorphism $B = \varprojlim_n B \otimes_{\RR_{\Tbf,\Fl}} \RR_{\Tbf,\Fl}/\mathfrak{m}^n$ and the second from the same isomorphism for $A$ and the fact that each $B \otimes_{\RR_{\Tbf,\Fl}} \RR_{\Tbf,\Fl}/\mathfrak{m}^n$ is discrete and thus a morphism from $A$ factors through one of its quotients. We apply $H^0$ to this isomorphism, there is a Milnor short exact sequence 
\begin{align*}
0 \rightarrow \mathrm{R}^1\varprojlim H^{-1}(\Hom(A,B \otimes_{\RR_{\Tbf,\Fl}} \RR_{\Tbf,\Fl}/\mathfrak{m}^n) &	\rightarrow H^0\varprojlim_{n} \Hom(A, B \otimes_{\RR_{\Tbf,\Fl}} \RR_{\Tbf,\Fl}/\mathfrak{m}^n)) \\
& \rightarrow \varprojlim_{n} H^0(\Hom(A, B \otimes_{\RR_{\Tbf,\Fl}} \RR_{\Tbf,\Fl}/\mathfrak{m}^n))) \rightarrow 0.
\end{align*} 
Note that, as a complex, $\Hom(A,B) \in \DD(\RR_{\Tbf, \Fl})$ is perfect and thus derived complete, indeed the category of derived complete objects is stable and contains $\RR_{\Tbf, \Fl}$ hence all perfect complexes. By \cite[Tag 091P]{Stacks} , all the cohomlogy groups of $\Hom(A,B)$ are derived complete hence $H^{-1}(\Hom(A,B))$ is derived complete. Since it is an $\RR_{\Tbf, \Fl}$-module of finite type, by Nakayama's lemma it is also $\mathfrak{m}$-adically separated and therefore $\mathfrak{m}$-adically complete by \cite[Tag 091T]{Stacks} hence $\varprojlim_n H^{-1}(\Hom(A,B))/\mathfrak{m}^n = H^{-1}(\Hom(A,B))$ and the first term of the above exact sequence vanishes. Hence $H^0$ commutes with the limit, since the colimit is filtered, it is exact and commutes with $H^0$.  The fully faithfulness now follows from the description of the morphisms in $\widehat{\mathrm{D}}^b_c(X \fatslash \Tbf)$ \cite[3.1]{BezrukRiche} and \cite[Section A.2]{BezYun} .

It remains to show that $\Psi$ is essentially surjective, note that we have a compatibility between free monodromic local systems as $\Psi(L_{\Tbf}) = \mathcal{L}_{\Tbf}$ where the second sheaf is the free monodromic local system of \cite[3.2]{BezrukRiche}. Let $A = "\varprojlim_n" A_n$ be an object in $\widehat{\mathrm{D}}^b_c(X \fatslash \Tbf)$, we can assume that for each $A_n$, Verdier's monodromy $\phi_{A_n} : \RR_{\Tbf,\Fl} \rightarrow \End(A_n)$-factors through $\RR_{\Tbf,\Fl}/\mathfrak{m}^n$. Consider now $\tilde{A} = \varprojlim_n \Phi_{\unip}(A_n) \in \DDc(X, \RR_{\Tbf,\Fl})_{\unip}$. By construction $\tilde{A}/\mathfrak{m}^n = \Phi_{\unip}(A_n)$ and forgetting the $\RR_{\Tbf,\Fl}/\mathfrak{m}^n$-structure yields back $A_n$ hence $\Psi(\tilde{A}) = A$ and $\Psi$ is essentially surjective.
\end{proof}

\begin{rque}
We have shown the comparison for modular coefficients, the same argument extends to the $\Qlb$-case. 
\end{rque}

\subsection{Computation of some examples}

Our formalism allows to define free monodromic sheaves on schemes or stacks with an action of $\Tbf$ where the action of $\Tbf$ is not free. Let us compute some example which will be relevant later. 

\begin{lem}[Case of the point]
There are natural equivalences
\begin{enumerate}
\item if $\chi \in \Ch_{\Lambda}(\Tbf)$ is non-trivial, then $\DDic(\pt, \RR_{\Tbf})_{\chi} = 0$.
\item if $\chi = 1$, then $\DDic(\pt, \RR_{\Tbf})_{\unip} = \DD(\Lambda)$.
\end{enumerate}
\end{lem}

\begin{proof}
Consider the pair of adjoint functors $\Av_{\chi, !} : \DDic(\pt, \RR_{\Tbf,\Lambda}) \leftrightarrows \DDic(\pt, \RR_{\Tbf,\Lambda})_{\chi} : \For$. By definition this pair of adjoint functors is monadic, hence the category $\DDic(\pt, \RR_{\Tbf,\Lambda})_{\chi}$ is equivalent to the category of algebras over the monad $\For\Av_{\chi, !}$. By proper base change, this functor is isomorphic to 
$$M \mapsto \RGamma_c(\Tbf, M \otimes_{\RR_{\Tbf, \Lambda}} L_{\Tbf} \otimes_{\Lambda} \Lcal_{\chi})[2\dim \Tbf] = M \otimes_{\RR_{\Tbf, \Lambda}} \RGamma_c(\Tbf, L_{\Tbf} \otimes_{\Lambda} \Lcal_{\chi})[2\dim \Tbf].$$
This last term is zero if $\chi$ is non trivial, hence if $\chi$ is non trivial then the monad itself is the zero monad and $\DDic(\pt, \RR_{\Tbf})_{\chi} = 0$. If $\chi$ is trivial, then as $\RGamma_c(\Tbf, L_{\Tbf}[2\dim \Tbf]) = \Lambda$, this monad is identified with the monad corresponding to the $\RR_{\Tbf, \Lambda}$-algebra $\Lambda$, where the map $\RR_{\Tbf, \Lambda} \rightarrow \Lambda$ is the augmentation. Hence $\DDic(\pt, \RR_{\Tbf})_{\unip} = \Lambda-\Mod(\DD(\RR_{\Tbf, \Lambda})) = \DD(\Lambda)$. 
\end{proof}

\begin{lem}[Two copies of $\Tbf$-acting on itself]\label{lemExampleTwoActions}
Let $f : \Tbf \rightarrow \Tbf$ be a smooth endomorphism of $\Tbf$. Consider the action of $\Tbf \times \Tbf$ on $\Tbf$ given by $(t_1,t_2).h = t_1f(t_2)h$, then for a pair $(\chi_1, \chi_2) \in \Ch_{\Lambda}(\Tbf \times \Tbf)$, we have 
\begin{enumerate}
\item $\DDic(\Tbf, \RR_{\Tbf \times \Tbf})_{\chi_1, \chi_2} = \DDic(\Tbf, \RR_{\Tbf})_{\chi_1} = \DDic(\Tbf, \RR_{\Tbf})_{\chi_2}$ if $f(\chi_2) = \chi_1$, where $\DDic(T, \RR_{\Tbf})_{\chi_1}$ denotes the category of $L_{\Tbf, \Lambda} \otimes_{\Lambda} \Lcal_{\chi_1}$-equivariant sheaves with respect to the action of the first copy of $T$ and similarly for $\DDic(\Tbf, \RR_{\Tbf})_{\chi_2}$ with respect to the action of the second copy of $\Tbf$,
\item $\DDic(\Tbf, \RR_{\Tbf \times \Tbf})_{\chi_1, \chi_2} = 0$ if $f(\chi_2) \neq \chi_1$.
\end{enumerate}
\end{lem}

\begin{proof}
Consider the map $m : \Tbf^2 \to \Tbf$ given by $(t_1, t_2) \mapsto t_1f(t_2)$. This map is equivariant for the natural translation of $\Tbf^2$ on itself and the action of $\Tbf^2$ described in the lemma. Since the map $m$ is smooth and surjective, the category $\DDic(\Tbf, \RR_{\Tbf \times \Tbf})_{\chi_1, \chi_2}$ is then equivalent to the category of modules in $\DDic(\Tbf^2, \RR_{\Tbf^2})_{\chi_1, \chi_2}$ over the monad $m^!m_!$. Since $\DDic(\Tbf, \RR_{\Tbf \times \Tbf})_{\chi_1, \chi_2} = \DD(\RR_{\Tbf^2})$, it is enough to compute the value of $m^!m_!$ on $L_{\Tbf^2, \chi_1 \times \chi_2}$. An easy computation yields
\begin{equation*}
m^!m_!L_{\Tbf^2, \chi_1 \times \chi_2} = L_{\Tbf^2, \chi_1\times \chi_2} \otimes_{\RR_{\Tbf^2}} \RR_{\Tbf} \otimes \RGamma_c(\Tbf, L_{\Tbf}[2\dim \Tbf] \otimes\Lcal_{\chi_1} \otimes \Lcal_{f(\chi_2^{-1})}). 
\end{equation*}
Hence this monad is non zero if and only if $\chi_1 = f(\chi_2)$. In this case, under the equivalence $\DDic(\Tbf, \RR_{\Tbf \times T})_{\chi_1, \chi_2} = \DD(\RR_{\Tbf^2})$, the monad is isomorphic to $\otimes_{\RR_{\Tbf^2}} \RR_{\Tbf}$ hence $\DDic(\Tbf, \RR_{\Tbf^2})_{[\chi_1, \chi_2]} \simeq \DD(\RR_{\Tbf})$. The claim about forgetting the left and right actions of $T^2$ is clear. 
\end{proof}

\begin{lem}[Equivariant and monodromic]\label{lemEquivariantvsMonodromic}
Assume that $\Tbf$ is the base change of a torus defined over $\mathbb{F}_q$, and denote by $\Frob : \Tbf \rightarrow \Tbf$ the corresponding $k$-linear Frobenius. Then there are equivalences 
\begin{equation*}
\DDic(\frac{\Tbf}{\Ad_{\Frob}\Tbf}, \Lambda) = \oplus_{\chi} \DDic(\frac{\Tbf}{\Ad_{\Frob}\Tbf}, \RR_{\Tbf})_{\chi},
\end{equation*}
where the sum ranges through all $\chi \in \Ch_{\Lambda}(\Tbf)$, for the right hand side, the twisted equivariant sheaves are taken with respect to the action of $\Tbf$ acting by left or right translations, finally the action $\Ad_{\Frob}$ of $\Tbf$ on itself is given by $t.x = t\Frob(t^{-1})x$.
\end{lem}

\begin{proof}
We proceed as in the previous lemma, consider the quotient map $\pi : \Tbf \to \frac{\Tbf}{\Ad_{\Frob}\Tbf} = \pt/\Tbf^F$. This map is equivariant for the natural translation action of $\Tbf$ on the source and on the target. Since it is smooth and surjective, we can write the category $\DDic(\frac{\Tbf}{\Ad_{\Frob}\Tbf}, \RR_{\Tbf})_{\chi}$ as the category of modules over the monad $\pi^!\pi_!$ in $\DDic(\Tbf, \RR_{\Tbf})_{\chi} =  \DD(\RR_{\Tbf})$. We compute its value on the sheaf $L_{\Tbf, \chi}$. Consider the diagram 
\[\begin{tikzcd}
	{\Tbf} & {\Tbf} \\
	\pt & {\frac{\Tbf}{\Ad_{\Frob}\Tbf}}
	\arrow["\alpha", from=2-1, to=2-2]
	\arrow["{\mathcal{L}}", from=1-1, to=1-2]
	\arrow["\pi"{description}, from=1-2, to=2-2]
	\arrow["\beta"{description}, from=1-1, to=2-1]
\end{tikzcd}\]
where $\beta : \Tbf \to \pt$ is the structure map, $\alpha : \pt \to \pt/\Tbf^F$ is the universal torsor and $\Lcal$ is the Lang map. 
 We have 
\begin{equation*}
\Lcal^*\pi^!\pi_!L_{\Tbf, \chi} = \beta^*\alpha^*\pi_!L_{\Tbf,\chi}[2\dim \Tbf] = \beta^*\RGamma_c(\Tbf, \Lcal^*L_{\Tbf,\chi}).
\end{equation*}
Now $\RGamma_c(\Tbf, \Lcal^*L_{\Tbf,\chi}) = 0$ if $\Lcal^*\Lcal_{\chi}$ is nontrival. If $\Lcal^*\Lcal_{\chi}$ is trivial then we have $\Lcal^*L_{\Tbf} = L_{\Tbf} \otimes_{\RR_{\Tbf}} \RR_{\Tbf}$ where the map $\RR_{\Tbf} \to \RR_{\Tbf}$ is induced by the Lang map. Taking cohomology yields
$$  \RGamma_c(\Tbf, \Lcal^*L_{\Tbf}) = \RR_{\Tbf} \otimes_{\RR_{\Tbf}} \otimes \Lambda[-2\dim \Tbf].$$
\begin{enumerate}
\item If $\Lambda \in \{\Zlb, \Flb\}$ then by Lemma \ref{lemFunctorialityOrderL}, this is isomorphic to $\Lambda[\Tbf^F[\ell^{\infty}]]$, 
\item if $\Lambda = \Qlb$, then by Lemma \ref{lemFunctorialityPrimeToL}, this is isomorphic to $\Qlb$. 
\end{enumerate}
In particular, all objects in $\DDic(\frac{\Tbf}{\Ad_{\Frob}\Tbf}, \RR_{\Tbf})_{\chi}$ are $\Lambda$-constructible and we get that $\DDic(\frac{\Tbf}{\Ad_{\Frob}\Tbf}, \RR_{\Tbf})_{\chi}$ is the category of $\RR_{\Tbf} \otimes_{\RR_{\Tbf}} \otimes \Lambda$-modules in $\DD(\Tbf, \RR_{\Tbf})_{\chi} = \DD(\RR_{\Tbf})$ which is then isomorphic to $\DD(\RR_{\Tbf} \otimes_{\RR_{\Tbf}} \otimes \Lambda)$. More canonically we get an equivalence
$$\DDic(\frac{\Tbf}{\Ad_{\Frob}\Tbf}, \RR_{\Tbf})_{\chi} = e_{\chi}\DD(\frac{\Tbf}{\Ad_{\Frob}\Tbf}, \Lambda)$$
where $e_{\chi}$ is the idempotent of $\Lambda[\Tbf^{\Frob}]$-projecting on the block corresponding to $\chi$. 
\end{proof}

\begin{rque}
The proof of Lemma \ref{lemEquivariantvsMonodromic} also shows the following. The Lang map $\mathcal{L}_{\Frob} : \Tbf \to \Tbf$ is a $\Tbf^{\Frob}$-cover and therefore defines a map $\pi_1(\Tbf) \to \Tbf^{\Frob}$. The summand $\DDic(\frac{\Tbf}{\Ad_{\Frob}\Tbf}, \RR_{\Tbf})_{\chi}$ is nonzero only if the character $\chi$ of $\pi_1(\Tbf)$ factors through $T^{\Frob}$. Hence, under the natural equivalence $\frac{\Tbf}{\Ad_{\Frob}\Tbf} = \pt/\Tbf^{\Frob}$, the decomposition 
$$\DDic(\frac{\Tbf}{\Ad_{\Frob}\Tbf}, \Lambda) = \oplus_{\chi} \DDic(\frac{\Tbf}{\Ad_{\Frob}\Tbf}, \RR_{\Tbf})_{\chi}$$
is exactly the block decomposition of $\DDic(\frac{\Tbf}{\Ad_{\Frob}\Tbf}, \Lambda)  = \DDic(\pt/\Tbf^{\Frob}, \Lambda)$.

\end{rque}

\subsection{Verdier duality for free monodromic sheaves}\label{sectionDuality}

In \cite{BezrukavnikovTomalchov}, the authors define a duality functor on completed unipotent monodromic categories extending the usual Verdier duality on $\Qlb$-constructible monodromic sheaves. We give a construction here that does not involve pro-objects, works for all schemes $X$ equipped with an action of $\Tbf$ and is valid in the non-unipotent setting. 

\begin{lem}\label{lemMonodDual}
Let $\Lambda$ be a coefficient ring. Let $X$ be a stack with a $\Tbf$ action and let $A$ be a $\Lambda$ constructible $\Tbf$-monodromic sheaf. The Verdier dual $\D_{\Lambda}(A)$ is monodromic and its canonical monodromy is given by 
\begin{equation*}
\phi^{\vee} : \Lambda[\pi_1^t(\Tbf)] \rightarrow H^0\End(A),
\end{equation*}
where $\phi$ is the canonical monodromy of $A$ and $\phi^{\vee} = \phi \circ \inv^*$ with $\inv : \Lambda[\pi_1^t(\Tbf)] \rightarrow \Lambda[\pi_1^t(\Tbf)]$ induced by $t \mapsto t^{-1}$. 
\end{lem}

\begin{proof}
We only need to check this on the fibers of $X \rightarrow X/\Tbf$ which are all isomorphic to $\Tbf$. This now follows from the fact that since $\Tbf$ is smooth, the Verdier dual of a lisse sheaf is lisse and corresponds to the dual representation of $\pi_1(\Tbf)$.
\end{proof}

\begin{defi}\label{defiTwist}
The map $\inv : \Tbf \rightarrow \Tbf$ induces a map $\inv_* : \RR_{\Tbf} \rightarrow \RR_{\Tbf}$. Given an $\RR_{\Tbf}$ module $M$, we denote by $M(\varepsilon) = M \otimes_{\RR_{\Tbf}, \inv_*} \RR_{\Tbf}$.
\end{defi}

\begin{rque}
Note that $L_{\Tbf}(\varepsilon) = L_{\Tbf}^{\vee}$ is the $\RR_{\Tbf}$-linear dual of $L_{\Tbf}$.
\end{rque}

\begin{defi}
Let $X$ be a stack with an action of $\Tbf$. We define 
\begin{align*}
\DDc(X, \RR_{\Tbf})_{\chi} &\rightarrow \DDc(X, \RR_{\Tbf})_{\chi^{-1}} \\
\D' &= \D_{\RR_{\Tbf}} (-) (\varepsilon),
\end{align*}
where $\D_{\RR_{\Tbf}}$ is the $\RR_{\Tbf}$-linear Verdier duality functor. 
\end{defi}

\begin{lem}
The functor $\D'$ satisfies
\begin{enumerate}
\item On the full subcategory $\DDc(X, \Lambda)_{\chi, \mon}$, we have a canonical isomorphism of functors 
\begin{equation*}
\D' = \D_{\Lambda}[-\dim \Tbf],
\end{equation*}
\item $\D'\D' = \id,$ 
\item For $A, B \in \DDc(X, \RR_{\Tbf})_{\chi}$ we have $\Hom(A,B) = \Hom(\D'(B), \D'(A))$. 
\item Given $f : X \rightarrow Y$ a morphism of $\Tbf$-scheme, we have $\D'f_! = f_*\D'$ and $\D'f^! = f^*\D'$. 
\end{enumerate}
\end{lem}

\begin{proof}
The last three points follow from the definition. We discuss the first one. We can work locally in the lisse topology and assume that we have a $\Tbf$-equivariant splitting $X = Y \times \Tbf$. We can then further assume that $Y$ is a point. Let $A$ be a $\chi$-monodromic sheaf on $T$, we can write $A = M \otimes_{\RR_{\Tbf}} (L_{\Tbf} \otimes \mathcal{L}_{\chi})[\dim \Tbf]$ for the $\RR_{\Tbf}$-module $M = 1^*[-\dim \Tbf]A$. Then by definition $\D'(A) = \Hom_{\RR_{\Tbf}}(M, \RR_{\Tbf}) \otimes_{\RR_{\Tbf}} (L_{\Tbf} \otimes \mathcal{L}_{\chi^{-1}})[\dim \Tbf]$. 
On the other hand, we have $1^*\D_{\Lambda}(A)[-\dim \Tbf] = \Hom_{\Lambda}(M, \Lambda)$, where $M$ is the $\Lambda$-module obtained by forgetting the $\RR_{\Tbf}$-stucture along the inclusion $\Lambda \rightarrow \RR_{\Tbf}$. 

We claim, that there is a natural $\RR_{\Tbf}$-linear isomorphism 
\begin{equation*}
\Hom_{\Lambda}(M, \Lambda) = \Hom_{\RR_{\Tbf}}(M, \RR_{\Tbf})(\varepsilon)[\dim \Tbf],
\end{equation*}
which is induced by local Serre duality for the pushforward along the map $\Spec(\RR_{\Tbf}) \rightarrow \Spec(\Lambda)$. 

Let us show this claim. Let $I \subset \RR_{\Tbf}$ be the augmentation ideal. Since $A$ is $\Lambda$-constructible, $M$ is of $I$-power torsion. Then we have 
\begin{align*}
\Hom_{\Lambda}(M, \Lambda) &= \Hom_{\RR_{\Tbf}}(M, \Hom_{\Lambda}(\RR_{\Tbf}, \Lambda)) \\
&= \Hom_{\RR_{\Tbf}}(M, \Gamma_I(\Hom_{\Lambda}(\RR_{\Tbf}, \Lambda)) \\
&= \Hom_{\RR_{\Tbf}}(M, \varinjlim_n(\Hom_{\RR_{\Tbf}}(\RR_{\Tbf}/I^n, \Hom_{\Lambda}(\RR_{\Tbf}, \Lambda)) \\
&= \Hom_{\RR_{\Tbf}}(M, \varinjlim_n(\Hom_{\Lambda}(\RR_{\Tbf}/I^n, \Lambda)) \\
&= \varinjlim_n\Hom_{\RR_{\Tbf}}(M, (\Hom_{\Lambda}(\RR_{\Tbf}/I^n, \Lambda)) \\
\end{align*}
where the first line comes from the adjunction between forgetful and $\Hom$, the second one from the fact that $M$ is of $I$-power torsion, the third one from the definition of local cohomology, the fourth one again from the adjunction and the last one from the compacity of $M$ as an $\RR_{\Tbf}$-module. 

On the other side, we have 
\begin{align*}
\Hom_{\RR_{\Tbf}}(M, \RR_{\Tbf}) &= \Hom_{\RR_{\Tbf}}(M, \Gamma_I(\RR_{\Tbf})) \\
&= \Hom_{\RR_{\Tbf}}(M, \varinjlim_n\Hom_{\RR_{\Tbf}}(\RR_{\Tbf}/I^n, \RR_{\Tbf})) \\
&= \varinjlim_n \Hom_{\RR_{\Tbf}}(M, \Hom_{\RR_{\Tbf}}(\RR_{\Tbf}/I^n, \RR_{\Tbf})).
\end{align*}

Let $\Tbf'$ be the split torus defined over $\Spec(\Z)$ equipped with an isomorphism $\Tbf'\times_{\Z} k = \Tbf$ and let $R = \mathcal{O}(\Tbf')$ and let $I' \subset R$ be the augmentation ideal. There is a natural flat map $R \rightarrow \RR_{\Tbf}$ and such that $I'\RR_{\Tbf} = I$. This implies in particular that $\Hom_{\RR_{\Tbf}}(\RR_{\Tbf}/I^n, \RR_{\Tbf}) = \RR_{\Tbf} \otimes \Hom_R(R/I'^n, R)$. Let $f : \Tbf' \rightarrow \Spec(\Z)$ be the structure map. Embedding the categories of $R$-modules and $\Z$-modules into the categories of solid $R$-modules and solid $\Z$-modules, built in \cite{Condensed1}. We get a pair of adjoint functors 
\begin{equation*}
f_! : \DD(R_{\blacksquare}) \leftrightarrows \DD(\Z_{\blacksquare}) : f^!.
\end{equation*}
Moreover by \cite[Observation 8.12]{Condensed1}, $f^!\Z = R[\dim \Tbf]$. A priori $f_!$ can be difficult to compute, but for the $R$-module $R/I'$, we have $f_!R/I'^n = f_*R/I'^n$ since $R/I'^n = i_*R/I'^n$ where $i : \Spec(R/I'^n) \rightarrow \Spec(R)$ is the closed embedding. Indeed, the formation of $i_!$ is compatible with composition and $i_! = i_*$ for proper maps,see \cite[ Theorem 11.1]{Condensed1} and the following discussion, but the map $\Spec(R/I'^n) \rightarrow \Spec(\Z)$ is finite hence proper. The adjunction and base change therefore provide a canonical $\RR_{\Tbf}$-linear isomorphism 
\begin{equation*}
\Hom_{\Lambda}(M, \Lambda) \rightarrow \Hom_{\RR_{\Tbf}}(M, \RR_{\Tbf})[\dim \Tbf].
\end{equation*}
\end{proof}

\begin{rque}
The setup of \cite{Condensed1} requires to consider rings that are of finite type over $\Z$ which is why we reduced everything to the ring $R$. 
\end{rque}

\begin{rque}\label{rmkRelativeTwists}
Let $\Tbf_1$ and $\Tbf_2$ be two tori and let $\Tbf_1 \to \Tbf_1 \times \Tbf_2 = \Tbf$ be the first inclusion. Denote by $\varepsilon_{\Tbf/\Tbf_1}$ the $\RR_{\Tbf}$-module 
$$\varepsilon_{\Tbf/\Tbf_1} = \RR_{\Tbf}(\varepsilon_{\Tbf}) \otimes_{\RR_{\Tbf_1}} \RR_{\Tbf_2}(\varepsilon_{\Tbf_{2}})$$
where the $\RR_{\Tbf_2}$-structure on $\RR_{\Tbf}$ is induced via the second inclusion. Then a variation of the previous argument yields an $\RR_{\Tbf}$-linear isomorphism of functors 
\begin{align*}
\DDc(X, \RR_{\Tbf})_{\chi}^{\op,\RR_{\Tbf_1}-c} &\to \DDc(X, \RR_{\Tbf_1})_{\chi^{-1}} \\
\For_{\RR_{\Tbf}}^{\RR_{\Tbf_1}} \D'_{\RR_{\Tbf}} &= \D_{\RR_{\Tbf_1}}'\For_{\RR_{\Tbf}}^{\RR_{\Tbf_1}}[\dim \Tbf_2](\varepsilon_{\Tbf/\Tbf_1})
\end{align*}
where $\For_{\RR_{\Tbf}}^{\RR_{\Tbf_1}}$ is the functor that forgets the $\RR_{\Tbf}$-structure down to an $\RR_{\Tbf_1}$-structure. 
\end{rque}

\section{Deligne-Lusztig theory and character sheaves}\label{sectionHorocycle}

We first review the formalism of Lusztig of horocycle correspondences, they were introduced in \cite{LusztigCS1} to define the character sheaves and various twisted version have been introduced since then. The ones we will need have been defined in \cite{LusztigCatCenter2}.

\subsection{Horocycle correspondences and their parabolic versions}

Let $\Gbf$ be a group scheme over $k$ such that $\Gbf^{\circ}$ is reductive. We let $\Frob : \Gbf \rightarrow \Gbf$ an endomorphism of $\Gbf$ that is either a purely inseparable isogeny such that a power of it is a Frobenius endomorphism (Frobenius case) or a finite order automorphism of $\Gbf$. We call morphisms of the first type Frobenius endomorphisms of $\Gbf$. We define the $\Frob$-twisted adjoint action of $\Gbf$ on itself as 
$$\Ad_{\Frob}(g)(x) = gx\Frob(g^{-1}).$$

Let $\Pbf = \Lbf\Vbf$ be a parabolic with Levi $\Lbf$ of $\Gbf^{\circ}$ and consider the following correspondence, usually called the horocycle correspondence
$$\frac{\Gbf}{\Ad_{\Frob} \Gbf} \xleftarrow{\mathfrak{r}_{\Pbf, \Frob}} \frac{\Gbf}{\Ad_{\Frob}\Gbf} \xrightarrow{\qfrak_{\Pbf, \Frob}} \frac{\Vbf\backslash \Gbf/\Frob(\Vbf)}{\Ad_{\Frob}\Lbf}.$$

The map $\qfrak_{\Pbf,\Frob}$ is a smooth $V$-fibration and the map $\rfrak_{\Pbf,\Frob}$ is $\Gbf/\Pbf$-fibration, in particular it is smooth an proper. We denote by 
\begin{enumerate}
\item $\hcfrak_{\Frob,\Pbf} = \qfrak_{\Pbf, \Frob, !}\rfrak_{\Frob,\Pbf}^*$, 
\item $\chfrak_{\Frob,\Pbf} = \rfrak_{\Frob, \Pbf, !}\rfrak_{\Frob,\Pbf}^*$.
\end{enumerate}
We have an adjunction
$$(\hcfrak_{\Frob,\Pbf}, \chfrak_{\Frob,\Pbf}[2\dim \Ubf]).$$

\subsection{The case of a Frobenius and Deligne--Lusztig functors}\label{sectionDLTheory}

We fix $\Frob : \Gbf \rightarrow \Gbf$ a Frobenius endomorphism of $\Gbf$. 

\begin{lem}
Assume that $\Gbf$ is connected then we have a canonical equivalence 
\begin{equation*}
\DDic(\frac{\Gbf}{\Ad_{\Frob}\Gbf}, \Lambda) = \DD(\Rep_{\Lambda}(\Gbf^{\Frob})). 
\end{equation*}
\end{lem}

\begin{proof}
Since $\Gbf$ is connected, by Lang's theorem we have $\frac{\Gbf}{\Ad_{\Frob}\Gbf} = \pt/\Gbf^{\Frob}$. Since the map $\pi : \pt \rightarrow \pt/\Gbf^{\Frob}$ is surjective and finite, the functor $\pi^!$ is monadic hence $\DDic(\frac{\Gbf}{\Ad_{\Frob}\Gbf}, \Lambda) = \pi^!\pi_!-\Mod(\DD(\Lambda))$. By proper base change, there is an isomorphism of monads $\pi^!\pi_! = \Lambda[\Gbf^{\Frob}] \otimes_{\Lambda} -$ where $\Lambda[\Gbf^{\Frob}]$ is equipped with its group algebra structure. Hence $\pi^!\pi_!-\Mod(\DD(\Lambda)) = \DD(\Lambda[\Gbf^{\Frob}])$ and the result follows.
\end{proof}

More generally there is a decomposition 
\begin{equation*}
\frac{\Gbf}{\Ad_{\Frob} \Gbf} = \bigsqcup_{\gamma \in H^1(\Frob, \pi_0(\Gbf))} \pt/(\Gbf^{\gamma\Frob}).  
\end{equation*}

\begin{corol}
There is a canonical decomposition 
\begin{equation*}
\DDic(\frac{\Gbf}{\Ad_{\Frob}\Gbf}, \Lambda) = \bigoplus_{\gamma \in H^1(\Frob, \pi_0(\Gbf))} \DD(\Rep_{\Lambda} \Gbf^{\gamma\Frob}). 
\end{equation*}
\end{corol}

Let $\Bbf = \Tbf\Ubf \subset \Gbf^{\circ}$ be a Borel pair of $\Gbf^{\circ}$. We assume that $(\Bbf,\Tbf)$ is $\Frob$-stable. We denote by $\Nbf$ the normalizer of $\Tbf$ in $\Gbf$, by $\Nbf^{\circ}$ the normalizer of $\Tbf$ in $\Gbf^{\circ}$, by $\Wbf = \Nbf/\Tbf$ the corresponding Weyl group and by $\Wbf^{\circ} = \Nbf^{\circ}/\Tbf$. There is a short exact sequence of groups
$$ 1 \rightarrow \Wbf^{\circ} \rightarrow \Wbf \rightarrow \pi_0(\Gbf) \rightarrow 1.$$  

\begin{rque}
If $\Gbf$ is connected, then we have $\Wbf = \Wbf^{\circ}$. In which case we will prefer the notation $\Wbf$. 
\end{rque}

\begin{lem}
Let $x \in \Wbf$ and let $\dot{x} \in \Nbf$ be a lift of this element. Consider the Bruhat stratum $\Bbf x\Bbf \subset \Gbf_x$. Then there is an isomorphism of stacks 
\begin{equation*}
\frac{\Ubf \backslash \Bbf x\Bbf/\Ubf}{\Ad_{\Frob}(\Tbf)} = \pt/(\Tbf^{\dot{x}\Frob} \ltimes (\Ubf \cap \Ad(\dot{x})\Ubf)). 
\end{equation*}
We denote by $k_{\dot{x}} : \frac{\Ubf \backslash \Bbf x\Bbf/\Ubf}{\Ad_{\Frob}(\Tbf)} \rightarrow \pt/(\Tbf^{\dot{x}\Frob})$, then 
\begin{equation*}
k_{\dot{x},*} : \DDic(\frac{\Ubf \backslash \Bbf x\Bbf/\Ubf}{\Ad_{\Frob}(\Tbf)}, \Lambda) \rightarrow \DDic(\Rep_{\Lambda}(\Tbf^{\dot{x}\Frob}), \Lambda)
\end{equation*}
is an equivalence.
\end{lem}

\begin{proof}
The choice of $\dot{x}$ induces an isomorphism 
\begin{equation*}
\Bbf x\Bbf = \Bbf \times (\Ubf \cap \Ad(\dot{x})(\Ubf^{-}))
\end{equation*}
where $\Ubf^{-}$ is the unipotent radical of the Borel opposite to $\Bbf$, and the map from right to left is $(b,u) \mapsto b\dot{x}u$. This map is compatible with the left and right actions of $\Bbf$ as $\Ubf \cap \Ad(\dot{x})(\Ubf^{-})$ is normalized by $\Bbf$. Taking quotients by $\Ubf \times \Ubf$ yields
$$\frac{\Ubf \backslash \Bbf x\Bbf/\Ubf}{\Ad_{\Frob}(\Tbf)} = \frac{\Tbf\dot{x}}{\Ad_{\Frob}(\Tbf) \ltimes (\Ubf \cap \Ad_{\dot{x}}\Ubf)} = \pt/(\Tbf^{\dot{x}\Frob} \ltimes (\Ubf \cap \Ad_{\dot{x}}(\Ubf)).$$ This proves the first claim. The second claim follows from the fact that the map $k_{\dot{x}}$ is a gerbe banded by a connected unipotent group.
\end{proof}

\begin{rque}
More generally, let $\Pbf = \Lbf\Vbf$ be a parabolic with $\Lbf$ stable under $\Frob$ then the stratum $\frac{\Vbf \backslash \Pbf\Frob(\Pbf)/\Frob(\Vbf)}{\Ad_{\Frob}\Lbf} \subset \frac{\Vbf \backslash \Gbf/\Frob(\Vbf)}{\Ad_{\Frob}\Lbf}$ is canonically isomorphic to $\pt/(\Lbf^{\Frob} \ltimes \Vbf \cap \Frob(\Vbf))$. Furthermore the map $l_{\Lbf} : \pt/(\Lbf^{\Frob} \ltimes \Vbf \cap \Frob(\Vbf)) \rightarrow \pt/\Lbf^{\Frob}$ induces an equivalence $l_{\Lbf,*} : \DDic(\pt/(\Lbf^{\Frob} \ltimes \Vbf \cap \Frob(\Vbf)), \Lambda) \rightarrow \DDic(\pt/\Lbf^{\Frob}, \Lambda)$.
\end{rque}

We now recall the definition of the Deligne--Lusztig varieties \cite{DeligneLusztig}. 
\begin{defi}\label{defiDeligneLusztigVarieties}
Assume that $\Gbf$ is connected and recall that $\Lcal(g) = g^{-1}\Frob(g)$ is the Lang map,
\begin{enumerate}
\item Let $\Pbf = \Lbf\Vbf$ be a parabolic with Levi $\Lbf$ that is stable under $\Frob$, then we define 
$$Y_{\Lbf \subset \Pbf} = \{g\Vbf, \Lcal_{\Frob}(g) \in \Vbf\Frob(\Vbf)\} \subset \Gbf/\Vbf.$$
\item Let $\Bbf = \Tbf\Ubf$ be a $\Frob$-stable Borel pair and let $w \in \Wbf$ be an element in the Weyl and let $n \in \Nbf$ be a lift of $w$, then we set $$Y(n) = \{g\Ubf, \Lcal_{\Frob}(g) \in \Ubf n\Ubf\} \subset \Gbf/\Ubf.$$
\end{enumerate}
\end{defi}

The following is well-known, see for instance \cite{DigneMichel} or \cite{Cabanes} for proofs, using the notations of Definition \ref{defiDeligneLusztigVarieties}
\begin{enumerate}
\item The varieties $Y_{\Lbf \subset \Pbf}$ and $Y(n)$ are smooth and $Y(n)$ is of dimension $\ell(n)$. 
\item The variety $Y_{\Lbf \subset \Pbf}$ is equipped with actions of $\Gbf^{\Frob}$ and $\Lbf^{\Frob}$ coming from the action left and right translations on $\Gbf/\Vbf$.
\item The variety $Y(n)$ is equipped with actions of $\Gbf^{\Frob}$ and $\Tbf^{n\Frob}$ coming from left and right translations on $\Gbf/\Ubf$. 
\end{enumerate}

We denote by $R_n = \RGamma_c(Y(n), \Lambda)$ and $R_{\Lbf \subset \Pbf} = \RGamma_c(Y_{\Lbf \subset \Pbf}, \Lambda)$ the cohomology with compact support of these varieties. These are complexes equipped with actions of $\Gbf^{\Frob} \times \Tbf^{n\Frob}$ and $\Gbf^{\Frob} \times \Lbf^{\Frob}$ respectively. They define the Deligne--Lusztig induction and restriction functors as follows 
\begin{align*}
\Rcal_{\Lbf \subset \Pbf} : \DD(\Rep_{\Lambda} \Lbf^{\Frob}) &\rightarrow \DD(\Rep_{\Lambda} \Gbf^{\Frob}) \\
M &\mapsto R_{\Lbf \subset \Pbf} \otimes_{\Lbf^{\Frob}} M \\
{^*}\Rcal_{\Lbf \subset \Pbf} : \DD(\Rep_{\Lambda} \Gbf^{\Frob}) &\rightarrow \DD(\Rep_{\Lambda} \Lbf^{\Frob}) \\
N &\mapsto \RHom_{\Gbf^{\Frob}}(R_{\Lbf \subset \Pbf}, N) \\
\Rcal(n) : \DD(\Rep_{\Lambda} \Tbf^{n\Frob}) &\rightarrow \DD(\Rep_{\Lambda} \Gbf^{\Frob}) \\
M &\mapsto R(n) \otimes_{\Tbf^{n\Frob}} M \\
{^*}\Rcal(n) : \DD(\Rep_{\Lambda} \Gbf^{\Frob}) &\rightarrow \DD(\Rep_{\Lambda} \Tbf^{n\Frob}) \\
N &\mapsto \RHom_{\Gbf^{\Frob}}(R(n), N).
\end{align*}
The functors $(\Rcal_{\Lbf \subset \Pbf}, {^*}\Rcal_{\Lbf \subset \Pbf})$ and $(\Rcal(n), {^*}\Rcal(n))$ are left and right adjoints. 

\begin{lem}\label{lemCompareDLHC}
Assume that $\Gbf$ is connected, then there are isomorphisms of functors 
\begin{enumerate}
\item $\DD(\Rep_{\Lambda} \Lbf^{\Frob}) \rightarrow \DD(\Rep_{\Lambda} \Gbf^{\Frob})$ between 
$$\Rcal_{\Lbf \subset \Pbf} = \chfrak_{\Pbf, \Frob}i_{\Lbf \subset \Pbf, !}l_{\Lbf}^*,$$
\item $\DD(\Rep_{\Lambda} \Tbf^{n\Frob}) \rightarrow \DD(\Rep_{\Lambda} \Gbf^{\Frob})$ between 
$$\Rcal(n) = \chfrak_{\Bbf, \Frob}i_{n,!}k_n^*.$$
\end{enumerate}
\end{lem}

\begin{proof}
We show the second one as the first one is shown in a similar way. Denote by $\tilde{Y}(w) = \{g \in \Gbf/\Ubf, \Lcal(g) \in \BwB \} \subset \Gbf/\Ubf$, it is a subscheme of $\Gbf/\Ubf$ that is stable under the left action of $\Gbf^{\Frob}$ and the right action of $\Tbf$. The quotient $X(w) = \tilde{Y}(w)/\Tbf = \{g \in \Gbf/\Bbf, \Lcal(g) \in \BwB\} \subset \Gbf/\Bbf$ is the classical Deligne--Lusztig variety of \cite{DeligneLusztig} and satisfies $Y(n)/\Tbf^{n\Frob} = X(w)$. By definition of $Y(n)$ there is a closed immersion $Y(n) \subset \tilde{Y}(w)$ that is stable under the actions of $\Gbf^{\Frob} \times \Tbf^{n\Frob}$. In particular, it induces an isomorphism $Y(n) \times^{\Tbf^{n\Frob}} \Tbf \to \tilde{Y}(w)$. 
The Lang map $\Lcal : \Gbf \to \Gbf$, it induces an isomorphism of stacks $\frac{\Gbf}{\Ad_{\Frob}\Bbf} = \Gbf^{\Frob}\backslash \Gbf/\Bbf$, and we therefore have a Cartesian diagram 
\[\begin{tikzcd}
	{\frac{\BwB}{\Ad_{\Frob}\Bbf}} & {\Gbf^{\Frob}\backslash X(w)} \\
	{\frac{\Gbf}{\Ad_{\Frob}\Bbf}} & {\Gbf^{\Frob}\backslash \Gbf/\Bbf.}
	\arrow[Rightarrow, no head, from=1-1, to=1-2]
	\arrow[hook, from=1-2, to=2-2]
	\arrow[hook, from=1-1, to=2-1]
	\arrow[Rightarrow, no head, from=2-1, to=2-2]
\end{tikzcd}\]
The same holds when we replace $\Bbf$ with $\Ubf$, we have a diagram where all three squares are Cartesian 
\[\begin{tikzcd}
	{\Gbf^{\Frob}\backslash \tilde{Y}(w)} & {\Ubf \backslash \BwB/\Ubf} \\
	{\Gbf^{\Frob}\backslash X(w)} & {\frac{\Ubf \backslash \BwB/\Ubf}{\Ad_{\Frob}\Tbf}} & {\Ubf \backslash \Gbf/\Ubf} \\
	& {\Gbf^{\Frob}\backslash \Gbf/\Bbf} & {\frac{\Ubf \backslash \Gbf/\Ubf}{\Ad_{\Frob}\Tbf}.}
	\arrow[from=2-1, to=3-2]
	\arrow[from=3-2, to=3-3]
	\arrow[from=2-2, to=3-3]
	\arrow[from=2-1, to=2-2]
	\arrow[from=2-3, to=3-3]
	\arrow[from=1-2, to=2-3]
	\arrow[from=1-2, to=2-2]
	\arrow[from=1-1, to=1-2]
	\arrow[from=1-1, to=2-1]
\end{tikzcd}\] 
The choice of $n$ yields an isomorphism $\frac{\Ubf \backslash \BwB/\Ubf}{\Ad_{\Frob}\Tbf} = \pt/(\Tbf^{n\Frob} \rtimes \Ubf_n)$. The map $\Gbf^{\Frob}\backslash X(w) \to \pt/(\Tbf^{n\Frob} \rtimes \Ubf_n) \to \pt/\Tbf^{n\Frob}$ then corresponds to a $\Tbf^{n\Frob}$-torsor and since $Y(n)$ is a $\Tbf^{n\Frob}$-reduction of $\tilde{Y}(w)$ this torsor is $Y(n)$. On the other hand the first map $\Gbf^{\Frob}\backslash X(w) \to \pt/(\Tbf^{n\Frob} \rtimes \Ubf_n)$ corresponds to a $\Tbf^{n\Frob} \rtimes \Ubf_n$-torsor $\Gbf^{\Frob}\backslash\tilde{X}(w) \to \Gbf^{\Frob}\backslash X(w)$ fitting into the following diagram 
\[\begin{tikzcd}
	{\Gbf^{\Frob}\backslash \pt} & {\Gbf^{\Frob}\backslash\tilde{X}(w)} \\
	& {\Gbf^{\Frob}\backslash Y(n)} \\
	{\Gbf^{\Frob}\backslash X(w)} && \pt \\
	& {\pt/(\Tbf^{n\Frob} \rtimes \Ubf_n)} \\
	& {\pt/\Tbf^{n\Frob}.}
	\arrow["a"{description}, from=4-2, to=5-2]
	\arrow["{s_1}"{description}, from=3-1, to=5-2]
	\arrow["{t_2}"{description}, from=3-3, to=4-2]
	\arrow["{t_1}"{description}, from=3-3, to=5-2]
	\arrow["{q_2}"{description}, from=1-2, to=3-3]
	\arrow["b"{description}, from=1-2, to=2-2]
	\arrow["{p_1}"{description}, from=2-2, to=3-1]
	\arrow["{q_1}"{description}, from=2-2, to=3-3]
	\arrow["\gamma"{description}, from=1-2, to=1-1]
	\arrow["\beta"{description}, from=2-2, to=1-1]
	\arrow["{p_2}"{description}, from=1-2, to=3-1]
	\arrow["\alpha"{description}, from=3-1, to=1-1]
	\arrow["{s_2}"{description}, from=3-1, to=4-2]
\end{tikzcd}\]
The two squares with maps labeled by $_1$ and by $_2$ respectively are Cartesian. We claim that we have
$$\Rcal(n) = \alpha_!s_1^*$$ 
and 
$$\chfrak_{\Bbf, \Frob}i_{n,!}k_n^* = \alpha_!s_2^*a^*$$
this will conclude the proof. 
Let us prove the first claim. We have a canonical isomorphism of functors $\id = (t_{1,!}t_1^*)^{\Tbf^{n\Frob}}$ to which we apply $\alpha_!s_1^*$. We have 
\begin{align*}
\alpha_!s_1^*(t_{1,!}t_1^*)^{T^{n\Frob}} &= (\alpha_!s_1^*t_{1,!}t_1^*)^{\Tbf^{n\Frob}} \\
&= (\alpha_!p_{1,!}q_1^*t_1^*)^{\Tbf^{n\Frob}} \\
&= (\beta_!q_1^*t_1^*)^{\Tbf^{n\Frob}} \\
&= R(n) \otimes_{\Tbf^{n\Frob}} -.
\end{align*}
The first line follows from the fact that since $\alpha$ is quasi-compact $\alpha_!$ commutes with limits hence with invariants and $s_1$ is smooth hence $s_1^*$ has a left adjoints and thus commutes with limits as well. The second line comes from the Cartesian square of the diagram above. The last one is just a reformulation. The second claim is immediate from the definition of $\chfrak_{\Bbf, \Frob}$. 
\end{proof}

\subsection{Convolution patterns, Springer sheaves and conservativity of the Deligne-Lusztig functors}\label{subsectionConvolutionPatterns}

We consider the category $\DDic(\frac{\Gbf}{\Ad(\Gbf)}, \Lambda)$. This category is equipped with a well-known convolution product given by 
\begin{align*}
\DDic(\frac{\Gbf}{\Ad(\Gbf)}, \Lambda) \times &\DDic(\frac{\Gbf}{\Ad(\Gbf)}, \Lambda) \rightarrow \DDic(\frac{\Gbf}{\Ad(\Gbf)}, \Lambda) \\
(A,B) &\mapsto m_! (\pr_1^*A \otimes \pr_2^*B) = A *^{\Gbf} B, 
\end{align*}
where $\pr_1$ and $\pr_2$ are induced by the first and second projection and $m$ by the multiplication in the following diagram,
\[\begin{tikzcd}
	& {\frac{\Gbf \times \Gbf}{\Ad(\Gbf)}} & {\frac{\Gbf}{\Ad(\Gbf)}} \\
	{\frac{\Gbf}{\Ad(\Gbf)}} && {\frac{\Gbf}{\Ad(\Gbf)}}
	\arrow["{\pr_1}", from=1-2, to=2-1]
	\arrow["{\pr_2}"', from=1-2, to=2-3]
	\arrow["m", from=1-2, to=1-3]
\end{tikzcd}\]
where the action of $\Gbf$ on $\Gbf \times \Gbf$ is given by simultaneous conjugation. 

\begin{rque}
We note that we have only defined a bifunctor $\DDic(\frac{\Gbf}{\Ad(\Gbf)}, \Lambda) \times \DDic(\frac{\Gbf}{\Ad(\Gbf)}, \Lambda) \rightarrow \DDic(\frac{\Gbf}{\Ad(\Gbf)}, \Lambda)$ and we have not checked the (higher) associativity constraints. To address this $\infty$-categorical issue, we make the following remark in the $(2,1)$-category of algebraic stacks, we consider the correspondence 
\begin{equation*}
\frac{\Gbf}{\Ad(\Gbf)} \times \frac{\Gbf}{\Ad(\Gbf)} \xleftarrow{(\pr_1,\pr_2)} \frac{\Gbf \times \Gbf}{\Ad(\Gbf)} \xrightarrow{m} \frac{\Gbf}{\Ad(\Gbf)}.
\end{equation*}
This correspondence is a morphism, in the category of correspondences on algebraic stacks, and it equips $\frac{\Gbf}{\Ad(\Gbf)}$ with the structure of a monoid object in the category of correspondences. Since our $6$-functors formalism $\DDic$ is a lax monoidal functor out of correspondences, it sends a monoid on a monoidal category. By definition the monoidal structure we have obtained is $*^{\Gbf}$.
\end{rque} 

Similarly the category $\DDic(\frac{\Gbf}{\Ad_{\Frob}\Gbf}, \Lambda)$ is a module over the category $\DDic(\frac{\Gbf}{\Ad(\Gbf)}, \Lambda)$. The action is given by 
\begin{align*}
\DDic(\frac{\Gbf}{\Ad(\Gbf)}, \Lambda) \times &\DDic(\frac{\Gbf}{\Ad_{\Frob}\Gbf}, \Lambda) \rightarrow \DDic(\frac{\Gbf}{\Ad_{\Frob}\Gbf}, \Lambda)\\
(A,B) &\mapsto A *^{\Gbf} B = m_!(\pr_1^*A \otimes \pr_2^*B),
\end{align*}
where $\pr_1,\pr_2$ and $m$ are induced by the projections and multiplication in the following diagram 
\[\begin{tikzcd}
	& {\frac{\Gbf \times \Gbf}{\Ad_{\Frob}(\Gbf)}} & {\frac{\Gbf}{\Ad_{\Frob}(\Gbf)}} \\
	{\frac{\Gbf}{\Ad(\Gbf)}} && {\frac{\Gbf}{\Ad_{\Frob}(\Gbf)}}
	\arrow["{\pr_1}", from=1-2, to=2-1]
	\arrow["{\pr_2}"', from=1-2, to=2-3]
	\arrow["m", from=1-2, to=1-3]
\end{tikzcd}\]
where the action of $\Gbf$ on $\Gbf \times \Gbf$ is given by $g.(x,y) = (gxg^{-1}, gy\Frob(g^{-1})$. As in the previous remark, this action comes from an action of the monoid $\frac{\Gbf}{\Ad(\Gbf)}$ in correspondences on the object $\frac{\Gbf}{\Ad_{\Frob}\Gbf}$, where the action map is given by the correspondence
$$ \frac{\Gbf}{\Ad(\Gbf)} \times \frac{\Gbf}{\Ad_{\Frob}(\Gbf)} \xleftarrow{(\pr_1,\pr_2)} \frac{\Gbf \times \Gbf}{\Ad_{\Frob}(\Gbf)} \xrightarrow{m} \frac{\Gbf}{\Ad_{\Frob}(\Gbf)}.$$

\begin{defi}
Let $\Pbf = \Lbf\Vbf \subset \Gbf^{\circ}$ be a parabolic subgroup. We denote by $\Spr_{\Pbf}$ the parabolic Springer sheaf corresponding to $\Pbf$, this is the sheaf on $\frac{\Gbf^{\circ}}{\Ad(\Gbf^{\circ})}$ given by $r_*\Lambda[\dim \Vbf - \dim \Pbf]$ where $r : \frac{\Vbf}{\Ad(\Pbf)} \rightarrow \frac{\Gbf^{\circ}}{\Ad(\Gbf^{\circ})}$ is the natural map. 
\end{defi}

\begin{lem}[\cite{BorhoMacPherson}]\label{lemSpringerContainsDirac}
The sheaf $\Spr_{\Pbf}$ is a perverse sheaf and it there is an inclusion in $\Perv(\frac{\Gbf^{\circ}}{\Gbf^{\circ}}, \Lambda), \delta_1 \hookrightarrow \Spr_{\Pbf}$, where $\delta_1$ is the skyscrapper sheaf at $1$.
\end{lem}

\begin{lem}[\cite{MirkovicVilonen}]\label{lemMirkovicVilonen}
Assume that $\Gbf$ is connected. There is an isomorphism of functors 
\begin{equation*}
\chfrak_{\Pbf,\Frob} \hcfrak_{\Pbf,\Frob} (-) = \Spr_{\Pbf} *^{\Gbf} -
\end{equation*}
\end{lem}

\begin{corol}\label{corolConservativity}
Assume that $\Gbf$ is connected then for all $\Pbf$ the functor $\hcfrak_{\Pbf,\Frob}$ is conservative on the bounded below part of the category. 
\end{corol}

\begin{proof}
This argument is an adaptation of an argument of \cite[Theorem 1.5]{BBM}. It is enough to show that the functor $\chfrak_{\Pbf,\Frob}\hcfrak_{\Pbf,\Frob}$ is conservative. By Lemma \ref{lemMirkovicVilonen}, this functor is isomorphic to $\Spr_{\Pbf} *^{\Gbf} -$. We show that for a nonzero $A \in \DDic^-(\frac{\Gbf}{\Ad_{\Frob}\Gbf}, \Lambda)$ then $\Spr_{\Pbf} *^{\Gbf} A$ is nonzero. First we assume that $A$ is constructible and up to shifting $A$, we assume that $A$ is concentrated in positive perverse degrees and that ${^p}H^0(A) \neq 0$. Then, as the map $m$ defining the convolution is affine, the functor $m_!$ is left perverse $t$-exact and therefore the convolution operation is also left perverse $t$-exact. Thus we have ${^p}H^0(\Spr_{\Pbf} *^{\Gbf} A) = {^p}H^0(\Spr_{\Pbf} *^{\Gbf} {^p}H^0(A))$. As $\delta_1 \hookrightarrow \Spr_{\Pbf}$ is an injection, by left $t$-exactness, it follows that the map ${^p}H^0(A) = {^p}H^0(\delta_1 *^{\Gbf}{^p}H^0(A)) \rightarrow {^p}H^0(\Spr_{\Pbf} *^{\Gbf} {^p}H^0(A))$ is nonzero, hence ${^p}H^0(\Spr_{\Pbf} *^{\Gbf} A)$ is nonzero. The general case can be treated by a colimit argument. 
\end{proof}

\begin{rque}
We have shown that there is a nonzero map from the identity functor into $\chfrak_{\Pbf,\Frob}\hcfrak_{\Pbf,\Frob}$. If $|\Wbf|$ is invertible in $\Lambda$, then the argument of \cite{BBM} shows even that the identity is a direct summand of $\chfrak_{\Pbf,\Frob}\hcfrak_{\Pbf,\Frob}$. 
\end{rque}

\begin{corol}\label{corolConservativityDL}
Assume that $\Gbf$ is connected and let $\Frob$ be a Frobenius. The collection of Deligne-Lusztig restriction functors ${^*}\Rcal_w : \DD(\Rep_{\Lambda} \Gbf^{\Frob}) \rightarrow \DD(\Rep_{\Lambda} \Tbf^{w\Frob})$ is conservative. 
\end{corol}

\begin{proof}
This is simply a matter of playing with the adjunctions. Since the functor $\hcfrak_{\Frob}$ is conservative, the collection of costalks $i_w^!\hcfrak_{\Frob}$ is also conservative, but by Lemma \ref{lemCompareDLHC}, these functors are isomorphic up to shifts to the functors ${^*}\Rcal_w$. 
\end{proof}

\begin{corol}[\cite{DeligneLusztig}, \cite{BonnafeRouquierCatDerivee}]\label{corolGenerationDL}
The Deligne-Lusztig cohomology complexes $\RGamma_c(Y(\dot{w}, \Lambda))$ generate $\Perf(\Lambda[\Gbf^{\Frob}])$ and for all irreducible $\rho \in \Irr_{\Lambda} \Gbf^{\Frob}$ there exists a pair $(w,j)$ where $j \in \Z$ such that $\Hom(H^j_c(Y(\dot{w}), \Lambda), \rho)$ is nonzero. 
\end{corol}

\begin{proof}
Arguing as in \cite[Section 9]{BonnafeRouquierCatDerivee}, it is enough to show the second statement. But this one follows immediately from Corollary \ref{corolConservativityDL}.
\end{proof}

\begin{rque}
This proof yields a more conceptual proof of one of the theorems of Deligne and Lusztig \cite{DeligneLusztig} that the cohomology of the Deligne-Lusztig varieties contain all irreducible representations of $\Gbf^{\Frob}$. The original proof relied on character computations to deduce the $\Qlb$-case, the $\Flb$-case was considered by \cite{BonnafeRouquierCatDerivee} and was deduced from the $\Qlb$-case using Brauer theory. We note that the argument we provide is entirely geometric and does not rely on character computations. 
\end{rque}

\section{Free monodromic Hecke categories}\label{sectionFreeMonoHecke}

\subsection{The Hecke category}\label{sectionDefiHecke}

Recall that $\Gbf$ is an algebraic group with reductive neutral component and that we have fixed $\Bbf = \Tbf\Ubf$ a Borel pair in $\Gbf^{\circ}$. We denote by $\Wbf = \Nbf(\Tbf)/\Tbf$ the Weyl group of $\Gbf$. We fix $\Lambda \in \{\Flb, \Zlb, \Qlb\}$ and we will denote $\RR_{\Tbf} = \RR_{\Tbf,\Lambda}$ and $\Ch(\Tbf) = \Ch_{\Lambda}(\Tbf)$. We also denote by $\Lambda_0 \in \{\Fl, \Zl, \Ql\}$ the corresponding subring of $\Lambda$. We use the formalism of monodromic sheaves that we have established to propose a simpler construction of the free monodromic Hecke categories of \cite{BezYun}, \cite{BezrukRiche} and \cite{Gouttard}. Our construction also simplifies the construction of the monoidal structure of these categories. 

Consider the stack $\Ubf \backslash \Gbf/\Ubf$. There are three actions of tori that we can consider : 
\begin{enumerate}
\item the action of $\Tbf$ induced by left translations $\Ubf \backslash \Gbf/\Ubf$, 
\item the action of $\Tbf$ on $\Ubf \backslash \Gbf/\Ubf$ induced by the action of $\Tbf$ on $\Gbf$ given by $t.x = xt$, we will refer to this action as the right action of $\Tbf$,
\item the action of $\Tbf \times \Tbf$ induced by left and right translations.
\end{enumerate}

We denote by 
\begin{enumerate}
\item $\HH^{\lef} = \bigoplus_{\chi \in \Ch(\Tbf)} \DDic(\Ubf \backslash \Gbf/\Ubf, \RR_{\Tbf})_{\chi}^{\lef}$ where the equivariance is relative to the action of $\Tbf$ on the left.
\item $\HH^{\rig} = \bigoplus_{\chi \in \Ch(\Tbf)} \DDic(\Ubf \backslash \Gbf/\Ubf, \RR_{\Tbf})_{\chi}^{\rig}$ where the equivariance is relative to the action of $\Tbf$ on the right.
\item $\HH^{\lef,\rig} = \bigoplus_{\chi, \chi' \in \Ch(\Tbf)} \DDic(\Ubf \backslash \Gbf/\Ubf, \RR_{\Tbf \times \Tbf})_{\chi, \chi'}^{\lef, \rig}$ where the equivariance is relative to the action of $\Tbf\times \Tbf$ on the right and the index $(\chi, \chi')$ refer to sheaves that are equivariant for $L_{\Tbf \times \Tbf}\otimes_{\Lambda} (\mathcal{L}_{\chi} \boxtimes_{\Lambda} \mathcal{L}_{\chi'})$. 
\end{enumerate}

We equip the space $\Ubf \backslash \Gbf/\Ubf$ with its Bruhat stratification. The strata are indexed by the Weyl group $\Wbf$, and the stratum corresponding to $w \in \Wbf$ is $\Ubf \backslash \BwB /\Ubf$. We denote by $i_w : \Ubf \backslash \BwB /\Ubf \rightarrow \Ubf \backslash \Gbf/\Ubf$. The next lemma is a direct application of Lemma \ref{lemExampleTwoActions}.

\begin{lem}\label{lemTorus}
The category $\DDic(\Tbf, \RR_{\Tbf \times \Tbf})^{\lef, \rig}_{\chi, \chi'}$ is zero unless $\chi' = \chi$. In this case, this category is equivalent to $\DD(\RR_{\Tbf})$.
\end{lem}

\begin{lem}\label{lemCompactGeneration}
All three categories $\HH^{\lef}, \HH^{\rig}$ and $\HH^{\lef, \rig}$ are compactly generated. The categories of compact objects are the categories 
\begin{enumerate}
\item $\bigoplus_{\chi \in \Ch(\Tbf)} \DDc(\Ubf \backslash \Gbf/\Ubf, \RR_{\Tbf})_{\chi}^{\lef}$, 
\item $\bigoplus_{\chi \in \Ch(\Tbf)} \DDc(\Ubf \backslash \Gbf/\Ubf, \RR_{\Tbf})_{\chi}^{\rig}$,
\item and $\bigoplus_{\chi, \chi' \in \Ch(\Tbf)} \DDc(\Ubf \backslash \Gbf/\Ubf, \RR_{\Tbf})_{\chi, \chi'}^{\lef, \rig}$,	
\end{enumerate}
respectively. 
\end{lem}

\begin{proof}
First since the inclusion $i_w$ are quasi-compact and schematic all functors $i_{w,!}, i_w^!, i_{w,*}, i_{w}^*$ between the categories $\HH^{?}$ and $\bigoplus_{\chi \in \Ch(\Tbf)} \DDic(\Ubf \backslash \BwB /\Ubf, \RR_{\Tbf})_{\chi}^{?}$ where $? \in \{\lef, \rig, (\lef,\rig)\}$ commute with arbitrary direct sums. For $i_{w,!}$ and $i_{w}^*$ this is clear since they are left adjoints. We show it for $i_{w,*}$, the case $i_{w}^!$ can be deduced from the case of $i_{w,*}$ using excision triangles. We only need to check that the canonical map $\oplus_i i_{w,*}A_i \rightarrow i_{w,*}\oplus_i A_i$ for $A_i \in \HH^{?}$ is an isomorphism, since the functor $i_{w,*}$ commutes with the forgetful functor and smooth pullbacks this can checked after pulling back to $\Gbf/\Ubf$, where this now follows from the statement on schemes by \cite[6.4.5]{BhattScholze}. 

Since $i_w^!$ and $i_{w,*}$ are continuous, their left adjoints preserve compact objects. We now show the lemma by induction on the strata. Denote by $V \xrightarrow{j} \Ubf \backslash \Gbf/\Ubf \xleftarrow{i} Z$ the inclusion of the open stratum $V$ and $Z$ its closed complement. Using the exicision triangles for $A \in \HH^{?}$, 
\begin{equation*}
j_!j^* A \rightarrow A \rightarrow i_*i^*A,
\end{equation*}
we see that $A$ is a colimit of compact objects if and only if $j^*A$ and $i^*A$ are so. By induction this reduces to showing that $\bigoplus_{\chi \in \Ch(\Tbf)} \DDic(\Ubf \backslash \BwB /\Ubf, \RR_{\Tbf})_{\chi}^{?}$ is compactly generated. But this category is equivalent to the category $\bigoplus_{\chi} \DD(\RR_{\Tbf})$. This is clear for $? \in \{\lef, \rig \}$ and by Lemma \ref{lemTorus} for the case of the action of $\Tbf \times \Tbf$. This proves the compact generation statement. 

We now identify the compact objects. Again by induction on the strata, and using the same triangle, we see that an object $A$ is compact if and only if for all $w, i_w^*A$ is compact. Hence the category of compact objects is the stable category generated by all objects of the form $i_{w,!}A$ for varying $w$ and $A \in \bigoplus_{\chi \in \Ch(\Tbf)} \DDic(\Ubf \backslash \BwB /\Ubf, \RR_{\Tbf})_{\chi}^{?}$ a compact object. On $\Ubf \backslash \BwB /\Ubf$ the category of compact objects is $\bigoplus_{\chi \in \Ch(\Tbf)} \DDc(\Ubf \backslash \BwB /\Ubf, \RR_{\Tbf})_{\chi}^{?}$. And the category generated by all $i_{w,!}A$ for varying $w$ and $A$ compact is then $\bigoplus_{\chi} \DDc(\Ubf \backslash \Gbf /\Ubf, \RR_{\Tbf})_{\chi}^{?}$.
\end{proof}

The inclusions of $\Tbf \xrightarrow{i_{\lef}} \Tbf \times \Tbf \xleftarrow{i_{\rig}} \Tbf$ given by $i_{\lef}(t) = (t,1)$ and $i_{\rig}(t) = (1,t)$ induce inclusions $\RR_{\Tbf} \xrightarrow{i_{\lef, *}} \RR_{\Tbf \times \Tbf} \xleftarrow{i_{\rig, *}} \RR_{\Tbf}$. 

\begin{lem}\label{lemConstructibility}
There are well defined continuous compact preserving functors 
\begin{equation*}
\HH^{\lef} \xleftarrow{\For^{\lef}} \HH^{\lef, \rig} \xrightarrow{\For^{\rig}} \HH^{\rig}
\end{equation*}
induced by forgetting the $(\Tbf \times \Tbf, L_{\Tbf \times \Tbf} \otimes_{\Lambda} (\mathcal{L}_{\chi} \boxtimes_{\Lambda} \mathcal{L}_{\chi'}))$-equivariance along $i_{\lef, *}$ and $i_{\rig, *}$ respectively. 
\end{lem}

\begin{proof}
To check that these functors are well defined we have to check that the functors $\For^{\lef}$ and $\For^{\rig}$ preserve constructibility. Let $A \in \DDc(\Ubf \backslash \Gbf/\Ubf, \RR_{\Tbf \times \Tbf})^{\lef, \rig}_{\chi, \chi'}$. As in the previous lemma, we can assume $A = i_{w,!}A_0$ for some object 
$$A_0 \in \bigoplus_{\chi, \chi' \in \Ch(\Tbf)} \DDc(\Ubf \backslash \BwB /\Ubf, \RR_{\Tbf\times T})_{\chi, \chi'}.$$ 
We can further assume that $A_0 \in  \DDc(\Ubf \backslash \BwB /\Ubf, \RR_{\Tbf\times T})_{\chi, \chi'}$ and that $w\chi' = \chi$ otherwise this category is zero. By Lemma \ref{lemTorus}, the category $\DDc(\Ubf \backslash \BwB /\Ubf, \RR_{\Tbf\times T})_{\chi, \chi'}$ is equivalent to $\DD_{\coh}(\RR_{\Tbf})$. We can assume that $A_0$ corresponds to $\RR_{\Tbf}$ as this objects generates the stable category $\DD_{\coh}(\RR_{\Tbf})$. Therefore $A_0 \simeq \nu_w^*((L_{\Tbf \times \Tbf} \otimes_{\RR_{\Tbf \times \Tbf}} \RR_{\Tbf}) \otimes_{\Lambda} \mathcal{L}_{\chi})[\ell(w) + \dim \Tbf]$ as an $\RR_{\Tbf \times \Tbf}$-sheaf. Since $L_{\Tbf \times \Tbf} \otimes_{\RR_{\Tbf \times \Tbf}} \RR_{\Tbf} \simeq L_{\Tbf}$ as an $\RR_{\Tbf}$-sheaf after forgetting along either the left of right inclusion, we get the desired constructibility statement. 
\end{proof}

\begin{lem}\label{lemEquivDefinition}
Both functors $\For^{\lef}$ and $\For^{\rig}$ are equivalences. 
\end{lem}

\begin{proof}
We first note that $\For^{\lef}$ and $\For^{\rig}$ induce equivalence on each strata. 
\begin{equation*}
\bigoplus_{\chi, \chi' \in \Ch(\Tbf)} \DDc(\Ubf \backslash \BwB /\Ubf, \RR_{\Tbf\times T})_{\chi, \chi'} \simeq \bigoplus_{\chi \in \Ch(\Tbf)} \DDc(\Ubf \backslash \BwB /\Ubf, \RR_{\Tbf})_{\chi}^{\lef},
\end{equation*}
and 
\begin{equation*}
\bigoplus_{\chi, \chi' \in \Ch(\Tbf)} \DDc(\Ubf \backslash \BwB /\Ubf, \RR_{\Tbf\times T})_{\chi, \chi'} \simeq \bigoplus_{\chi' \in \Ch(\Tbf)} \DDc(\Ubf \backslash \BwB /\Ubf, \RR_{\Tbf})_{\chi'}^{\rig}.
\end{equation*}
This is an immediate application of Lemma \ref{lemTorus}.
To conclude that the functor $\For^{\lef}$ is an equivalence, we proceed by induction on the strata. Let $V \subset \Ubf \backslash \Gbf/\Ubf$ be a stratum and let $Z = \overline{V} - V$ be the closed complementary of the closure of $V$. Denote by $i$ and $j$ the inclusions $Z \subset \overline{V}$ and $V \subset \overline{V}$ respectively. Assume by induction that $\For^{\lef}$ induces an equivalence on the full subcategory of $\HH^{\lef}$ and $\HH^{\lef, \rig}$ supported on $Z$. Let $A,B \in \HH^{\lef, \rig}$ be supported on $\overline{V}$. Using excision triangles, we can assume that $A = i_*A_0$ or $A = j_!A_0$ and that $B = j_!B_0$ or $B = i_*B_0$. We now have 
\begin{enumerate}
\item if $A = i_*A_0$ and $B = i_*B_0$, then $\Hom(A,B) = \Hom(\For^{\lef}(A), \For^{\lef}(B))$ by induction,
\item if $A = j_!A_0$ and $B = j_!B_0$, then $\Hom(A,B) = \Hom(\For^{\lef}(A), \For^{\lef}(B))$ using the stratum case, 
\item if $A = j_!A_0$ and $B = i_*B_0$, then $\Hom(A,B) = 0$ and $\Hom(\For^{\lef}(A), \For^{\rig}(B)) = 0$, 
\item finally if $A = i_*A_0$ and $B = j_!B_0$, then as the forgetful functor commutes with $i^!$ and $j_!$, we have 
\begin{align*}
\Hom(A,B) = \Hom(A_0, i^!j_!B_0) 
&= \Hom(\For^{\lef}A_0, \For^{\lef}i^!j_!B_0) \\
&= \Hom(\For^{\lef}A_0, i^!j_!\For^{\lef}B_0) \\
&= \Hom(\For^{\lef}A, \For^{\lef}B).
\end{align*}
This establishes that $\For^{\lef}$ is fully faithful, as the subcategories of $\HH^{\lef, \rig}$ and $\HH^{\lef}$ of sheaves supported on $\overline{V}$ are generated by the sheaves of the form $i_*A_0$ and $j_!A_0$, we get the essential surjectivity.
\end{enumerate}

\end{proof}

\begin{rque}
Note that the functor $\For^{\rig, -1} \For^{\lef}$ is an equivalence that is \emph{not} $\RR_{\Tbf}$-linear. 
\end{rque}

From now on we denote by $\HH$ either of the categories $\HH^{\lef}, \HH^{\rig}$ or $\HH^{\lef, \rig}$ which are identified through $\For^{\lef}$ and $\For^{\rig}$. This category is equipped with its perverse $t$-structure. For most constructions, we will work with $\HH^{\rig}$. We denote by $\HH^{\RR}$ the full subcategory of compact objects and by $\HH^{\RR}_{\chi}$ the category $\DDc(\Ubf \backslash \Gbf/\Ubf, \RR_{\Tbf})_{\chi}$.

\begin{rque}
Since $\For^{\rig}$ is an equivalence, the category $\HH$ is equipped with an $\RR_{\Tbf \times \Tbf}$-linear structure.
\end{rque}

\subsection{The monoidal structure}\label{sectionMonoidalStructures}

We define the convolution structure on $\HH^{\RR}$. We do it in several steps. We will use the model of $\HH$ given by $\HH^{\lef, \rig}$.

First let $X$ be a stack. We define a functor 
\begin{align*}
\widehat{\otimes}_{\Lambda} : \DDc(X, \RR_{\Tbf}) \times \DDc(X, \RR_{\Tbf}) &\rightarrow \DDc(X, \RR_{\Tbf \times \Tbf}) \\
(A,B) \mapsto (A \widehat{\otimes}_{\Lambda} B).  
\end{align*}
For $A_0,B_0 \in \DDc(X, \RR_{\Tbf, \Lambda_0})$ we first construct a sheaf $(A_0 \widehat{\otimes}_{\Lambda_0} B_0) \in \DDc(X,\RR_{\Tbf \times T, \Lambda_0})$. First consider $A_0 \otimes_{\Lambda_0} B_0$, this is naturally an $\RR_{\Tbf, \Lambda_0} \otimes_{\Lambda_0} \RR_{\Tbf, \Lambda_0}$-sheaf on $X_{\proet}$. Let $I_{\Tbf \times \Tbf}$ be the ideal of $\RR_{\Tbf, \Lambda_0} \otimes_{\Lambda_0} \RR_{\Tbf, \Lambda_0}$ given by $\RR_{\Tbf, \Lambda_0} \otimes_{\Lambda_0} I  + I \otimes_{\Lambda_0} \RR_{\Tbf, \Lambda_0}$ where $I$ is as before the augmentation ideal of $\RR_{\Tbf, \Lambda_0}$. The ring $\RR_{\Tbf \times \Tbf}$ is then the completion of $\RR_{\Tbf, \Lambda_0} \otimes_{\Lambda_0} \RR_{\Tbf, \Lambda_0}$ along $I$. We then denote by $(A_0 \widehat{\otimes}_{\Lambda_0} B_0)$ the derived completion of $A_0 \otimes_{\Lambda_0} B_0$ along the ideal $I$ in $\DD(X_{\proet}, \RR_{\Tbf, \Lambda_0} \otimes_{\Lambda_0} \RR_{\Tbf, \Lambda_0})$ in the sense of \cite[Section 3.5]{BhattScholze}. This derived completion is the functor $\DD(X_{\proet}, \RR_{\Tbf, \Lambda_0} \otimes_{\Lambda_0} \RR_{\Tbf, \Lambda_0})\rightarrow \DD(X_{\proet}, \RR_{\Tbf \times T, \Lambda_0})$ given by
\begin{equation*}
C_0 \mapsto \lim_n (C_0 \otimes_{(\RR_{\Tbf, \Lambda_0} \otimes_{\Lambda_0} \RR_{\Tbf, \Lambda_0})} (\RR_{\Tbf, \Lambda_0} \otimes_{\Lambda_0} \RR_{\Tbf, \Lambda_0})/I^n).
\end{equation*}

\begin{lem}
For $A_0, B_0 \in \DDc(X, \RR_{\Tbf, \Lambda_0})$ the $\RR_{\Tbf \times T, \Lambda_0}$-sheaf $(A_0 \widehat{\otimes}_{\Lambda_0} B_0)$ is constructible. 
\end{lem}

\begin{proof}
Recall that an $\RR_{\Tbf \times T, \Lambda_0}$-complete sheaf $A_0$ is constructible if $A_0 \otimes_{\RR_{\Tbf \times T, \Lambda_0}} \RR_{\Tbf \times T, \Lambda_0}/I_{\Tbf \times \Tbf}$ is constructible as a $\Lambda_0$-sheaf. But we have 
\begin{equation*}
(A_0 \widehat{\otimes}_{\Lambda_0} B_0) \otimes_{\RR_{\Tbf \times T, \Lambda_0}} \RR_{\Tbf \times T, \Lambda_0}/I_{\Tbf \times \Tbf} = (A_0 \otimes_{\RR_{\Tbf, \Lambda_0}} \RR_{\Tbf, \Lambda_0}/I) \otimes_{\Lambda_0} (B_0 \otimes_{\RR_{\Tbf, \Lambda_0}} \RR_{\Tbf, \Lambda_0}/I). 
\end{equation*}
But $ (A_0 \otimes_{\RR_{\Tbf, \Lambda_0}} \RR_{\Tbf, \Lambda_0}/I)$ and $(B_0 \otimes_{\RR_{\Tbf, \Lambda_0}} \RR_{\Tbf, \Lambda_0}/I)$ are constructible $\Lambda_0$-sheaf therefore this tensor product is also constructible.
\end{proof}

\begin{defi}\label{defiConvolCompleted}
The functor 
\begin{align*}
\widehat{\otimes}_{\Lambda} : \DDc(X, \RR_{\Tbf}) \times \DDc(X, \RR_{\Tbf}) &\rightarrow \DDc(X, \RR_{\Tbf \times \Tbf}) \\
(A,B) \mapsto (A \widehat{\otimes}_{\Lambda} B).  
\end{align*}
is defined as the $\Lambda$-extension of the functor
\begin{equation*}
\widehat{\otimes}_{\Lambda_0} : \DDc(X, \RR_{\Tbf, \Lambda_0}) \times \DDc(X, \RR_{\Tbf, \Lambda_0}) \rightarrow \DDc(X, \RR_{\Tbf \times T, \Lambda_0}).
\end{equation*}
\end{defi}

Then we define a functor 
\begin{equation*}
\DDc(\Ubf \backslash \Gbf /\Ubf, \RR_{\Tbf \times \Tbf}) \times \DDc(\Ubf \backslash \Gbf/\Ubf, \RR_{\Tbf \times \Tbf})  \rightarrow \DDc(\Ubf \backslash \Gbf /\Ubf, \RR_{\Tbf\times \Tbf \times \Tbf \times \Tbf}) 
\end{equation*}
as follows. Consider the following diagram 
\[\begin{tikzcd}
	& {\Ubf \backslash \Gbf \times^{\Ubf} \Gbf/\Ubf} & {\Ubf \backslash \Gbf/\Ubf} \\
	{\Ubf \backslash \Gbf/\Ubf} && {\Ubf \backslash \Gbf/\Ubf}
	\arrow["{\pr_1}", from=1-2, to=2-1]
	\arrow["{\pr_2}"', from=1-2, to=2-3]
	\arrow["m", from=1-2, to=1-3]
\end{tikzcd}\]
where $m$ is induced by the multiplication map. Then we set for $A, B \in \DDc(\Ubf \backslash \Gbf /\Ubf, \RR_{\Tbf \times \Tbf})$
\begin{equation*}
A * B = m_!(A \widehat{\boxtimes}_{\Lambda} B)[\dim \Tbf].
\end{equation*}

\begin{lem}\label{lemConstructibleConvol}
Assume that $A \in \HHo_{[\chi_1, \chi_2]}$ and $B \in \HHo_{[\chi_3, \chi_4]}$, if $\chi_3 \neq \chi_2$ then $A * B = 0$ and in general the $\RR_{\Tbf \times \Tbf \times \Tbf \times \Tbf}$-structure on $A * B$ is constructible as an $\RR_{\Tbf \times \Tbf}$-sheaf after forgetting along the inclusion $\RR_{\Tbf \times \Tbf} \rightarrow \RR_{\Tbf \times \Tbf \times \Tbf \times \Tbf}$ induced by the outer inclusions.
\end{lem}

\begin{proof}
We argue as in \cite[4.3]{BezYun}. We decompose the map $m$ in two steps 
\begin{equation*}
\Ubf \backslash \Gbf \times^{\Ubf} \Gbf /\Ubf \xrightarrow{q} \Ubf \backslash \Gbf \times^{\Bbf} \Gbf /\Ubf \xrightarrow{\tilde{m}} \Ubf \backslash \Gbf/\Ubf,
\end{equation*}
where $\tilde{m}$ is the map induced by the multiplication in $\Gbf$ and $q$ is the quotient by the $\Tbf$-action $t.(g,g') = (gt^{-1}, tg)$, in particular the map $q$ is a $\Tbf$-torsor. It is enough to check the triviality of $q_! (A \widehat{\boxtimes}_{\Lambda} B)$ if $\chi_2 \neq \chi_3$. The triviality can be checked after reducing modulo $I_{\Tbf \times \Tbf \times \Tbf  \times \Tbf}$. We have an isomorphism 
\begin{equation*}
q_!(A \widehat{\boxtimes}_{\Lambda} B)/I_{\Tbf \times \Tbf \times \Tbf  \times \Tbf} = q_!((A/I_{\Tbf \times \Tbf}) \boxtimes_{\Lambda} (B/I_{\Tbf \times \Tbf})).
\end{equation*}
The sheaf $A/I_{\Tbf \times \Tbf}$ is $(\Tbf,\mathcal{L}_{\chi_2, \Flb})$-equivariant for the right action of $\Tbf$ and $B/I_{\Tbf \times \Tbf}$ is $(\Tbf,\mathcal{L}_{\chi_3, \Flb})$-equivariant for the left action of $\Tbf$. Hence their tensor product is $(\Tbf, \mathcal{L}_{\chi_2^{-1}\chi_3, \Flb})$-equivariant for the action of $\Tbf$ given by $t.(x,y) = (xt^{-1},ty)$. The pushforward along $q_!$ is therefore $0$ if $\chi_2 \neq \chi_3$. 
The constructibility assertion follows from the fact that $A \widehat{\boxtimes}_{\Lambda} B$ is already $\RR_{\Tbf \times \Tbf}$-constructible after forgetting along the outer inclusions. This follows from Lemma \ref{lemEquivDefinition}, indeed $A$ is $\RR_{\Tbf}$-constructible after forgetting the right action of $\RR_{\Tbf}$ and $B$ is constructible after forgetting the left action of $\RR_{\Tbf}$. 
\end{proof}

We can now define the convolution functor
\begin{align*}
\HHo \times \HHo &\rightarrow \HHo \\
(A,B) &\mapsto \For_{\mathrm{ext}}m_!(A \widehat{\boxtimes}_{\Lambda} B), 
\end{align*}
where $\For_{\mathrm{ext}}$ is the forgetful functor induced by the map $\RR_{\Tbf \times \Tbf} \rightarrow \RR_{\Tbf \times \Tbf \times \Tbf \times \Tbf}$ induced by the outer inclusions.

\begin{corol}
The convolution defines a monoidal structure on $\HHo$ and $\HHo_{[\chi, \chi]}$.
\end{corol}

As $\HH = \Ind(\HHo)$, we can extend the monoidal structure to $\HH$ using the universal property of ind-completions. In particular we have a monoidal structure on $\HH$ that is continuous on both variables.

\begin{rque}\label{rquePushForward}
In the definition of the convolution, for a pair of objects $A \in \HH_{[\chi_1, \chi_2]}$ and $B \in \HH_{[\chi_2, \chi_3]}$ the sheaf $A \widehat{\boxtimes}_{\Lambda} B$ is $\Tbf$-unipotent monodromic for the $\Tbf$-action defining the quotient $\Ubf \backslash \Gbf \times^{\Ubf} \Gbf/\Ubf \to \Ubf \backslash \Gbf \times^{\Bbf} \Gbf/\Ubf$ hence by Lemma \ref{lemPushForwardUnipMono}, we can replace $m_!$ by $m_*$ up to shift. 
\end{rque}

\subsection{Standard and costandard}

Let $n \in \Nbf(\Tbf)$ be an element in the normalizer of $\Tbf$, the choice of this element defines a splitting $\Bbf n\Bbf = \Ubf \times \Tbf \times (\Ubf \cap \Ad(n)\Ubf)$, we denote by $p_n : \Bbf n\Bbf \rightarrow \Tbf$ the projection on $\Tbf$. We also denote by $i_n$ the inclusion of the stratum $\Ubf \backslash \Bbf n\Bbf /\Ubf \subset \Ubf \backslash \Gbf/\Ubf$ it depends only on the image of $n$ in the $\Wbf$. 

\begin{defi}
We define the standard and costandard sheaves 
\begin{enumerate}
\item $\Delta_{n, \chi} = i_{n,!} p_n^*L_{\Tbf, \chi}[\dim \Tbf + \ell(n)]$, 
\item $\nabla_{n, \chi} = i_{n,*} p_n^*L_{\Tbf, \chi}[\dim \Tbf + \ell(n)]$.
\end{enumerate} 
\end{defi}

Since the map $i_n$ is affine, all the sheaves $\Delta_{n, \chi}$ and $\nabla_{n, \chi}$ are perverse. The next lemma is classical, see for instance \cite[7.7]{BezrukRiche} and \cite[8.4.1]{Gouttard}.

\begin{lem}
There are isomorphisms 
\begin{enumerate}
\item if $\Bbf n\Bbf n'\Bbf = \Bbf nn'\Bbf$ then $\Delta_{n,n'^{-1}\chi} * \Delta_{n', \chi} = \Delta_{nn', \chi}$
\item if $\Bbf n\Bbf n'\Bbf = \Bbf nn'\Bbf$ then $\nabla_{n,n'^{-1}\chi} * \nabla_{n', \chi} = \nabla_{nn', \chi}$
\item there are isomorphisms $\Delta_{n^{-1}, n\chi} * \nabla_{n,\chi} = \Delta_{1,\chi}$ and $\nabla_{n^{-1}, n\chi} * \Delta_{n,\chi} = \Delta_{1,\chi}$.
\end{enumerate}
\end{lem}

\begin{rque}
If both $n$ and $n'$ are in $G^{\circ}$ then the condition $$\Bbf n\Bbf n'\Bbf = \Bbf nn'\Bbf$$ is equivalent to the condition that $\ell(n) + \ell(n') = \ell(nn')$. 
\end{rque}

\begin{rque}
It follows from these identities that all the sheaves $\Delta_{n, \chi}$ are dualizable. 
\end{rque} 

\subsection{Rigidity of the Hecke category}

We recall the notion of rigid and quasi-rigid categories.

\begin{defi}[\cite{BenZviNadlerComplexGroup}, \cite{GaiRoz}] 
Let $\mathcal{C} \in \Pr_{\Lambda}$ be a monoidal category and assume that $\mathcal{C}$ is compactly generated. 
\begin{enumerate}
\item The category $\mathcal{C}$ is quasi-rigid if it is compactly generated by left and right dualizable objects. 
\item The category $\mathcal{C}$ is rigid if it is quasi-rigid and the monoidal unit is compact.
\end{enumerate}
\end{defi}

\begin{thm}\label{thmQuasiRigidity}
The category $\HH$ is quasi-rigid. Let $\mathfrak{o}$ be a $\Wbf$-orbit in $\Ch(\Tbf)$, then the category $\HH_{\mathfrak{o}}$ is rigid. 
\end{thm}

\begin{rque}
The only obstruction to $\HH$ being rigid is the fact that the unit is $\oplus_{\chi \in \Ch(\Tbf)} \Delta_{1, \chi}$ which is not compact. But once  we restrict a single $\Wbf$-orbit, this sum becomes finite and the unit is compact.
\end{rque}

\begin{proof}
The category $\HH$ is generated by the $\Delta_{w, \chi}$ and they are all compact and dualizable.
\end{proof}

\begin{rque}
The category of dualizable objects in $\HH$ is stable under all finite colimits and therefore all compact objects are dualizable. 
\end{rque}

We now study the functor of taking left and right adjoints.
\begin{defi}\label{defiDualityFunctor}
We denote by 
\begin{align*}
\D^{-} : \HHo \rightarrow \HHo, 
M \mapsto \inv^*\D'(M)(\varepsilon)[-\dim \Tbf],
\end{align*}
where $\inv : \Gbf \rightarrow \Gbf$ denotes the inversion map.
\end{defi}

\begin{lem}\label{lemDuality}
There is a canonical isomorphism 
\begin{equation*}
\D'(- * -) = (\D'(-) * \D'(-))[-2\dim \Tbf].
\end{equation*}
\end{lem}

\begin{proof}
Recall that the convolution \ref{defiConvolCompleted} was defined as 
\begin{equation*}
A * B = \For_{\RR_{\Tbf \times \Tbf}}^{\RR_{\Tbf}}m_!(A \widehat{\boxtimes}_{\Lambda} B)[\dim \Tbf],
\end{equation*}
where the forgetful functor is induced by the second inclusion $\RR_{\Tbf} \rightarrow \RR_{\Tbf \times \Tbf}$. Consider the $\RR_{\Tbf \times \Tbf}$-module $\RR_{\Tbf \times \Tbf}(\epsilon_{\Tbf \times \Tbf}) \otimes_{\RR_{\Tbf}} \RR_{\Tbf}(\epsilon_{\Tbf})$ where $\RR_{\Tbf} \rightarrow \RR_{\Tbf \times \Tbf}$ is induced via the second inclusion. Tensoring by this module defines a twist $M \mapsto M(\varepsilon_{\Tbf \times \Tbf/\Tbf})$ for $M \in \DD(\RR_{\Tbf})$. We claim that there are natural $\RR_{\Tbf \times \Tbf}$-linear isomorphism of functors 
\begin{enumerate}
\item $\D'_{\RR_{\Tbf}}(\For_{\RR_{\Tbf \times \Tbf}}^{\RR_{\Tbf}}(-)) = \For_{\RR_{\Tbf \times \Tbf}}^{\RR_{\Tbf}}(\D'_{\RR_{\Tbf \times \Tbf}}(-))[-\dim \Tbf](\varepsilon_{\Tbf \times \Tbf/\Tbf})$, where the index $\RR_{\Tbf}$ or $\RR_{\Tbf \times \Tbf}$ specifies where we use the version of the duality $\D'$ we use, this follows from Remark \ref{rmkRelativeTwists}.
\item $m_! = m_*[\dim \Tbf]$ on objects that are $L_{\Tbf \times \Tbf} \otimes \mathcal{L}_{\chi, \chi'}$-equivariant, this follows from Lemma \ref{lemPushForwardUnipMono}.
\end{enumerate}
Let us assume both claims. And let us show how this implies the theorem
\begin{align*}
\D'_{\RR_{\Tbf}}(A * B) &= \For_{\RR_{\Tbf \times \Tbf}}^{\RR_{\Tbf}}\D'_{\RR_{\Tbf \times \Tbf}}m_!(A \widehat{\boxtimes}_{\Lambda} B)(\varepsilon_{\Tbf \times \Tbf/\Tbf}) \\
&= \For_{\RR_{\Tbf \times \Tbf}}^{\RR_{\Tbf}}m_*\D'_{\RR_{\Tbf \times \Tbf}}(A \widehat{\boxtimes}_{\Lambda} B)(\varepsilon_{\Tbf \times \Tbf/\Tbf}) \\
&= \For_{\RR_{\Tbf \times \Tbf}}^{\RR_{\Tbf}}m_!(\D'_{\RR_{\Tbf}}(A) \widehat{\boxtimes}_{\Lambda} \D'_{\RR_{\Tbf}}(B))(\varepsilon_{\Tbf \times \Tbf/\Tbf}) [-\dim \Tbf]. 
\end{align*}
The first line follow from the first point and the last one from the compatibility between $\otimes, \Hom$ and the Kunneth formula. Note that once we forget along $\RR_{\Tbf} \rightarrow \RR_{\Tbf \times \Tbf}$ the second inclusion,  the twist $(\varepsilon_{\Tbf \times \Tbf/\Tbf})$ becomes trivial. 
\end{proof}

\begin{lem}\label{lemDuality2}
Let $X$ be a $\Tbf$-scheme, there is a canonical isomorphism 
\begin{equation*}
\Hom_{\DD(X, \RR_{\Tbf})_{\chi}}(A,B)[\dim \Tbf] = \RGamma(X, \D'(A) \widehat{\otimes}^!_{\Lambda} B)
\end{equation*}
where $\otimes^!_{\Lambda}$ denotes $\Delta^!(-\boxtimes_{\Lambda}-)$ and the completion as in Definition \ref{defiConvolCompleted}.  
\end{lem}

\begin{proof}
By descent, we can assume that $X = Y \times \Tbf$. We assume that $\chi = 1$. Let $A,B  \in \DD(X, \RR_{\Tbf})_{\unip}$. Then $A = A' \boxtimes_{\RR_{\Tbf}} L_{\Tbf}$ and $B = B' \boxtimes_{\RR_{\Tbf}} L_{\Tbf}$. We can thus compute
\begin{align*}
\Hom_{\DD(X, \RR_{\Tbf})_{\unip}}(A,B) &= \Hom_{\DDc(Y, \RR_{\Tbf})}(A',B') \\
&= \RGamma(Y, \D'_{\RR_{\Tbf}}(A) \otimes^!_{\RR_{\Tbf}} B). 
\end{align*}
On the other side, we have
\begin{equation*}
\D'(A) \widehat{\otimes}^!_{\Lambda} B = (\D_{\RR_{\Tbf}}(A) \widehat{\otimes}^!_{\Lambda} B) \boxtimes_{\RR_{\Tbf \times \Tbf}} \Delta^*L_{\Tbf \times \Tbf}.
\end{equation*}
Applying the functor $\RGamma(X, -)$ we get
\begin{align*}
\RGamma(X, \D'(A) \widehat{\otimes}^!_{\Lambda} B) &= \RGamma(Y, \D_{\RR_{\Tbf}}(A) \widehat{\otimes}^!_{\Lambda} B) \otimes_{\RR_{\Tbf \times \Tbf}} \RGamma(T, \Delta^*L_{\Tbf \times \Tbf}) \\
&= \RGamma(Y, \D_{\RR_{\Tbf}}(A) \widehat{\otimes}^!_{\Lambda} B) \otimes_{\RR_{\Tbf \times \Tbf}} \RR_{\Tbf}[\dim \Tbf] \\
&= \RGamma(Y, \D_{\RR_{\Tbf}}(A) \otimes_{\RR_{\Tbf}} B)[\dim \Tbf]. 
\end{align*}
\end{proof}

\begin{thm}\label{thmPivotal}
All objects $A \in \HHo$ are left and right dualizable with left and right duals canonically identified with $\D^-(A)$. 
\end{thm}

\begin{proof}
We want to show that there are canonical isomorphisms for all $A,B,C \in \HHo$,
\begin{equation*}
\Hom(A * B, C) = \Hom(A, C * \D^-(B)) = \Hom(B, \D^-(A) * C).
\end{equation*}
By symmetry we will only show the first one. We follow the construction of \cite{BenZviNadlerComplexGroup}. Assume that $A,B,C \in \HHo$. Then we have by Lemma \ref{lemDuality2}
\begin{equation*}
\Hom(A * B, C)  = \RGamma(\Ubf \backslash \Gbf/\Ubf, \D'(A * B) \widehat{\otimes}^!_{\Lambda} C) = \RGamma(\Ubf \backslash \Gbf/\Ubf, \D'(A) * \D'(B) \widehat{\otimes}^! C)
\end{equation*}
and 
\begin{equation*}
\Hom(A, C *\D^{-}(B)) = \RGamma(\Ubf \backslash \Gbf/\Ubf, \D'(A) \widehat{\otimes}^!_{\Lambda} (C *\D^{-}(B)).
\end{equation*}

Replacing $A,B$ by $\D'(A)$ and $\D'(B)$, it is enough to show that 
\begin{equation}\label{eq:Dualizability}
\RGamma(\Ubf \backslash \Gbf/\Ubf, A \widehat{\otimes}^!_{\Lambda}(C * \inv^* B)) = \RGamma(\Ubf \backslash \Gbf/\Ubf, (A * B) \widehat{\otimes}^!_{\Lambda} C).
\end{equation}

Consider the following diagram 
\[\begin{tikzcd}
	{\Ubf \backslash \Gbf \times^{\Ubf} \Gbf/\Ubf} & {\Ubf \backslash \Gbf \times^{\Ubf} \Gbf/\Ubf \times \Ubf \backslash \Gbf/\Ubf} & {\Ubf \backslash \Gbf/\Ubf \times \Ubf \backslash \Gbf/\Ubf \times \Ubf \backslash \Gbf/\Ubf} \\
	{\Ubf \backslash \Gbf/\Ubf} & {\Ubf \backslash \Gbf/\Ubf \times \Ubf \backslash \Gbf/\Ubf}
	\arrow["m"', from=1-1, to=2-1]
	\arrow["\Delta"', from=2-1, to=2-2]
	\arrow["{\Delta_1}", from=1-1, to=1-2]
	\arrow["{\id  \times m}", from=1-2, to=2-2]
	\arrow["q_1", from=1-2, to=1-3]
\end{tikzcd}\]
where $m$ is the multiplication and $\Delta_1 = \id \times m$ and $q_1$ is simply the projection. The square in this diagram is Cartesian. We have
\begin{align*}
(A * B) \widehat{\otimes}^!_{\Lambda} C &= \Delta^!(m \times \id)_!q_1^*(A \widehat{\boxtimes}_{\Lambda} B \widehat{\boxtimes}_{\Lambda} C)) \\
&= m_!\Delta_1^!q_1^*(A \widehat{\boxtimes}_{\Lambda} B \widehat{\boxtimes}_{\Lambda} C).
\end{align*}
This follows from the fact that $m_!(A \boxtimes B) \simeq m_*(A \boxtimes B)[\dim \Tbf]$ as in the proof of Lemma \ref{lemDuality}. 

Similarly, we have  
\[\begin{tikzcd}
	{\Ubf \backslash \Gbf \times^{\Ubf} \Gbf/\Ubf} & {\Ubf \backslash \Gbf \times^{\Ubf} \Gbf/\Ubf \times \Ubf \backslash \Gbf/\Ubf} & {\Ubf \backslash \Gbf/\Ubf \times \Ubf \backslash \Gbf/\Ubf \times \Ubf \backslash \Gbf/\Ubf} \\
	{\Ubf \backslash \Gbf/\Ubf} & {\Ubf \backslash \Gbf/\Ubf \times \Ubf \backslash \Gbf/\Ubf}
	\arrow["m"', from=1-1, to=2-1]
	\arrow["\Delta"', from=2-1, to=2-2]
	\arrow["{\Delta_2}", from=1-1, to=1-2]
	\arrow["{\id  \times m}", from=1-2, to=2-2]
	\arrow["q_2", from=1-2, to=1-3]
\end{tikzcd}\]
where $\Delta_2(x,y) = (m\times \id)$ and $q_2$ is induced by the maps $\Gbf^3 \rightarrow \Gbf^3, (a,b,c) \mapsto (a,\inv(c),b)$ and the square is cartesian. Then we have
\begin{align*}
A \widehat{\otimes}^!_{\Lambda}(C * \inv^*B) &= \Delta^!(\id \times m)_!q_2^*(A \widehat{\boxtimes}_{\Lambda} B \widehat{\boxtimes}_{\Lambda} C) \\
&= m_!\Delta_2^!q_2^*(A \widehat{\boxtimes}_{\Lambda} B \widehat{\boxtimes}_{\Lambda} C)
\end{align*}

Consider now the diagram 
\[\begin{tikzcd}
	{\Ubf \backslash \Gbf \times^{\Ubf} \Gbf/\Ubf} & {\Ubf \backslash \Gbf/\Ubf \times \Ubf \backslash \Gbf/\Ubf \times \Ubf \backslash \Gbf/\Ubf} \\
	{\Ubf \backslash \Gbf/\Ubf} & {\Ubf \backslash \Gbf \times^{\Ubf} \Gbf/\Ubf}
	\arrow["m"', from=1-1, to=2-1]
	\arrow["q_1", from=2-2, to=2-1]
	\arrow["{q_1\Delta_1}"', from=2-2, to=1-2]
	\arrow["{q_2\Delta_2}", from=1-1, to=1-2]
	\arrow["r"{description}, from=1-1, to=2-2]
\end{tikzcd}\]
where $r$ is the map induced by the map $\Gbf \times \Gbf \rightarrow \Gbf \times \Gbf, (x,y) \mapsto (xy, \inv(y))$. This diagram is commutative. We therefore have 
\begin{align*}
 \RGamma(\Ubf \backslash \Gbf/\Ubf, m_!\Delta_2^!q_2^*(A \widehat{\boxtimes}_{\Lambda} B \widehat{\boxtimes}_{\Lambda} C)) &= \RGamma(\Ubf \backslash \Gbf/\Ubf, m_!r^!\Delta_1^!q_1^*(A \widehat{\boxtimes}_{\Lambda} B \widehat{\boxtimes}_{\Lambda} C)) \\
&= \RGamma(\Ubf \backslash \Gbf/\Ubf, q_{1,!}r_!r^!\Delta_1^!q^*_1(A \widehat{\boxtimes}_{\Lambda} B \widehat{\boxtimes}_{\Lambda} C)) \\
&= \RGamma(\Ubf \backslash \Gbf/\Ubf, q_{1,!}\Delta_1^!q_1^*(A \widehat{\boxtimes}_{\Lambda} B \widehat{\boxtimes}_{\Lambda} C)).
\end{align*}
The first line comes from the remark that $q_1$ and $q_2$ are smooth of relative dimension $\dim \Ubf$ hence since $\Delta_2^!q_2^! = r^!\Delta_1^!q_1^!$ after shifting by $[-2\dim \Ubf]$ we get $\Delta_2^!q_2^* = r^!\Delta_1^!q_1^*$. The passage from the third to the fourth line follows from the fact that $r$ is an isomorphism. Now both sides of Equation \ref{eq:Dualizability} are isomorphic to $\RGamma(\Ubf \backslash \Gbf \times^{\Ubf} \Gbf/\Ubf, \Delta_1^!q_1^*(A \widehat{\boxtimes}_{\Lambda} B \widehat{\boxtimes}_{\Lambda} C))$.
\end{proof}

The proof of this theorem yields a finer information. Since all objects of $\HHo$ are dualizable, there are well defined functors $L, R : \HH^{\omega, \op} \rightarrow \HH$ such that for all $x \in \HHo$, the sheaf $L(x)$ (resp. $R(x)$) is the left dual of $x$ (resp. the right dual of $x$). The next corollary is then also a consequence of the proof Theorem \ref{thmPivotal}. 
\begin{corol}\label{corolPivotalStructure}
There is a monoidal isomorphism of functors $\HH^{\omega, \op} \rightarrow \HHo, L \rightarrow R$. 
\end{corol}

We now pass to the $\Ind$-extensions. Recall that we extended the convolution product to all of $\HH$ by continuity. Since the category $\HH$ is compactly generated, it is dualizable and its dual is canonically identified with $\HH^{\vee} = \Ind(\HH^{\omega, \op})$. By extending by continuity the functors $L$ and $R$, we get continuous functor $L,R : \HH^{\vee} \rightarrow \HH$ defined on compact objects by taking left and right duals.

\section{Free monodromic character sheaves}\label{sectionFreeMonocharacterSheaves}

In this section we assume that $G$ is connected, has connected center and that $\ell$ is good for $G$, see \cite[Section 4.3]{SteinbergSpringer}. 

\subsection{Equivariant parabolic induction and restriction functors}

We first recall the equivariant parabolic induction and restrictions. Let $\Lbf\Vbf = \Pbf \subset \Gbf$ be a parabolic subgroup with Levi factor $\Lbf$ and consider the correspondence 
\begin{equation*}
\frac{\Gbf}{\Ad(\Gbf)} \xleftarrow{a} \frac{\Pbf}{\Ad(\Pbf)} \xrightarrow{b} \frac{\Lbf}{\Ad(\Lbf)}.
\end{equation*}

There are well defined functors 
\begin{equation*}
\ind_{\Pbf} : \DDic(\frac{\Lbf}{\Ad(\Lbf)}, \Lambda) \rightarrow \DDic(\frac{\Gbf}{\Ad(\Gbf)}, \Lambda)
\end{equation*}
and 
\begin{equation*}
\res_{\Pbf} : \DDic(\frac{\Gbf}{\Ad(\Gbf)}, \Lambda) \rightarrow \DDic(\frac{\Lbf}{\Ad(\Lbf)}, \Lambda),
\end{equation*}
given by $\ind_{\Pbf} = a_!b^*$ and $\res_P = b_!a^*$. 
\begin{rque}
These functors are linked to the horocycle transform in the following way, let $i : \frac{\Lbf}{\Ad(\Lbf)} \to \frac{\Vbf \backslash \Gbf/\Vbf}{\Ad(\Lbf)}$ be the inclusion of the closed stratum. There are canonical isomorphisms $\ind_{\Pbf} = \chfrak_{\Pbf}i_!$ and $\res_{\Pbf} = i^*\hcfrak_{\Pbf}$. 
\end{rque}

\begin{thm}[\cite{BezYomDin}]\label{thmExactnessParabInd}
Both functors $\ind_{\Pbf}$ and $\res_{\Pbf}$ are perverse $t$-exact. 
\end{thm}

\subsection{The sheaf $\chfrak(\Delta_1)$}

Recall that we denote $\Delta_1 = \oplus_{\chi \in \Ch(\Tbf)}(\Delta_{1, \chi})$. We define a $\Wbf$-action on $\chfrak(\Delta_1)$ and more generally for any $W$-orbit $\mathfrak{o} \subset \chfrak_{\Lambda}(\Tbf)$, we also define a $\Wbf$-action on $\chfrak(\oplus_{\chi \in \mathfrak{o}} \Delta_{1,\chi})$. 

The following construction is due to \cite[Section 3.3]{THChen} and naturally generalizes the Springer action. Let $A \in \Perv(\Tbf, \RR_{\Tbf})^{\Wbf}$ be a $\Wbf$-equivariant sheaf on $\Tbf$. There is a canonical $\Wbf$-action on $\chfrak(A)$ which is constructed as follows. Let $\pi : \tilde{\Gbf} \to \Gbf$ be the Grothendieck-Springer resolution and let $t : \tilde{\Gbf} \to \Tbf$ be the natural projection. We denote by $\Gbf_{rs}, \tilde{\Gbf}_{rs}, \Tbf_{rs}$ the regular semisimple loci. The map $\tilde{\Gbf}_{rs} \to \Gbf_{rs}$ is naturally a $\Wbf$-torsor and the map $t_{rs} : \tilde{\Gbf}_{rs} \to \Gbf_{rs}$ is $\Wbf$-equivariant. We also denote by $j : \Gbf_{rs} \to \Gbf$ the inclusion. As the map $\tilde{\Gbf} \to \Gbf$ is small, for any sheaf $A$ on $\Tbf$ such that $A = j_{!*}A$ we have 
$$\ind_{\Bbf} A = j_{!*}q_{rs,!}t_{rs}^*{A}$$
which is now equipped with a $\Wbf$-action. But since perverse monodromic sheaves on $\Tbf$ are lisse they satisfy $A = j_{!*}A$.

\begin{lem}\label{lemPerversityWInvariants}
The sheaf $\chfrak(\oplus_{\chi}\Delta_{1,\chi})^{\Wbf}$ is perverse. 
\end{lem}

\begin{proof}
We argue by studying the $!$ and $*$-fibers of this sheaf. By \cite[Corollary 6.1.5]{BhattScholze} taking $*$-fibers commutes with small limits, in particular it commutes with the functor $(-)^{\Wbf}$. Since $!$-pullback is a right adjoint it also commutes with limits. Let $s \in \Tbf$ be a semisimple element and denote by $\Zbf^{\circ} = \Zbf_{\Gbf}^{\circ}(s)$ the neutral component of the centralizer of $s$. We denote by $\Spr^{\Zbf^{\circ}}$ the Springer sheaf of $\Zbf^{\circ}$ (define with respect to the Borel $\Bbf \cap \Zbf^{\circ})$. It is a sheaf with an action of $\Wbf_{s}$ the Weyl group of $\Zbf^{\circ}$. The construction of \cite[Corollaire 1.3.2]{LaumonFaisceauxCharactere} shows that, for $u \in \Zbf^{\circ}$ a unipotent element, there is a $\RR_{\Tbf}^{\Wbf}[\Wbf]$-linear isomorphism 
\begin{equation*}
\chfrak(\oplus_{\chi} \Delta_{1, \chi})_{su} = \ind_{\Wbf_s}^{\Wbf} \oplus_{\chi} (\Delta_{1, \chi})_s \otimes \Spr_u^{\Zbf^{\circ}}.
\end{equation*}
After taking $\Wbf$ invariants, we have 
\begin{equation*}
\chfrak(\oplus_{\chi} \Delta_{1, \chi})_{su}^{\Wbf} = (\oplus_{\chi} (\Delta_{1, \chi})_s \otimes \Spr_u^{\Zbf^{\circ}})^{\Wbf_s}.
\end{equation*}
Finally, the hypothesis on $\Gbf$, along with the Pittie-Steinberg theorem \cite{PittieSteinberg} and its completed version \cite[Theorem 8.1]{BezrukRiche} shows that the first factor of this tensor product is a projective $\Wbf_s$-module hence the $\Wbf_s$-invariant sit in the same degrees as the complex without taking $\Wbf_s$-invariant. But, as the parabolic induction functor is $t$-exact by Theorem \ref{thmExactnessParabInd}, the starting complex is perverse hence the $*$-fiber  sits in the expected degrees. The same argument holds for $!$-fibers in place of the $*$-fibers.
\end{proof}

\subsection{The completed category of pre-character sheaves}

As explained in Section \ref{subsectionConvolutionPatterns}, the category $\DDic(\frac{\Gbf}{\Ad(\Gbf)}, \Lambda)$ is equipped with a natural monoidal structure. We want to equip the category $\DDic(\frac{\Gbf}{\Ad(\Gbf)}, \RR_{\Tbf \times \Tbf, \Lambda})$ with a monoidal structure that is comparable with the monoidal structure of the Hecke category. 

We first define a bifunctor 
\begin{align*}
\DDic(\frac{\Gbf}{\Ad(\Gbf)}, \RR_{\Tbf \times \Tbf, \Lambda}) \times \DDic(\frac{\Gbf}{\Ad(\Gbf)}, \RR_{\Tbf \times \Tbf, \Lambda}) &\rightarrow \DDic(\frac{\Gbf}{\Ad(\Gbf)}, \RR_{\Tbf \times \Tbf \times \Tbf \times \Tbf, \Lambda}) \\
(A,B) \mapsto A *^{\Gbf}_{\Lambda} B &= m_!(A \widehat{\boxtimes}_\Lambda B),
\end{align*}
where the map $m$ is the map used in the definition of the convolution in Section \ref{subsectionConvolutionPatterns} and the completed tensor product is the one of Definition \ref{defiConvolCompleted}. 

\begin{defi}
The category of pre-characters sheaves is the full subcategory of sheaves $A \in \DDic(\frac{\Gbf}{\Ad(\Gbf)}, \RR_{\Tbf \times \Tbf, \Lambda})$ such that $\hcfrak(A) \in \DDic(\frac{\Ubf \backslash \Gbf/\Ubf}{\Ad(\Tbf)}, \RR_{\Tbf \times \Tbf})$ lies in the category generated by essential image of the forgetful functor 
$$\oplus_{\chi,\chi' \in \CH_{\Lambda}(T)}(T) \DDic(\frac{\Ubf \backslash \Gbf/\Ubf}{\Ad(\Tbf)}, \RR_{\Tbf \times \Tbf})_{[\chi, \chi']} \to \DDic(\frac{\Ubf \backslash \Gbf/\Ubf}{\Ad(\Tbf)}, \RR_{\Tbf \times \Tbf}).$$ We denote this category by $\Pre\CS^{\wedge}$.
\end{defi}

Similarly there is a bifunctor 
\begin{align*}
\DDic(\frac{\Ubf \backslash \Gbf/\Ubf}{\Ad(\Tbf)}, \RR_{\Tbf \times \Tbf, \Lambda}) \times \DDic(\frac{\Ubf \backslash \Gbf/\Ubf}{\Ad(\Tbf)}, \RR_{\Tbf \times \Tbf, \Lambda}) &\to \DDic(\frac{\Ubf \backslash \Gbf/\Ubf}{\Ad(\Tbf)}, \RR_{\Tbf \times \Tbf \times \Tbf \times \Tbf, \Lambda}) \\
(A,B) \mapsto A *^{\Tbf}_{\Lambda} B = m_!(A \widehat{\boxtimes}_{\Lambda} B)
\end{align*}
where the convolution in a similar way to the convolution on $\frac{\Gbf}{\Ad(\Gbf)}$. 

The next lemma is classical, see for instance \cite[3.3]{BFOCharSheaves}.
\begin{lem}
The functor $\hcfrak$ is monoidal and compatible with both operations $ -*^{\Gbf}_{\Lambda} -$ and $- *^{\Tbf}_{\Lambda} -$. 
\end{lem}

We denote by $\For_{\Tbf^2}^{\Tbf^4} : \DDic(X, \RR_{\Tbf^4, \Lambda})^{\RR_{\Tbf^2}-\cons} \rightarrow \DDic(X, \RR_{\Tbf^2, \Lambda})$ the functor induced by the forgetful functor $\RR_{\Tbf^2} \rightarrow \RR_{\Tbf^4}$ induced by the inclusions of the outer copies of $\Tbf$ in $\Tbf^4$. 

\begin{lem}
The functors $\For_{\Tbf^2}^{\Tbf^4}(- *^{\Tbf}_{\Lambda} -)$ and $\For_{\Tbf^2}^{\Tbf^4}(- *^{\Gbf}_{\Lambda} -)$ define monoidal structures on the full subcategories of $\DDic(\frac{\Ubf \backslash \Gbf/\Ubf}{\Ad(\Tbf)}, \RR_{\Tbf \times \Tbf})$ and $\DDic(\frac{\Gbf}{\Ad(\Gbf}), \RR_{\Tbf \times \Tbf})$ generated by the essential image of 
$$\oplus_{\chi,\chi' \in \CH_{\Lambda}(T)}(T) \DDic(\frac{\Ubf \backslash \Gbf/\Ubf}{\Ad(\Tbf)}, \RR_{\Tbf \times \Tbf})_{[\chi, \chi']}$$
and on $\Pre\CS^{\wedge}$ respectively. 
\end{lem}

\begin{proof}
It is enough to show that for $A,B$ which are $\RR_{\Tbf^2}$-constructible in either of those categories then $A *^{\Tbf}_{\Lambda} B$ is again $\RR_{\Tbf^2}$-constructible. For sheaves on $\DDic(\frac{\Ubf \backslash \Gbf/\Ubf}{\Ad(\Tbf)}, \RR_{\Tbf \times \Tbf})$ this reduces down to the generators of the subcategory and then to the Hecke category where the constructibility assertion is Lemma \ref{lemConstructibleConvol}. For the category of pre character sheaves, we apply the next Lemma \ref{lemHCPreserveConstructibility}.
\end{proof}

\begin{lem}\label{lemHCPreserveConstructibility}
Let $A \in \DDc(\frac{\Gbf}{\Ad(\Gbf)}, \RR_{\Tbf, \Lambda})$ then $A$ is $\Lambda$-construtible if and only is $\hcfrak(A)$ is $\Lambda$-constructible. 
\end{lem}

\begin{proof}
A sheaf $A  \in \DDc(\frac{\Gbf}{\Ad(\Gbf)}, \RR_{\Tbf, \Lambda})$ is $\Lambda$-constructible if and only if it is of $I^{\infty}$-torsion, where $I \subset \RR_{\Tbf, \Lambda}$ is as before the augmention ideal. Being of $I$-torsion is equivalent to asking that $\otimes_{\RR_{\Tbf}} \RR_{\Tbf}[\frac{1}{i}] = 0$ for all $i \in I$. This can be checked after applying the functor $\hcfrak$ by conservativity \ref{corolConservativity}.
\end{proof}

\begin{lem}\label{lemSelfDualityPreCS}
The category $\Pre\CS^{\wedge}$ is compactly generated. Furthermore the functor 
$$A \mapsto \D^{-}(A) = \inv^*\D_{\RR_{\Tbf \times \Tbf}}(A)(\varepsilon_{\Tbf \times \Tbf})$$
defined on constructible sheaves in $\DDic(\frac{\Gbf}{\Ad(\Gbf)}, \RR_{\Tbf \times \Tbf})$, and $(\varepsilon_{\Tbf \times \Tbf})$ is the twisting of the $\RR_{\Tbf \times \Tbf}$-structure introduced in Definition \ref{defiTwist} is a self duality on the category of compact objects in $\Pre\CS^{\wedge}$. 
\end{lem}

\begin{proof}
Since both maps in the correspondence 
$$\frac{\Gbf}{\Ad(\Gbf)} \leftarrow \frac{\Gbf}{\Ad(\Bbf)} \rightarrow \frac{\Ubf \backslash \Gbf/\Ubf}{\Ad(\Tbf)}$$
are representable by finite type schemes, pull-push along these maps preserve compact objects. In particular an object $A \in \Pre\CS^{\wedge}$ is compact if an only if $\hcfrak(A)$ is compact. Conversely if $B \in \DDic(\frac{\Ubf \backslash \Gbf/\Ubf}{\Ad(\Tbf)}, \RR_{\Tbf \times \Tbf})$ is compact, then $\chfrak(B)$ is compact. Hence is it enough to show that $ \DDic(\frac{\Ubf \backslash \Gbf/\Ubf}{\Ad(\Tbf)}, \RR_{\Tbf \times \Tbf})$ is compactly generated. This can be checked on a single Bruhat stratum. Let $w \in \Wbf$, then the stratum corresponding to $w$ is isomorphic to $\frac{\Tbf.w}{\Ad(\Bbf \cap {^w}\Bbf)} \simeq \frac{\Tbf}{\Ad_w(\Tbf)} \times \pt/(\Ubf \cap {^w\Ubf})$. Hence it is enough to check that $\DDic(\frac{\Tbf}{\Ad_w(\Tbf)}, \RR_{\Tbf \times \Tbf})$ is compactly generated which is well known fact.

For the second part of the lemma, first note that as in Definition \ref{defiDualityFunctor} the functor $\D^-$ is the natural duality functor. Hence since $\hcfrak$ is monoidal it commutes with $\D^-$. We only need to check that $\D^-$ preserve the category of compact objects, but this can be seen on $\frac{\Ubf \backslash \Gbf/\Ubf}{\Ad(\Tbf)}$. On this category, this reduces down to Lemma \ref{lemCompactGeneration}. 
\end{proof}

\begin{rque}
The reader should note that even if the category $\Pre\CS^{\wedge}$ is compactly generated, all constructible objects need not be compact. This already happens on a torus, see for instance \cite[Section 7.2]{DrinfeldGaitsgory}. 
\end{rque}

In the next section we will define the category of free monodromic character sheaves as the center of the category $\HH$. The technical part then will be to compute what this category is.

\section{Categorical traces of the free monodromic category}\label{sectionCategoricalTraces}

Recall the following definition of categorical traces, we refer to \cite{BenZviNadlerComplexGroup}, \cite{HoyoisScherotzkeSibilla} and  \cite{ToyModel} for a general presentation of the formalism. 
Let $(\Ccal, \Frob)$ be a pair where $\Ccal \in \Pr_{\Lambda}^L$ is a monoidal category and $\Frob : \Ccal \to \Ccal$ is a monoidal endofunctor. 

\begin{defi}\label{defi:CategoricalTrace}
\begin{enumerate}
\item The categorical trace of $\Frob$ on $\Ccal$ is 
$$\Tr(\Frob, \Ccal) = \Ccal \otimes_{\Ccal \otimes \Ccal^{\rev}} \Ccal_{\Frob}$$
where $\Ccal^{\rev}$ is the monoidal category with underlying category $\Ccal$ and tensor structure given by $x \otimes^{\rev} y = y \otimes x$, and $\Ccal_{\Frob}$ denotes the $\Ccal$-bimodule with right module structure twisted by $\Frob$. 
\item The $\Frob$-categorical center of $\Ccal$ is 
$$\mathcal{Z}_{\Frob}(\Ccal) = \Fun^L_{\Ccal \otimes \Ccal^{\rev}}(\Ccal, \Ccal_{\Frob}).$$
\end{enumerate}
\end{defi}

\subsection{Statements}

We let $\Frob$ be a Frobenius of $\Gbf$. We consider the category $\HH \in \Pr_{\Lambda}^L$ and we equip it with the endomorphism $\Frob_* : \HH \rightarrow \HH$. This morphism is a monoidal equivalence, note that it is a priori not $\RR_{\Tbf \times \Tbf}$-linear. We are interested in computing the $\Frob_*$-trace and $\Frob_*$-center of $\HH$. By Lemma \ref{corolPivotalStructure}, the category $\HH$ is equipped with a pivotal structure and therefore its $\Frob$-trace and $\Frob$-center are canonically identified by \cite[Proposition 3.13]{BenZviNadlerComplexGroup}.

\begin{thm}\label{thmTraceFrobenius}
There is are equivalences $\Zcal_{\Frob_*}(\HH) = \DDic(\pt/\Gbf^{\Frob}, \Lambda) = \Tr(\Frob_*, \HH)$ making the following diagram commutative 
\[\begin{tikzcd}
	& \HH \\
	{\Tr(\Frob_*, \HH)} && {\Zcal_{\Frob_*}(\HH)} \\
	& {\DD(\pt/\Gbf^{\Frob}, \Lambda).}
	\arrow[Rightarrow, no head, from=2-1, to=3-2]
	\arrow[Rightarrow, no head, from=3-2, to=2-3]
	\arrow[from=2-3, to=1-2]
	\arrow[from=1-2, to=2-1]
	\arrow["{\chfrak_{\Frob}}"{description}, from=1-2, to=3-2]
\end{tikzcd}\]
\end{thm}

We also study the categorical trace/center of the identity. 
\begin{defi}
We define $\CS^{\vee}$ the category of free monodromic character sheaves to be 
$$\CS^{\vee} = \Zcal(\HH) = \Tr(\id, \HH).$$
\end{defi}

\begin{thm}\label{thmTraceIdentity}
Assume $G$ is connected, has connected center and $\ell$ is good for $\Gbf$. 
There is an adjunction 
\begin{equation*}
T : \CS^{\vee} \leftrightarrows \Pre\CS^{\wedge} : Z, 
\end{equation*}
such that the map $Z$ is the left adjoint and is monoidal and the right adjoint $T$ is monadic. 

Furthermore, the monad $ZT$ is isomorphic to the convolution along the object $\chfrak(\oplus_{\chi} \Delta_{1,\chi})^{\Wbf}$ which is the unit object in $\CS^{\wedge}$. 
\end{thm}

We will prove Theorem \ref{thmTraceIdentity} in two steps. First we will setup the pair of adjoint functors and show that $ZT$ that is isomorphic to the convolution against some object which has to be unit in $\Zcal(\HH)$. We then compute this object in Section \ref{sectionUnit}. 

\subsection{Categorical Kunneth formulas}

\begin{defi}
A stack $X$ with a $\Tbf$-action is called spherical if there exists a finite $\Tbf$-equivariant stratification of $X$ such that the action of $\Tbf$ on each strata is transitive. 
\end{defi}

\begin{rque}
The link with spherical varieties is the following one. A $\Gbf$-variety $X$ is called spherical if there is an open $\Bbf$-orbit in $X$. Given a spherical variety $X$, the stack $X/\Ubf$ is a spherical stack. The stratification exhibiting it as spherical is the stratification in $\Tbf$-orbits.
\end{rque}

We denote by $\HH_{\Tbf} $ the Hecke category for $\Tbf$, that is, 
$$\HH_{\Tbf}  = \bigoplus_{\chi \in \Ch_{\Lambda}(\Tbf)} \DD(\Tbf, \RR_{\Tbf})_{\chi}.$$

\begin{lem}
Let $X,Y$ be spherical stacks and let $\chi, \chi' \in \Ch(\Tbf)$. The functor $\widehat{\boxtimes}_{\Lambda}$ induces an equivalence 
\begin{equation*}
\DDic(X, \RR_{\Tbf})_{\chi} \otimes_{\DD(\Lambda)} \DDic(Y, \RR_{\Tbf})_{\chi'} \simeq \DDic(X \times Y, \RR_{\Tbf \times \Tbf})_{\chi \times \chi'}.
\end{equation*}
\end{lem}

\begin{proof}
Since both sides are stratified by the $\Tbf$-orbits and there are finitely many strata on both sides, we can reduce to the case where both $X$ and $Y$ have a single stratum. After choosing base points $x$ and $y$ of $X$ and $Y$, by descent along the orbit maps $\Tbf \to X$ and $\Tbf \to Y$ induced by $x$ and $y$, it is enough to treat the case $X = Y = \Tbf$. But now the statement is obvious. 
\end{proof}

\begin{lem}\label{lemKunneth}
Let $X, Y$ be $\Tbf \times \Tbf$-spherical stacks that are also $\Tbf$-spherical after restricting the action along either the first or second inclusion separately.
\begin{enumerate}
\item The stack $X \times^{\Tbf} Y$ is $\Tbf$-spherical with respect to either the left or right action of $\Tbf$, where $\times^{\Tbf}$ denotes the quotient of $X \times Y$ by the action $t.(x,y) = (xt^{-1}, ty)$. 
\item The functor $\widehat{\boxtimes}_{\Lambda}$ induces an equivalence 
\begin{equation*}
\DDic(X, \RR_{\Tbf \times \Tbf})_{[\chi_1, \chi_2]} \otimes_{\DDic(T, \RR_{\Tbf \times \Tbf})_{[\chi_2, \chi_3]}} \DDic(Y, \RR_{\Tbf \times \Tbf})_{[\chi_3, \chi_4]} \simeq \DDic(X \times^{\Tbf} Y, \RR_{\Tbf \times \Tbf})_{[\chi_1, \chi_4]}.
\end{equation*}
This equivalence is $\RR_{\Tbf \times \Tbf}$-linear where the $\RR_{\Tbf \times \Tbf}$-action on the source is given by the outer actions of $\Tbf \times \Tbf$. 
\end{enumerate}
\end{lem}

\begin{proof}
For the point $(i)$, consider the stratification of $X \times Y$ induced by the stratifications of $X$ and $Y$, since all strata are stable under the action of $\Tbf^4$ it is enough to consider the case where $X$ and $Y$ are reduced to a single strata. After choosing base points $x \in X$ and $y \in Y$, we get a surjective map $\Tbf \times \Tbf \to X \times Y$ which is $\Tbf^4$-equivariant and after contracting we get a surjective map $\Tbf = \Tbf \times^{\Tbf} \Tbf \to X \times^{\Tbf} Y$. 

For the point $(ii)$, ad before we can reduce to the case where there is only a single stratum and then to the case where $X = Y = \Tbf$ be the statement is now trivial. 
\end{proof}

\begin{corol}\label{corolKunneth}
The functor 
\begin{align*}
\HH \otimes_{\HH_{\Tbf} } \HH \otimes_{\HH_{\Tbf} } \dots \otimes_{\HH_{\Tbf} } \HH &\rightarrow \bigoplus_{\chi, \chi'} \DDic(\Ubf \backslash \Gbf/\Ubf \times^{\Tbf} \dots \times^{\Tbf} \Ubf \backslash \Gbf/\Ubf)_{[\chi, \chi']} \\
A_1 \otimes A_2 \dots \otimes A_n &\mapsto A_1 \widehat{\boxtimes}_{\Lambda} A_2 \widehat{\boxtimes}_{\Lambda} \dots \widehat{\boxtimes}_{\Lambda} A_n,
\end{align*}
is an equivalence. 
\end{corol}

\begin{proof}
It is clear that the stack $\Ubf \backslash \Gbf/\Ubf$ is spherical. Hence the corollary follows from \ref{lemKunneth}, after taking a direct sum over all $\chi, \chi'$. 
\end{proof}

\subsection{Computation of the trace of Frobenius}\label{subsectionCalculTraceFrobenius}

In this section, we prove Theorem \ref{thmTraceFrobenius}. We follow closely the proof of \cite[Theorem 6.6]{BenZviNadlerComplexGroup}. Recall that we denote by $\HH_{\Tbf} $ the Hecke category for $\Tbf$. Then $\HH$ is an $\HH_{\Tbf} $-bimodule, we write the relative Bar resolution of $\HH$ as an $\HH$-bimodule, that is 
\begin{equation}\label{equationbarResol}
\HH = \varinjlim_{\Delta^{\op}}(\HH \otimes_{\HH_{\Tbf} } \HH \leftleftarrows \HH \otimes_{\HH_{\Tbf} } \HH \otimes_{\HH_{\Tbf} } \HH \dots  )
\end{equation}
where all the maps are the obvious partial multiplications.  We denote by $d_0, d_{n}$ the $n+1$ partial multiplications from $\HH^{\otimes_{\HH_{\Tbf} } n} \to \HH^{\otimes_{\HH_{\Tbf} } n-1}$. This allows us to rewrite the trace and center as 
\begin{equation*}
\Tr(\Frob_*, \HH) = \varinjlim_{[n] \in \Delta^{\op}} \HH^{\otimes_{\HH_{\Tbf} } n+2} \otimes_{\HH \otimes \HH^{\rev}} \HH_{\Frob},
\end{equation*}
and 
\begin{equation*}
\Zcal_{\Frob_*}(\HH) = \varprojlim_{[n] \in \Delta}  \Fun_{\HH \otimes \HH^{\rev}}(\HH^{\otimes_{\HH_{\Tbf} } n+2}, \HH_{\Frob}).
\end{equation*}
Using the adjunction obtained by scalar extension/forgetful functors 
$$(\HH \otimes \HH) \otimes_{\HH_{\Tbf}  \otimes \HH_{\Tbf} } - : (\HH_{\Tbf}  \otimes \HH_{\Tbf}  )-\Mod \leftrightarrows (\HH \otimes \HH)-\Mod : \For,$$
we get 
\begin{enumerate}
\item $\Tr(\Frob_*, \HH)  = \varinjlim_{[n] \in \Delta^{\op}} \HH^{\otimes_{\HH_{\Tbf} } n} \otimes_{\HH_{\Tbf}  \otimes \HH^{\rev}_T} \HH_{\Frob},$
\item $\Zcal_{\Frob_*}(\HH) = \varprojlim_{[n] \in \Delta}  \Fun_{\HH_{\Tbf}  \otimes \HH^{\rev}_T}(\HH^{\otimes_{\HH_{\Tbf} } n}, \HH_{\Frob}).$
\end{enumerate}

The colimit computing the trace is now 
\begin{align*}
\Tr(\Frob_*, \HH) = \varinjlim_{[n] \in \Delta^{\op}} \HH_{\Tbf}  \otimes_{\HH_{\Tbf}  \otimes \HH_{\Tbf} ^{\rev}} \HH_{\Frob} &\leftleftarrows \HH \otimes_{\HH_{\Tbf}  \otimes \HH_{\Tbf} ^{\rev}} \HH_{\Frob} \leftarrow \dots \\ \leftarrow (\HH \otimes_{\HH_{\Tbf} } \HH \otimes_{\HH_{\Tbf} }
\dots \otimes_{\HH_{\Tbf} } &\HH \otimes_{\HH_{\Tbf}  \otimes \HH_{\Tbf} ^{\rev}} \HH_{\Frob}) \dots .
\end{align*}
In the term indexed by $n$, there are in total $n+1$ copies of $\HH$, the last one carrying the index $_{\Frob}$. Using the categorical Künneth formulas of Corollary \ref{corolKunneth}, we have
$$\HH \otimes_{\HH_{\Tbf} } \HH \otimes_{\HH_{\Tbf} } \dots \otimes_{\HH_{\Tbf} } \HH = \oplus_{\chi, \chi'} \DDic(\Ubf \backslash \Gbf/\Ubf \times^{\Tbf} \dots \times^{\Tbf} \Ubf \backslash \Gbf/\Ubf, \RR_{\Tbf \times \Tbf})_{[\chi, \chi']}.$$
Adding the last factor, we get 
\begin{align*}
\HH \otimes_{\HH_{\Tbf} } \HH \otimes_{\HH_{\Tbf} } \dots \otimes_{\HH_{\Tbf} } \HH \otimes_{\HH_{\Tbf}  \otimes \HH_{\Tbf} ^{\rev}} \HH_{\Frob}
&= \oplus_{\chi, \chi'} \DDic(\Ubf \backslash \Gbf/\Ubf \times^{\Tbf} \dots \times^{\Tbf} \Ubf \backslash \Gbf/\Ubf, \RR_{\Tbf \times \Tbf})_{[\chi, \chi']} \\ &\otimes_{\DD(\Tbf \times \Tbf, \RR_{\Tbf \times \Tbf})_{[\chi, \chi']}} \DDic(\Ubf \backslash \Gbf/\Ubf, \RR_{\Tbf \times \Tbf})_{[\chi, \Frob(\chi')]}.
\end{align*}
We now apply the categorical Künneth formula of Lemma \ref{lemKunneth} and we get
$$\HH \otimes_{\HH_{\Tbf} } \HH \otimes_{\HH_{\Tbf} } \dots \otimes_{\HH_{\Tbf} } \HH \otimes_{\HH_{\Tbf}  \otimes \HH_{\Tbf} ^{\rev}} \HH_{\Frob} = \oplus_{\chi, \chi'} \DDic(\frac{\Ubf \backslash \Gbf/\Ubf \times^{\Tbf} \dots \times^{\Tbf} \Ubf \backslash \Gbf/\Ubf}{\Ad_{\Frob}\Tbf}, \RR_{\Tbf \times \Tbf})_{[\chi, \chi']},$$
where there are now $(n+1)$-copies of $G$. 

\begin{lem}\label{lemCalculTrace1}
The forgetful functor induces an equivalence 
$$\oplus_{\chi, \chi'} \DDic(\frac{\Ubf \backslash \Gbf/\Ubf \times^{\Tbf} \dots \times^{\Tbf} \Ubf \backslash \Gbf/\Ubf}{\Ad_{\Frob}\Tbf}, \RR_{\Tbf \times \Tbf})_{[\chi, \chi']} = \DDic(\frac{\Ubf \backslash \Gbf/\Ubf \times^{\Tbf} \dots \times^{\Tbf} \Ubf \backslash \Gbf/\Ubf}{\Ad_{\Frob}\Tbf}, \Lambda).$$
\end{lem} 

\begin{proof}
This is a direct application of Lemma \ref{lemEquivariantvsMonodromic}. 
\end{proof}

Under the equivalence of Lemma \ref{lemCalculTrace1}, from degree $n$ to $n-1$ there are $n+1$-maps 
$$\DDic(\frac{\Ubf \backslash \Gbf/\Ubf^{\times^{\Tbf} n+1}}{\Ad_{\Frob}\Tbf}, \Lambda) \rightarrow \DDic(\frac{\Ubf \backslash \Gbf/\Ubf^{\times^{\Tbf} n}}{\Ad_{\Frob}\Tbf}, \Lambda)$$
denoted by $d_0, \dots, d_n$ which as before are induced by the Bar resolution \ref{equationbarResol}. They are given by the following functors
\begin{enumerate}
\item The first $n$th arrows are induced by partial convolutions diagrams. Namely, for $0 \leq i \leq n-1$ they are given by pull-push along the correspondence 
$$\frac{\Ubf \backslash \Gbf/\Ubf^{\times^{\Tbf} n+1}}{\Ad_{\Frob}\Tbf} \xleftarrow{q_i} \frac{\Ubf \backslash \Gbf/\Ubf \times^{\Tbf} \dots \times^{\Bbf} \dots \times^{\Tbf} \Ubf \backslash \Gbf/\Ubf}{\Ad_{\Frob}\Tbf} \xrightarrow {m_i} \frac{\Ubf \backslash \Gbf/\Ubf^{\times^{\Tbf} n}}{\Ad_{\Frob}\Tbf}$$
given by partial convolution in degree $i+1$. 
\item The last one is induced by the correspondence 
$$\frac{\Ubf \backslash \Gbf/\Ubf^{\times^{\Tbf} n+1}}{\Ad_{\Frob}\Tbf} \xleftarrow{q_0} \frac{\Ubf \backslash \Gbf \times^{\Bbf} \dots \times^{\Tbf} \Ubf \backslash \Gbf/\Ubf}{\Ad_{\Frob}\Tbf} \xrightarrow{m_0} \frac{\Ubf \backslash \Gbf/\Ubf^{\times^{\Tbf} n}}{\Ad_{\Frob}\Tbf}$$
where the first map is induced by $(x_0, \dots, x_n) \mapsto (x_1, \dots, x_n, \Frob(x_0))$ and the second one by partial multiplication on the first two coordinates. 
\end{enumerate}

Consider the simplicial stack obtained as the Cech nerve of the map $\frac{\Gbf}{\Ad_{\Frob}(\Bbf)} \rightarrow \frac{\Gbf}{\Ad_{\Frob}(\Gbf)}$. This is the following (augmented) simplicial stack 
$$\frac{\Gbf}{\Ad_{\Frob}(\Gbf)} \leftarrow \frac{\Gbf}{\Ad_{\Frob}(\Bbf)} \leftleftarrows \frac{\Gbf \times^{\Bbf} \Gbf}{\Ad_{\Frob}(\Bbf)} \leftarrow \dots \leftarrow \frac{\Gbf^{\times^{\Bbf} n}}{\Ad_{\Frob}(\Bbf)} \dots $$
where the transition maps are induced by the partial multiplications. By $!$-descent we have 
$$\DDic(\frac{\Gbf}{\Ad_{\Frob}(\Gbf)}, \Lambda) = \varinjlim_{\Delta^{\op}, !} \DDic(\frac{\Gbf^{\times^{\Bbf} n}}{\Ad_{\Frob}(\Bbf)}, \Lambda)$$ 
and taking right adjoints we get 
$$\DDic(\frac{\Gbf}{\Ad_{\Frob}(\Gbf)}, \Lambda) = \varprojlim_{\Delta, !} \DDic(\frac{\Gbf^{\times^{\Bbf} n + 1}}{\Ad_{\Frob}(\Bbf)}, \Lambda)$$
where the transition maps are given by $!$-pullbacks. 
For all $n$ there is a natural quotient map $p_n : \frac{G^{\times^{\Bbf} n + 1}}{\Ad_{\Frob}(\Bbf)} \to \frac{\Ubf\backslash \Gbf/\Ubf^{\times^{\Tbf} n+1}}{\Ad_{\Frob}}$ which is a $\Ubf^{n+1}$-fibration. Furthermore there is a canonical factorization 
$$\frac{\Gbf^{\times^{\Bbf} n + 1}}{\Ad_{\Frob}(\Bbf)} \xrightarrow{\tilde{q}_i}  \frac{\Ubf \backslash \Gbf/\Ubf \times^{\Tbf} \dots \times^{\Bbf} \dots \times^{\Tbf} \Ubf \backslash \Gbf/\Ubf}{\Ad_{\Frob}\Tbf} \xrightarrow{q_i} \frac{\Ubf\backslash \Gbf/\Ubf^{\times^{\Tbf} n+1}}{\Ad_{\Frob}\Tbf}$$
where $\times^{\Bbf}$ appears in the $i$-position.

\begin{lem}\label{lemCalculTraceFrob2}
Taking $*$-pullback along the maps $p_n$ induces a morphism of simplicial objects 
$$\HH^{\otimes_{\HH_{\Tbf} } n} \otimes_{\HH_{\Tbf}  \otimes \HH_{\Tbf} ^{\rev}} \HH_{\Frob} \rightarrow \DDic(\frac{\Gbf^{\times^{\Bbf} n + 1}}{\Ad_{\Frob}(\Bbf)}, \Lambda).$$
Similarly, consider the two cosimplicial objects $\HH^{\otimes_{\HH_{\Tbf} } n} \otimes_{\HH_{\Tbf}  \otimes \HH_{\Tbf} ^{\rev}} \HH_{\Frob}$ and $\DDic(\frac{G^{\times^{\Bbf} n + 1}}{\Ad_{\Frob}(\Bbf)}, \Lambda)$ obtained from the previous simplicial objects by passing to right adjoints. Then $*$-pullback along the map $p_n$ also defines a morphism of cosimplicial objects. 
\end{lem}

\begin{proof}
For clarity of exposition, we show it in the first degree. 
Consider the diagram 
\[\begin{tikzcd}
	{\frac{\Ubf \backslash \Gbf/\Ubf \times^{\Tbf} \Ubf\backslash \Gbf/\Ubf}{\Ad_{\Frob}\Tbf}} & {\frac{\Ubf \backslash \Gbf \times^{\Bbf} \Gbf/\Ubf}{\Ad_{\Frob}\Tbf}} & {\frac{\Ubf \backslash \Gbf /\Ubf}{\Ad_{\Frob}\Tbf}} \\
	{\frac{\Gbf \times^{\Bbf}\Gbf}{\Ad_{\Frob}(\Bbf)}} && {\frac{\Gbf}{\Ad_{\Frob}(\Bbf)}}
	\arrow["{p_0}"', from=2-3, to=1-3]
	\arrow["{m'}"', from=2-1, to=2-3]
	\arrow["{p_1}", from=2-1, to=1-1]
	\arrow["\pr"', from=1-2, to=1-1]
	\arrow["q_0", from=1-2, to=1-3]
	\arrow["{\tilde{q}_0}"{description}, from=2-1, to=1-2]
\end{tikzcd}\]
where $\pr$ is the quotient map and $q_0$ induced by the multiplication. The right slanted square is Cartesian and the left triangle is commutative. Then we have 
\begin{align*}
p_0^*q_{0,!}\pr^* &= m'_!\tilde{q}_0^*\pr^* \\
&= m'_!p_1^*,
\end{align*}
as desired. 
For the second point, we have 
\begin{align*}
p_1^*\pr_*q_0^! &= \tilde{q}_0^*\pr^*\pr_*q_0^! \\
&=\tilde{q}_0^*q_0^! \\
&=\tilde{q}_0^*q_0^*[2\dim \Gbf/\Bbf] \\
&= m'^*p_0^*[2\dim \Gbf/\Bbf] \\
&= m'^!p_0^*,
\end{align*}
where the second line follows from the fact that $\pr$ is a $\Ubf$-gerbe thus $\pr^*\pr_* = \id$.
\end{proof}

Taking limits yields a functor 
$$\Tr(\Frob_*, \HH) \rightarrow \DDic(\pt/\Gbf^{\Frob}, \Lambda).$$ 
Furthermore the composition $\HH \rightarrow \Tr(\Frob_*, \HH) \rightarrow \DDic(\pt/\Gbf^{\Frob}, \Lambda)$ is nothing else than the functor $\chfrak_{\Frob}$. 

\begin{lem}
The functor $\Tr(\Frob_*, \HH) \rightarrow \DDic(\pt/\Gbf^{\Frob}, \Lambda)$ is an equivalence. 
\end{lem}

\begin{proof}
Since all the arrows $p_n$ are fibrations under unipotent groups the functors $p_n^*$ are fully faithful. Since a limit of fully faithful functors is fully faithful, the resulting functor $\Tr(\Frob_*, \HH) \rightarrow \DDic(\pt/\Gbf^{\Frob}, \Lambda)$ is fully faithful. 
Since the composition $\HH \rightarrow \Tr(\Frob_*, \HH) \rightarrow \DDic(\pt/\Gbf^{\Frob}, \Lambda)$ is the functor $\hcfrak_{\Frob}$, the essential image of this functor contains all the cohomology complexes of the Deligne-Lusztig varieties. By Corollary \ref{corolGenerationDL},  these generate the target category which implies that this functor is essentially surjective.  
\end{proof}

This completes the proof for the equivalence between the trace and the category $\DDic(\pt/\Gbf^{\Frob}, \Lambda)$. The statement about the center follow from the general nonsense of \cite[Proposition 3.13]{BenZviNadlerComplexGroup}, the quasi-rigidity of the category $\HH$, see Theorem \ref{thmQuasiRigidity} and its  pivotal structure, see Corollary \ref{corolPivotalStructure}. 

\subsection{Computation of the trace of the identity}\label{sectionTraceIdentity}

In this section, we prove the Theorem \ref{thmTraceIdentity}. We use the same construction as in Section \ref{subsectionCalculTraceFrobenius}. Everything holds up until Lemma \ref{lemCalculTrace1}, hence we have 
$$\Tr(\id, \HH) = \varinjlim_{[n] \in \Delta^{\op}} \bigoplus_{\chi,\chi'} \DDic(\frac{\Ubf \backslash \Gbf/\Ubf^{\times^{\Tbf} n+1}}{\Ad(\Tbf)},\RR_{\Tbf \times \Tbf})_{[\chi, \chi']}.$$ The analog of Lemma \ref{lemCalculTraceFrob2} holds if $\Frob = \id$.

Let $n \geq 0$ be an integer, we consider the functor 
\begin{align*}
T_n : \bigoplus_{\chi,\chi'} \DDic(\frac{\Ubf \backslash \Gbf/\Ubf^{\times^{\Tbf} n+1}}{\Ad(\Tbf)},\RR_{\Tbf \times \Tbf})_{[\chi, \chi']} &\rightarrow \DDic(\frac{\Ubf \backslash \Gbf/\Ubf^{\times^{\Tbf} n+1}}{\Ad(\Tbf)},\RR_{\Tbf \times \Tbf}) \\
&\xrightarrow{p_n^*} \DDic(\frac{G^{\times^{\Bbf} n+1}}{\Ad(\Bbf)}, \RR_{\Tbf \times \Tbf}),
\end{align*}
where the first map is the forgetful functor. By the same argument of Lemma \ref{lemCalculTraceFrob2}, the maps $T_n$ define a map cosimplicial objects.

Passing to the limit, we get a functor $\Tr(\id, \HH) \to \DD(\frac{\Gbf}{\Ad(\Gbf)}, \RR_{\Tbf \times \Tbf})$. The composition 
$$\HH \to \Tr(\id, \HH) \to \DD(\frac{\Gbf}{\Ad(\Gbf)}, \RR_{\Tbf \times \Tbf})$$
clearly factors through $\Pre\CS^{\wedge}$. Hence we get functor 
\begin{equation*}
T : \Tr(\id, \HH) \to \Pre\CS^{\wedge}.
\end{equation*}
Since all functors $T_n$ have left adjoints so does  the functor $T$. Let us denote by $Z$ the left adjoint to these functors. 

\begin{lem}
The functor $Z : \Pre\CS^{\wedge} \to \Tr(\id, \HH) = \Zcal(\HH)$ is monoidal. 
\end{lem}

\begin{proof}
From the compatibility of the functor $T$ with $\chfrak$ already established when computing the trace of Frobenius, passing to left adjoints shows that the functor $\Pre\CS^{\wedge} \to \Zcal(\HH) \to \HH$ is the functor $\hcfrak$ which is monoidal. It is also central, see for instance by \cite{BenZviNadlerComplexGroup}. Hence we get that the factorization to $\Zcal(\HH)$ is monoidal.
\end{proof}

\begin{proof}[Proof of \ref{thmTraceIdentity}, part 1]
We now prove the existence part in the statement of the theorem. By Lemma \ref{lemSelfDualityPreCS}, the category $\Pre\CS^{\wedge}$ is compactly generated and is thus dualizable as presentable category. It is even self dual and a self duality is given by the functor $\D^-$. Hence the endofunctor $TZ$ can be represented as $TZ = A * -$ for some object $A \in \Pre\CS^{\wedge}$. 
\end{proof}

\begin{rque}
The object $A \in \Pre\CS^{\wedge}$ has a canonical algebra structure coming from the monad structure on $TZ$. 
\end{rque}

\begin{lem}
Consider the category $\CS^{\wedge, \Lambda-c}$ of characters sheaves that are constructible as $\Lambda$-modules. The forgetful functor $\CS^{\wedge, \Lambda-c} \to \DD(\frac{\Gbf}{\Ad(\Gbf)}, \Lambda)$ is fully faithful. 
\end{lem}

\begin{proof}
The forgetul functor has a natural left adjoint given by $B \mapsto A*(\oplus_{\chi} \Av_{\chi} (\RR_{\Tbf \times \Tbf} \otimes_{\Lambda}  B))$. Hence it is enough to show that the adjunction map $B \to A*(\oplus_{\chi} \Av_{\chi} (\RR_{\Tbf \times \Tbf} \otimes_{\Lambda}  B))$ is an isomorphism. This can be checked after applying the functor $\hcfrak$. But now we have $\hcfrak(A*(\oplus_{\chi} \Av_{\chi} (\RR_{\Tbf \times \Tbf} \otimes_{\Lambda}  B))) = \oplus_{\chi} \Delta_{1,\chi} * B = B$. 
\end{proof}

\begin{rque}
It follows that our category $\CS^{\wedge, \Lambda-c}$ then coincides with the pro-completion of the category spanned by character sheaves as done in \cite{BezrukavnikovTomalchov}. In particular our construction recovers the usual notion of character sheaves of Lusztig. 
\end{rque}

\subsection{The free monodromic unit in character sheaves}\label{sectionUnit}

From now on we assume the hypothesis of Lemma \ref{lemPerversityWInvariants}, that is $\Gbf$ is connected, has connected center and that $\ell$ is not a torsion prime for $\Gbf$.
Theorem \ref{thmTraceIdentity} equips the category $\CS^{\wedge}$ with a canonical monoidal structure. We now wish to identify the unit in this category. We will recover a construction of \cite{BezrukavnikovTomalchov}, we denote by $\delta_1^{\CS^{\wedge}}$ the unit in $\Zcal(\HH)$, this is the object $A$ of the previous section. 

\begin{thm}\label{thmUnitCenter}
There is a canonical isomorphism $\delta_{1}^{\CS^{\wedge}}$ between the unit in character sheaves and $\chfrak(\oplus_{\chi} \Delta_{1, \chi})^{\Wbf}$ in $\DDic(\frac{\Gbf}{\Ad(\Gbf)}, \RR_{\Tbf\times \Tbf}^{\Wbf \times \Wbf})$. Moreover both objects are perverse sheaves on $\frac{\Gbf}{\Ad(\Gbf)}$. 
\end{thm}

\begin{proof}
Firstly since the object $\delta_{1}^{\CS^{\wedge}}$ is the unit in the categorical center of $\HH$, its image in $\HH$ is the unit of $\HH$ hence $\hcfrak(\delta_{1}^{\CS^{\wedge}}) = \oplus_{\chi} \Delta_{1, \chi}$, by adjunction this defines a map 
$$\delta_{1}^{\CS^{\wedge}} \to \chfrak(\oplus_{\chi} \Delta_{1, \chi}).$$
We want to see that this map induces an isomorphism 
$$ \delta_{1}^{\CS^{\wedge}} = \chfrak(\oplus_{\chi} \Delta_{1, \chi})^{\Wbf}.$$
From Theorem \ref{lemMirkovicVilonen}, there is an isomorphism 
\begin{equation}\label{equationSpringerUnitCS}
\chfrak(\oplus_{\chi} \Delta_{1, \chi}) = \Spr * \delta_{1}^{\CS^{\wedge}}. 
\end{equation}
We will prove in Lemma \ref{lemEquivarianceUnderW} that the map (\ref{equationSpringerUnitCS}) is $\Wbf$-equivariant. 

Since $\hcfrak(\delta_{1}^{\CS^{\wedge}}) = \res_{\Tbf}(\delta_{1}^{\CS^{\wedge}})$ is isomorphic to the parabolic restriction of $\Tbf$, we get that $\delta_{1}^{\CS^{\wedge}}$ is perverse as $\res_{\Tbf}$ is $t$-exact and $\hcfrak$ is conservative. We also know that $\Spr * \delta_{1}^{\CS^{\wedge}}$ and its $\Wbf$-invariants are perverse by Lemma \ref{lemPerversityWInvariants}, hence we can compute the $\Wbf$-invariants in the abelian category of perverse sheaves. Hence we have
$$ (\Spr * \delta_{1}^{\CS^{\wedge}})^{\Wbf} = (^pH^0(\Spr^{\Wbf})) * \delta_{1}^{\CS^{\wedge}} = \delta_{1}^{\CS^{\wedge}}.$$
This proves the claim.
\end{proof}

\begin{lem}\label{lemEquivarianceUnderW}
The map (\ref{equationSpringerUnitCS}) is $\Wbf$-equivariant. 
\end{lem}

\begin{proof}
The proof of this lemma is done as in \cite[Proposition 5.4.6]{BezrukavnikovTomalchov}. We highlight the main steps.
\begin{enumerate}
\item Both the source and the target of map (\ref{equationSpringerUnitCS}) are perverse hence the $\Wbf$-equivariance is a property and not a structure. It can be checked after applying $\res_{\Tbf}$ as both objects lie in the category generated by its adjoint. 
\item After applying $\res_{\Tbf}$, the source of the map can be filtered by the Mackey filtration. As in \cite{THChen}, this filtration can be split in the category of perverse sheaves compatibly with the $\Wbf$-actions. Note that this still holds in the modular settings. This is a statement that is true only in the category of perverse sheaves (not in the derived category) and holds in the category of perverse sheaves for dimensional reasons. 
\end{enumerate}
The rest of the proof goes as in \cite{BezrukavnikovTomalchov}.
\end{proof}

\appendix

\section{Recollections on $6$-functors formalism and twisted equivariant sheaves}

\subsection{$6$-functors formalisms}

In this appendix, we recall the definition of a $6$-functors formalism as in \cite{SixFunctors} and \cite[Appendix]{MannThesis}. Let $\Ccal$ be a category with all finite limits and $E$ be a class of maps stable under composition, pullbacks and containing all isomorphisms. There is an $(\infty, 1)$-symmetric monoidal category $\Corr(\Ccal, E)$ of correspondences. The objects of $\Corr(\Ccal, E)$ are the same as the objects of $\Ccal$. The morphisms $X \rightarrow Y$ are given by the diagrams 
$$X \leftarrow C \rightarrow Y$$ where the maps $C \rightarrow Y$ belongs to $E$. The symmetric monoidal structure is given by Cartesian product in $\Ccal$. 

\begin{defi}
A $3$-functor formalism is a lax-symmetric monoidal functor
\begin{equation*}
\DD : \Corr(\Ccal, E) \rightarrow 1-\Cat,
\end{equation*}
where $1-\Cat$ denotes the $(\infty,1)$-category of small $(\infty, 1)$-categories. 
\end{defi}

Let $\DD$ be a $3$-functors formalism. Let $f : X \rightarrow Y$. 
\begin{enumerate}
\item the correspondence $Y \xleftarrow{f} X = X$ is sent to the functor $f^* : \DD(Y) \rightarrow \DD(X)$, 
\item the correspondence $X = X \xrightarrow{f} Y$ is sent to the functor $f_! : \DD(X) \rightarrow \DD(Y)$, 
\item the lax-symmetric monoidality of the functor $\DD$ equips $\DD(X)$ with a symmetric monoidal structure given by 
\begin{equation*}
\DD(X) \times \DD(X) \rightarrow \DD(X \times X) \xrightarrow{\Delta^*} \DD(X), 
\end{equation*}
where the first map is the lax-symmetric monoidal constraint and the second map is the pullback under the diagonal $\Delta : X \rightarrow X \times X$. 
\end{enumerate}

\begin{defi}
A $6$-functors formalism is the data of a $3$-functors formalism such that for all $f$ and all $A \in \DD(X)$, the functors $f^*, f_!$ and $A \otimes -$ admit right adjoints.
\end{defi}

Given a $6$-functor formalism, all adjunctions maps, Künneth maps, projection formulas and base change formulas are compatible in a very strong sense. We refer to \cite{SixFunctors} for a discussion. 

\begin{rque}
We will often consider the following variant of six functors formalism. Namely, we consider a lax-monoidal functor 
$$\DD : \Corr(\Ccal, E) \to \mathrm{Pr}^L_{\Lambda}$$
the category of presentable $\Lambda$-linear categories. In this context the existence of right adjoints is automatic by the adjoint functor theorem \cite[Proposition 5.5.2.2]{LurieA}. 
\end{rque}

\begin{defi}\label{defiDescentFor6Functors}
Let $\DD$ be a $6$-functors formalism on $(\Ccal, E)$. And let $f : X \rightarrow Y$ be a morphism in $E$, 
\begin{enumerate}
\item the morphism $f$ is of universal $*$-descent if after any pullback we have 
\begin{equation*}
\DD(Y) = \varprojlim_{\Delta,*} \DD(X^{\times_Y n}),
\end{equation*}
where $X^{\times_Y n}$ denotes the $n$-fold fiber product $X \times_Y X \times_Y \dots \times_Y X$ and the transition maps are given by pullback. 
\item The morphism $f$ is of universal $!$-descent if after any pullback, we have 
\begin{equation*}
\DD(Y) = \varprojlim_{\Delta,!} \DD(X^{\times_Y n}),
\end{equation*}
where the transition maps are given by $!$-pullback. 
\end{enumerate}
\end{defi}

\begin{thm}[\protect{\cite[Appendix A.5.16]{MannThesis}, \cite[Proposition 4.17]{SixFunctors}}] 
Let $\DD$ be a $6$-functors formalism on $(\Ccal, E)$ and assume that $\Ccal$ is equipped with a subcanonical Grothendieck topology $\tau$ and that all $\tau$-covers are of universal $!$ and $*$-descent. Then there exists a unique extension of $\DD$ to $(\Stk_{\tau}(\Ccal), \tilde{E})$ where $\Stk_{\tau}(\Ccal)$ denotes the category of stacks for the $\tau$-topology and $\tilde{E}$ is the class of maps such that after any pullback to an object in $\Ccal$, it belongs in $E$.  
\end{thm}

\subsection{Twisted equivariant sheaves}\label{subsectionTwistedEquivSheaves}

We let $\DD$ be a $6$-functors formalism on $(\Ccal, E)$. Let $\Lambda$ be a base ring such that $\DD(\pt) = \DD(\Lambda)$, we assume that for all $X \in \Ccal$ the category $\DD(X)$ lies in $\Pr_{\DD(\Lambda)}^L$. Let $G$ be a group object in $\Ccal$ such that the multiplication and unit maps of $G$ are in $E$. Consider the functor 

\begin{align*}
\Ccal &\rightarrow \Corr(\Ccal, E) \\
X &\mapsto X \\
(f : X \rightarrow Y) &\mapsto (X = X \rightarrow Y).
\end{align*}
This functor is monoidal with respect to the Cartesian product on $\Ccal$, hence the group object $G$ can also be seen as a group object in the category of correspondences. The lax monoidality of the functor $\DD$ then equips the category $\DD(G)$ with a monoidal structure, usually called the $!$-convolution. By definition the monoidal structure is 
\begin{align*}
\DD(G) \times \DD(G) &\rightarrow \DD(G) \\
(A,B) &\mapsto A*B = m_!(A \boxtimes B), 
\end{align*}
where $m$ is the multiplication map of $G$.

We denote by $\Pr_{\DD(G)}^L$ the category of modules over $\DD(G)$ in $\Pr^L$. Similarly, any object $X$ in $\Ccal$ equipped with an action of $G$ such that the action map $G \times X \rightarrow X$ lies in $E$ yields an object in $\Pr_{\DD(G)}^L$. That is $\DD(X)$ is naturally a module over $\DD(G)$. We denote by $\DD(\Lambda)_{\triv}$ the corresponding $\DD(G)$-module when $X = \pt$.

Given a category $C \in \Pr_{\DD(G)}^L$, its categorical invariants and coinvariants are defined as follows 
\begin{enumerate}
\item the category of invariants $C^G := \Fun^L_{\DD(G)}(\DD(\Lambda)_{\triv}, C)$, 
\item the category of coinvariants $C_G := \DD(\Lambda)_{\triv} \otimes_{\DD(G)} C$. 
\end{enumerate}

We take the six functor formalism of étale sheaves of $\Lambda$-modules on schemes of finite type over $k$ as defined in Section \ref{subsectionCategoriesOfScheaves}. Let $G$ be a finite dimensional group scheme. 

\begin{thm}[\protect{\cite[Theorem B.1.2]{GaiWhittCats}}]
Let $C$ be a $G$-category. The two categories $C^G$ and $C_G$ are canonically identified moreover the natural forgetful map $\For : C^G \to C$ has a left adjoint denoted by $\Av_! : C \to C^G$ such that the composition
$$\For \Av_!(c) = a(\omega_G \otimes c)$$
where $a : \DD(G) \otimes C \to C$ is the action map and $\omega_G$ is the dualizing sheaf of $G$. 
\end{thm}

\begin{rque}
The proof of \cite{GaiWhittCats} is done when $\Lambda = \Qlb$, however the argument relies on the fact that $G$ is quasi-compact and that $\RGamma(G, \Lambda) \in \DD(\Lambda)$ is a compact object. Both conditions are satisfied in our situation since all our coefficient rings $\Lambda$ are (filtered colimits of) regular rings and the group $G$ has finite dimensional cohomology. 
\end{rque} 

\begin{rque}
In \cite{GaiWhittCats}, the averaging functor is right adjoint to the forgetful functor when ours is a left adjoint. This difference comes from the fact that we use $!$-convolution instead of $*$-convolution. 
\end{rque}

\begin{rque}
The natural map $\delta_1 = 1_!1^!\omega_G \to \omega_G$ induces after convolving with $\omega_G$ a natural coalgebra structure on $\omega_G$. For a $G$-category, the category of $G$-invariants is then identified with the category of comodules over this coalgebra. 
\end{rque} 

We now want to discuss multiplicative sheaves on $G$. Let $\psi$ denote an action of $\DD(G)$ on $\DD(\Lambda)$ and denote by $a_{\psi} : \DD(G) \to \DD(\Lambda)$ the action map. Since all categories in sight are dualizable, we can pass to duals and using Verdier duality, we identify $\DD(G)$ and $\DD(\Lambda)$ with their own duals. We denote by $\Lcal \in \DD(G) $ the object such that $a_{\psi}^{\vee}(\Lambda) = \D(\Lcal)$. Consider the associativity axiom for the action of $\DD(G)$ on $\DD(\Lambda)$, namely the commutativity of the following diagram 
\[\begin{tikzcd}
	{\DD(G) \otimes \DD(G)} & {\DD(G)} \\
	{\DD(G)} & {\DD(\Lambda)}
	\arrow["{*}", from=1-1, to=1-2]
	\arrow["{\id \otimes a_{\psi}}"', from=1-1, to=2-1]
	\arrow["{a_{\psi}}", from=1-2, to=2-2]
	\arrow["{a_{\psi}}"', from=2-1, to=2-2]
\end{tikzcd}\]
where $*$ denotes is the convolution. Passing to duals and taking the image of $\Lambda \in \DD(\Lambda)$ yields an isomorphism 
\begin{equation}\label{eqMultiplicative}
m^*\Lcal = \Lcal \boxtimes \Lcal
\end{equation}
together with a trivialization of of $\Lcal$ at $1$ and higher associativity constraints coming from the higher associativity constraints in the diagram. We note that this forces $\Lcal$ to be lisse, concentrated in degree $0$ and locally free of rank one on $G$. 
\begin{lem}
Assume that $G$ is connected and normal, then the isomorphism of Equation \ref{eqMultiplicative} is canonically determined by the data of $\Lcal$ and its trivialization at $1$. 
\end{lem}
\begin{proof}
Since $\Lcal$ is lisse, the data of $\Lcal$ equipped with its trivialization is equivalent to the data of a morphism $\chi : \pi_1(G,1) \to \Lambda^{\times}$. The data of the isomorphism of of Equation \ref{eqMultiplicative} comes from the commutativity of the diagram 
\[\begin{tikzcd}
	{\pi_1(G,1)} & {\Lambda^{\times}} \\
	{\pi_1(G^2,1)} & {\pi_1(G,1) \times \pi_1(G,1)}
	\arrow["m", from=2-1, to=1-1]
	\arrow[from=2-1, to=2-2]
	\arrow["{\chi \otimes \chi}"', from=2-2, to=1-2]
	\arrow["\chi", from=1-1, to=1-2]
\end{tikzcd}\]
\end{proof}
Hence from an action of $\DD(G)$ on $\DD(\Lambda)$ we have produced a rank one locally constant sheaf $\Lcal$ equipped with a canonical isomorphism $m^*\Lcal = \Lcal \boxtimes \Lcal$. Conversely assume we are given such a sheaf $\Lcal$. Then the projection formula yields a canonical isomorphism 
$$(A * B) \otimes \Lcal = (A \otimes \Lcal) * (B \otimes \Lcal)$$
for all $A,B \in \DD(G)$. In particular the functor 
\begin{align*}
\DD(G) &\to \DD(\Lambda) \\
A &\mapsto \RGamma_c(G, A \otimes \Lcal)
\end{align*}
extends to an action of $\DD(G)$ on $\DD(\Lambda)$. 
\begin{defi}\label{defiMultiplicativeSheaf}
Let $G$ be a connected normal group, a multiplicative sheaf on $G$ is one of the following equivalent data
\begin{enumerate}
\item A morphism $\pi_1(G,1) \to \Lambda^{\times}$ such that the diagram 
\[\begin{tikzcd}
	{\pi_1(G,1)} & {\Lambda^{\times}} \\
	{\pi_1(G^2,1)} & {\pi_1(G,1) \times \pi_1(G,1)}
	\arrow["m", from=2-1, to=1-1]
	\arrow[from=2-1, to=2-2]
	\arrow["{\chi \otimes \chi}"', from=2-2, to=1-2]
	\arrow["\chi", from=1-1, to=1-2]
\end{tikzcd}\]
commutes.
\item A rank one locally constant sheaf $\Lcal$ on $G$, equipped with a trivialization $1^*\Lcal = \Lambda$ and an isomorphism $m^*\Lcal = \Lcal \boxtimes \Lcal$ compatible with this trivialization. 
\item An action of $\DD(G)$ on $\DD(\Lambda)$. 
\end{enumerate}
\end{defi}

\begin{defi}[\cite{GaiWhittCats}]\label{defiTwistedEquivSheaves}
Let $X$ be a stack with an action of $G$ and let $\Lcal$ be a multiplicative sheaf. We twist the natural action of $\DD(G)$ on $\DD(\Lambda)$ by $\Lcal$, so that the action of given by the natural action on $\DD(X) \otimes_{\DD(\Lambda)} \DD(\Lambda)_{\Lcal^{-1}}$. The category of $\Lcal$-equivariant sheaves on $X$ is defined to be the category of $G$-invariants for this action. 
\end{defi}

\begin{rque}\label{rqueEquivAsModules}
The natural map $\delta_1 = 1_!1^!(\Lcal \otimes \omega_G) \to \Lcal \otimes \omega_G$ equips $\Lcal \otimes \omega_G$ with a natural algebra structure. The category of $(G, \Lcal)$-equivariant sheaves is then the category of modules for this algebra. 
\end{rque}

\bibliographystyle{alpha}
\bibliography{biblio}
\end{document}


In this section we want to recall known facts about the formalism of categorical traces and categorical centers. We fix $\Lambda$ a base ring. Let $(\Ccal, \star) \in \Pr_{\Lambda}$ be a monoidal category. We denote by $\otimes$ Lurie's tensor product in $\Pr_{\Lambda}^L$.  We denote by $\Ccal^{\rev}$ the category $\Ccal$ with its monoidal structure reverse. Let $\Mcal \in \Pr_{\Lambda}$ be a $\Ccal$-bimodule. Its Hochschild homoloogy and cohomology are defined as 

\begin{enumerate}
\item $HH_*(\Ccal, \Mcal) = \Ccal \otimes_{\Ccal \otimes \Ccal^{\rev}} \Mcal \in \Pr_{\Lambda}$. 
\item $HH^*(\Ccal, \Mcal) = \Fun^L_{\Ccal \otimes \Ccal^{\rev}}(\Ccal, \Mcal) \in \Pr_{\Lambda}$.
\end{enumerate}
where $\Ccal$ is equipped with its usual $\Ccal$-bimodule structure. These are the usual categorical generalizations of Hochschild (co)homology. Let $\Frob : \Ccal \to \Ccal$ be a monoidal endomorphism of $\Ccal$. We denote by $\Ccal_{\Frob}$ the $\Ccal$-bimodule where the right action of $\Ccal$ is twisted by $\Frob$, that is for $x \in \Ccal_{\Frob}$ and $c \in \Ccal$ the action is given by 
$$x \cdot c = x \star \Frob(c)$$. We will denote by 
\begin{enumerate}
\item $\Tr(\Frob, \Ccal) = \HH_*(\Ccal, \Ccal_{\Frob})$ and we call it the categorical trace of $\Frob$ on $\Ccal$, 
\item $\Zcal_{\Frob}(\Ccal) = \HH^*(\Ccal, \Ccal_{\Frob})$ and we call it the $\Frob$-center of $\Frob$ on $\Ccal$. 
\end{enumerate}

From now on we assume that $\Ccal$ is compactly generated and denote by $\Ccal^{\omega}$ the subcategory of compact objects. 
\begin{defi}
The category $\Ccal$ is quasi-rigid if every compact object is left and right dualizable. It is rigid if the unit of $\Ccal$ is compact. 
\end{defi}

Assume that $\Ccal$ is quasi-rigid. We denote by $L,R : \Ccal^{\omega, \op} \to \Ccal^{\omega}$ the functors taking left and right duals. After passing to ind-categories, we obtain functors $\Ccal^{\wedge} \to \Ccal$ where $\Ccal^{\wedge}$ is the category dual to $\Ccal$. 

\begin{thm}
\todo{cite BenZvi Nadler}
Assume that $\Ccal$ is quasi-rigid, then for all bimodules $\Mcal$ there is a canonical equivalence 
$$HH^*(\Ccal, \Mcal) = HH_*(\Ccal, \Mcal_{LL})$$
where $\Mcal_{LL}$ is the $\Ccal$-bimodule $\Mcal$ with its right module structure twisted by $LL$. 
\end{thm} 

\begin{defi}
The data of a monoidal isomorphism of functors $\id \to LL$ is called a pivotal structure. 
\end{defi}

\begin{corol}
Assume that $\Ccal$ is equipped with a pivotal structure, then there is an isomorphism (induced by the pivotal structure) 
$$HH^*(\Ccal, \Mcal) = HH_*(\Ccal, \Mcal)$$
for all $\Ccal$-bimodule $\Mcal$. 
\end{corol}